\documentclass[bj,preprint]{imsart2}

\RequirePackage[OT1]{fontenc}
\RequirePackage{amsthm,amsmath}
\RequirePackage[numbers]{natbib}
\RequirePackage[colorlinks,citecolor=blue,urlcolor=blue]{hyperref}

\newcommand{\ip}{\mathop{\longrightarrow}\limits^{\footnotesize\mbP^\star}}

\newcommand{\norm}[1]{\lVert#1\rVert}
\arxiv{1702.07801}

\startlocaldefs
\numberwithin{equation}{section}
\theoremstyle{plain}

\endlocaldefs

\newcommand{\size}{L}
\newcommand{\bzs}{\bz^\star}
\newcommand{\bthetas}{\btheta^\star}
\newcommand{\bthetah}{\widehat\btheta}

\usepackage[mathscr]{eucal}
\usepackage{afterpage}
\usepackage{amsmath}
\usepackage{tikz}

\usepackage{centernot}
\usepackage{dlfltxbcodetips,bbm,mathtools,subfigure,amsfonts,amsmath,amssymb,fancybox,graphicx,bm,latexsym}

\usepackage{tikz}

\newcommand{\zs}{\bz^\star}

\newcommand{\bthetad}{\dot\btheta}
\newcommand{\truth}{(\bthetas, \zs)}
\newcommand{\zh}{\widehat{\bz}}

\newcommand{\estimate}{(\bthetah, \zh)}

\usetikzlibrary{fit,positioning}

\newcommand{\lip}{{\tiny\mbox{Lip}}}

\newcommand{\wh}{\widehat}

\newcommand{\uu}{a_1}
\newcommand{\vv}{a_2}
\newcommand{\logit}{\mbox{logit}}
\newcommand{\tv}{{\tiny\mbox{TV}}}

\usepackage[tikz]{bclogo}

\newcommand{\su}{\sup\limits}
\newcommand{\rint}{\mathop{\mbox{rint}}}

\newcommand{\bt}{\begin{bclogo}[couleur={rgb:orange,0;yellow,0;white,1},arrondi=0.1,logo=\bcplume,ombre=true]}
\newcommand{\et}{\end{bclogo}\s}
\newcommand{\btt}{\begin{box}}
\newcommand{\ett}{\end{box}}

\newcommand{\btheorem}{\begin{bclogo}[couleur={rgb:orange,0;yellow,0;white,1},arrondi=0.1,logo=\bcplume,ombre=true]{Theorem}}
\newcommand{\ettheorem}{\end{bclogo}}

\newcommand{\bst}{\begin{bclogo}[couleur={rgb:orange,1;yellow,1;white,0.5},arrondi=0.1,logo=\bcpanchant]}
\newcommand{\est}{\end{bclogo}}

\DeclareMathOperator*{\argmin}{arg\,min}
\DeclareMathOperator*{\argmax}{arg\,max}
\DeclareMathOperator*{\ri}{\mbox{ri}}

\newcommand{\benum}{\begin{enumerate}}
\newcommand{\eenum}{\end{enumerate}}

\newcommand{\tends}{\to}

\newcommand{\bq}{\begin{quote}\em}
\newcommand{\eq}{\end{quote}}
\newcommand{\bbq}{\begin{quote}\bf\em}
\newcommand{\ebq}{\end{quote}}

\renewcommand{\=}{&=&}

\newcommand{\lte}{&\leq&}
\newcommand{\gte}{&\geq&}
\newcommand{\mR}{\mathbb{R}}

\newcommand{\mbR}{\mathbb{R}}

\newcommand{\mbX}{\mathbb{X}}
\newcommand{\mbA}{\mathbb{A}}

\newcommand{\mbL}{\mathbb{L}}
\newcommand{\one}{\mathbbm{1}}

\newcommand{\mbE}{\mathbb{E}}

\newcommand{\mbS}{\mathbb{S}}
\newcommand{\mbP}{\mathbb{P}}

\newcommand{\hide}[1]{}
\newcommand{\ghost}[1]{}
\newcommand{\ba}{\begin{array}{llllllllll}}
\newcommand{\ea}{\end{array}}
\newcommand{\bea}{\begin{equation}\begin{array}{llllllllll}}
\newcommand{\eea}{\end{array}\end{equation}}
\newcommand{\be}{\begin{equation}\begin{array}{lllllllllllllllll}}
\newcommand{\beno}{\begin{equation}\begin{array}{lllllllllllll}}
\newcommand{\ee}{\end{array}\end{equation}}
\newcommand{\bel}{\begin{equation}\begin{array}{lllllllllllll}\nonumber}
\newcommand{\eel}{\Box\end{array}\end{equation}}
\newcommand{\bi}{\begin{itemize}}
\newcommand{\ei}{\end{itemize}}
\newcommand{\ben}{\begin{enumerate}}
\newcommand{\een}{\end{enumerate}}
\newcommand{\alert}{\textcolor{black}}
\newcommand{\dis}{\displaystyle}
\newcommand{\dsum}{\displaystyle\sum\limits}

\newcommand{\dprod}{\displaystyle\prod\limits}

\newcommand{\mM}{\mbox{$\mathbb{M}$}}

\newcommand{\mS}{\mathbb{S}}

\newcommand{\mG}{\mathbb{G}}
\newcommand{\mA}{\mathscr{A}}

\newcommand{\mE}{\mathscr{E}}
\newcommand{\s}{\vspace{0.175cm}}
\newcommand{\bx}{\bm{x}}

\newcommand{\bA}{\bm{A}}

\newcommand{\bX}{\bm{X}}
\newcommand{\mX}{\mathbb{X}}

\newcommand{\mC}{\mathbb{C}}

\newcommand{\mZ}{\mathbb{Z}}

\newcommand{\mB}{\mathscr{B}}
\newcommand{\BB}{\mathscr{B}}

\newcommand{\mL}{\mathscr{L}}

\newcommand{\mU}{\mathbb{U}}

\newcommand{\bz}{\bm{z}}

\newcommand{\bmu}{\mbox{\boldmath$\mu$}}

\newcommand{\bTheta}{\bm\Theta}
\newcommand{\bta}{\boldsymbol{\eta}}

\newcommand{\btheta}{\boldsymbol{\theta}}

\newcommand{\mbV}{\mathbb{V}}

\newcounter{comment}
\newenvironment{comment}[1][]{\refstepcounter{comment}\par\noindent%
\textbf{Comment~\thecomment #1}:\\ \rmfamily}{}
\newcommand{\bc}{\begin{comment}\em}
\newcommand{\ec}{\end{comment}}

\newcounter{ex}

\newcounter{counterexample}

\newcounter{definition}
\newenvironment{definition}[1][]{\refstepcounter{definition}\par\medskip\indent%
\textbf{Definition.} \rmfamily}{\medskip}

\newcounter{theorem}
\newenvironment{theorem}[1][]{\refstepcounter{theorem}\par\smallskip\indent%
\textbf{Theorem~\thetheorem #1}.\em \rmfamily}{}

\newcounter{proposition}
\newenvironment{proposition}[1][]{\refstepcounter{proposition}\par\smallskip\indent%
\textbf{Proposition~\theproposition #1}.\em \rmfamily}{}

\newcounter{ttproof}
\newcommand{\ttproof}{%
\addtocounter{ttproof}{1}%
{\medskip}\textbf{Proof of Theorem}{}
}

\newcounter{ttproofs}
\newcounter{pproof}
\newcommand{\pproof}{%
\addtocounter{pproof}{1}%
{\medskip}\textbf{Proof of Proposition}{}
}

\newcounter{corollary}
\newenvironment{corollary}[1][]{\refstepcounter{corollary}\par\smallskip\indent%
\textbf{Corollary~\thecorollary #1}.\em \rmfamily}{}

\newcounter{ccproof}
\newcommand{\ccproof}{%
\addtocounter{ccproof}{1}%
{\medskip}\textbf{Proof of Corollary}{}
}

\newcounter{ccproofs}
\newcounter{llproof}
\newcommand{\llproof}{%
\addtocounter{llproof}{1}%
{\smallskip}\textbf{Proof of Lemma}{}
}

\newcounter{lemma}
\newenvironment{lemma}[1][]{\refstepcounter{lemma}\par\smallskip\noindent%
\textbf{Lemma~\thelemma #1}.\em \rmfamily}{}

\newcounter{example}

\newcounter{com}
\newenvironment{com}[1][]{\refstepcounter{com}\par\s\indent%
\textit{Remark~\thecom #1}. \rmfamily}{\medskip}

\newcounter{assumption}
\newenvironment{assumption}[1][]{\refstepcounter{assumption}\par\medskip\noindent%
\textit{Condition C.\theassumption #1}. \rmfamily}{\medskip}

\newcommand{\czero}{[C.1]}
\newcommand{\cone}{[C.2]}
\newcommand{\ctwo}{[C.3]}
\newcommand{\cthree}{\mbox{[C.4]}}
\newcommand{\cfour}{[C.5]}

\renewcommand{\mM}{\mathbb{M}}
\newcommand{\mMs}{\mM(\alpha)}
\newcommand{\etaspace}{\bm{\Xi}}
\newcommand{\fullspace}{\mathbb{N}}
\renewcommand{\uu}{x_i}
\renewcommand{\vv}{x_i^\star}

\renewcommand{\bA}{\bm{A}}
\usepackage{multicol}

\newcommand{\conv}{\mbox{conv}}
\renewcommand{\bx}{\bm{x}}
\renewcommand{\bX}{\bm{X}}

\renewcommand{\bz}{\bm{z}}

\renewcommand{\=}{&=&}

\renewcommand{\mG}{\mathbb{X}(\alpha)}
\newcommand{\mBs}{\mB^\star}
\newcommand{\GG}{\mathscr{G}}
\newcommand{\uuushort}{u(n)}
\newcommand{\uuu}{|\ell(\bthetas, \bzs; \bmu^\star)|}

\newcommand{\uuulong}{C\, K^{1+\beta}\, \norm{\mA}_\infty}

\newcommand{\longtitle}{Consistent Structure Estimation of Exponential-Family Random Graph Models With Block Structure}
\newcommand{\shorttitle}{Consistent Structure Estimation}

\begin{document}

\begin{frontmatter}

\title{\longtitle}
\runtitle{\shorttitle}

\begin{aug}
\author{\fnms{Michael} \snm{Schweinberger}\ead[label=e1]{m.s@rice.edu}}
\affiliation{Rice University}
\address{
Michael Schweinberger\\
Department of Statistics\\
Rice University\\
6100 Main St, MS-138\\
Houston, TX 77005-1827\\
E-mail:\ m.s@rice.edu
}
\end{aug}

\runauthor{Michael Schweinberger}

\begin{abstract}
We consider the challenging problem of statistical inference for exponential-family random graph models based on a single observation of a random graph with complex dependence. To facilitate statistical inference, we consider random graphs with additional structure in the form of block structure. We have shown elsewhere that when the block structure is known, it facilitates consistency results for $M$-estimators of canonical and curved exponential-family random graph models with complex dependence, such as transitivity. In practice, the block structure is known in some applications (e.g., multilevel networks), but is unknown in others. When the block structure is unknown, the first and foremost question is whether it can be recovered with high probability based on a single observation of a random graph with complex dependence. \alert{The main consistency results of the paper show that it is possible to do so under weak dependence and smoothness conditions.} These results confirm that exponential-family random graph models with block structure constitute a promising direction of statistical network analysis.
\end{abstract}


\end{frontmatter}

\section{Introduction}
\label{sec:introduction}

Exponential-family random graph models \citep[][]{FoSd86,WsPp96,HuHa04,SnPaRoHa04,Kr11} are models of network data,
such as disease transmission networks,
insurgent and terrorist networks,
social networks, 
and the World Wide Web \citep[][]{ergm.book}.
Such models can be viewed as generalizations of Bernoulli random graphs with independent edges \citep{Gi59,ErRe60} to random graphs with dependent edges.
Exponential-family random graph models are popular among network scientists \citep[][]{ergm.book},
because network data are dependent data and exponential-family random graph models enable network scientists to model a wide range of dependencies found in network data.

Exponential-family random graph models of dependent network data were pioneered by \citep{FoSd86}.
The models of \citep{FoSd86} and more general models \citep[][]{WsPp96,HuHa04,SnPaRoHa04,Kr11} are discrete exponential families of densities with countable support $\mX$---the set of possible graphs with $n$ nodes and binary or non-binary, count-valued edges---of the form
\be
\label{local.ergm}
p_{\bta}(\bx)
&=& \exp\left(\langle\bta,\, s(\bx)\rangle - \psi(\bta)\right),
& \bx \in \mX,
\ee
where $\langle\bta,\, s(\bx)\rangle$ denotes the inner product of a vector of natural parameters $\bta \in \{\bta \in \mR^{\dim(\bta)}: \psi(\bta) < \infty\}$ and a vector of sufficient statistics $s: \mX \mapsto \mR^{\dim(\bta)}$ and $\psi(\bta)$ ensures that $\sum_{\bx^\prime \in \mX} p_{\bta}(\bx^\prime) = 1$.

\alert{In general,
statistical inference for exponential-family random graph models is challenging \citep{Ha03,BaBrSl11,Sc09b,ChDi11,ShRi11},
because exponential-family random graph models induce complex dependence \citep[e.g., transitivity,][]{ergm.book} and many network data sets either consist of a single observation of a population graph or subgraphs sampled from the population graph.
For example,
epidemiologists studying the spread of infectious diseases (e.g., HIV, Ebola) may be able to observe whether population members were in contact during an epidemic,
but may not be able to obtain independent or repeated observations of contacts over time.
As a result,
the epidemiologists may have to be content with a single observation of the population contact network of interest or subgraphs sampled from the population contact network.
The fact that many network data sets consist of a single observation of a population graph or sampled subgraphs means that concentration and consistency results cannot be obtained along the lines of classical and high-dimensional statistics,
which rely on independent observations from the same source (in a well-defined sense).

In addition,
the complex dependence induced by these models implies that establishing concentration, consistency, and weak convergence results for estimators requires concentration-of-measure results for dependent random variables,
which are more challenging than concentration-of-measure results for independent random variables \citep[e.g.,][]{KoRa08}.}

\subsection{Advantages of block structure}

\alert{While statistical inference for exponential-family random graph models is challenging,
statistical inference for models with additional structure has advantages.}

To demonstrate the advantages of additional structure,
we consider a natural form of additional structure known as block structure.
Block structure is popular in the large and growing body of literature on stochastic block models \citep[e.g.,][]{NkSt01,BiCh09,ABFX08,ChWoAi12,CeDaLa11,RoChYu11,BiChLe11,zhao2012,AmChBiLe13,MoNeSl15,LeRi13,RoQiFa13,Gao15,Jin15,ZhZh16,BiVoRo17}.
We focus here on exponential-family random graph models with block structure,
which allow edges within blocks to be dependent \citep{ScHa13}.
Such models are less restrictive than stochastic block models \citep{NkSt01},
which assume that edges within blocks are independent Bernoulli random variables.
Indeed, 
sensible specifications of exponential-family random graph models can capture excesses in transitivity and many other interesting features of random graphs that induce complex dependence among edges within blocks \citep{ScHa13}.
We have shown elsewhere that when the block structure is known,
exponential-family random graph models with block structure have important advantages:
\bi
\item If edges depend on other edges within the same block but do not depend on edges outside of the block,
models induce local dependence within blocks.
Local dependence makes sense in applications,
because network data are dependent data but network dependence is more local than global \citep{PpRg02,ScHa13}.
\item Models with block structure are weakly projective in the sense that the probability mass function of a random graph with block structure is consistent with the probability mass function of a larger random graph with more blocks \citep{ScHa13,ScSt16},
whereas many models without block structure are not projective \citep{Sn10,ShRi11,CrDe15,LaRiSa17}.
\item Local dependence induces weak dependence as long as the blocks are not too large.
Weak dependence facilitates concentration and consistency results for $M$-estimators,
including maximum likelihood estimators \citep{ScSt16}.
These results are of fundamental importance,
because they are the first consistency results for models with transitivity and other interesting features of random graphs that induce complex dependence.
Transitivity is interesting in practice \citep{WsFk94},
but is challenging from a theoretical point of view \citep[e.g.,][]{ChDi11,ShRi11},
and indeed no other consistency results exist for transitivity.
\ei
In other words, 
block structure is not only useful for community detection in social networks,
for which stochastic block models can be used,
but also facilitates statistical inference for random graphs with complex dependence induced by transitivity and many other interesting features of random graphs.

\subsection{Recovery of unknown block structure}

In some applications,
the block structure is known.
An example is multilevel networks,
which are popular in network science \citep[e.g.,][]{Wa13,ZaLo15,Lo16,slaughter2016multilevel,hollway2016multilevel,hollway2017multilevel}:
e.g.,
the blocks may correspond to school classes in schools,
units of armed forces,
and departments of universities.

While the block structure is known in some applications,
it is unknown in others.
When the block structure is unknown,
the first and foremost question is whether it can be recovered with high probability.
A large and growing body of consistency results for stochastic block models shows that it is possible to recover the block structure of stochastic block models with high probability \citep[e.g.,][]{NkSt01,BiCh09,ABFX08,ChWoAi12,CeDaLa11,RoChYu11,BiChLe11,zhao2012,AmChBiLe13,MoNeSl15,LeRi13,RoQiFa13,Gao15,Jin15,ZhZh16,BiVoRo17}.
While it is encouraging that the block structure of stochastic block models can be recovered with high probability,
these results are restricted to models with independent edges within and between blocks. 
It is not at all obvious whether the block structure of the much more complex exponential-family random graph models can be recovered with high probability.

We show here that consistent recovery of block structure is not limited to stochastic block models,
but is possible for the much more complex exponential-family random graph models.
\alert{The main consistency results of the paper show that it is possible to recover the block structure with high probability under weak dependence and smoothness conditions.
}
Among other things,
these consistency results demonstrate that the conditional independence assumptions underlying stochastic block models are not necessary for consistent recovery of block structure. 
In other words,
these results suggest that it is possible to obtain consistency results for many interesting models with block structure,
both stochastic block models with independent edges within blocks and richer models with dependent edges within blocks,
such as the models and methods proposed by \citep{ScHa13} and \citep{Waetal18}.
\alert{An indepth investigation of all of these models and methods is beyond the scope of a single paper:
each of them is challenging,
owing to the complex dependence within blocks and the wide range of model terms and canonical and curved exponential-family parameterizations.
}
However,
the main consistency results reported here suggest that statistical inference for these models and methods is possible and worth exploring in more depth.

The paper is structured as follows.
Section \ref{sec:model} introduces exponential-family random graph models with additional structure in the form of block structure.
Section \ref{sec:structure} discusses the main consistency results.
Section \ref{sec:simulations} presents simulation results. 
Section \ref{sec:proofs} proves the main consistency results.

\subsection{Other, related literature}

It is worth noting that two broad classes of exponential-family random graph models can be distinguished based on the underlying dependence assumptions:
one class of models assumes that edges or pairs of directed edges are independent \citep[e.g., the $\beta$-model and the $p_1$-model,][]{HpLs81,DiChSl11,RiPeFi13,YaLeZh11,YaZhQi15,YaWaQi16}, 
while the other class of models allows edges or pairs of directed edges to be dependent \citep[][]{FoSd86,SnPaRoHa04,HuHa04}.
The independence assumptions of the first class of models are restrictive,
because it is known that edges in real-world networks tend to depend on other edges \citep[][]{HpLs76}.
The dependence assumptions of the second class of models are problematic,
because some of these models allow edges to depend on many other edges: 
e.g.,
the conditional independence assumptions of \citep{FoSd86} allow the conditional distribution of each edge variable to depend on $2\, (n - 2)$ other edge variables.
Some---but not all---of these models induce strong dependence in large random graphs and therefore have undesirable properties,
such as model near-degeneracy \citep{Ha03,BaBrSl11,Sc09b,ChDi11,ShRi11,CrDe15}. 
Exponential-family random graph models with block structure strike a middle ground between models with independence assumptions and models with strong dependence assumptions,
because sensible specifications of these models induce weak dependence.
As a consequence, 
sensible specifications of these models have desirable properties, 
as explained above.

\section{Exponential-family random graph models with additional structure}
\label{sec:model}

\alert{In general,
statistical inference for exponential-family random graph models is challenging,
as discussed in Section \ref{sec:introduction}.}
We facilitate statistical inference by endowing exponential-family random graph models with additional structure that induces weak dependence and hence facilitates consistency results.

Throughout,
we consider random graphs with a set of nodes $\mA = \{1, \dots, n\}$ and a set of edges $\mE \subseteq \mA \times \mA$,
where edges between pairs of nodes $(i, j) \in \mA \times \mA$ are regarded as random variables $X_{i,j}$ with countable sample spaces $\mX_{i,j}$. 
We focus on undirected graphs without self-edges---i.e.,
$X_{i,i} = 0$ and $X_{i,j} = X_{j,i}$ with probability $1$---but extensions to directed random graphs are straightforward.
We write $\bX = (X_{i,j})_{i<j}^n$ and $\mX = \bigtimes_{i<j}^n \mX_{i,j}$.

\alert{To facilitate statistical inference,
we assume that the random graph is endowed with additional structure in the form of a partition of the set of nodes $\mA$ into $K \geq 2$ subsets of nodes $\mA_1, \dots, \mA_K$,
called blocks.}
\alert{To obtain concentration and consistency results,
it is important that the additional structure induces weak dependence,
because strong dependence can make concentration results impossible \citep[e.g.,][]{KoRa08}.
}
We induce weak dependence by restricting dependence to within-block subgraphs $\bX_{k,k} = (X_{i,j})_{i\in\mA_k\, <\, j\in\mA_k}$ ($k = 1, \dots, K$).
The resulting exponential families induce a form of local dependence defined as follows \citep{ScHa13}.
\definition
{\bf Exponential families with local dependence.}
{\em
An exponential family of densities of the form \eqref{local.ergm} with countable support $\mX$ satisfies local dependence as long as its densities satisfy
\be
\label{example1}
p_{\bta}(\bx) 
&=& \dprod_{k = 1}^K p_{\bta}(\bx_{k,k})\; \dprod_{l=1}^{k-1}\; \dprod_{i\in\mA_k,\; j\in\mA_l} p_{\bta}(x_{i,j})
& \mbox{for all } \bx \in \mX.
\ee
}

We give examples of canonical and curved exponential families with local dependence in Sections \ref{class2} and \ref{class1},
respectively.
We discuss the well-known, but restrictive special case of stochastic block models in Section \ref{class3} and demonstrate the added value of exponential families with local dependence relative to stochastic block models in Section \ref{sec:comparison}.

\subsection{Example: canonical exponential families with local dependence}
\label{class2}

An example of canonical exponential families with local dependence and support $\mX = \{0, 1\}^{{n \choose 2}}$ is given by exponential families with block-dependent edge and transitive edge terms of the form
\beno
\label{geo0}
p_{\bta}(\bx)
&\propto& \exp\left(\dsum_{k \leq l}^K \eta_{1,k,l} \dsum_{i \in \mA_k,\; j \in \mA_l} x_{i,j} + \dsum_{k=1}^K \eta_{2,k,k}\, s_{k,k}(\bx)\right),
\ee
where
\beno
s_{k,k}(\bx)
\= \dsum_{i \in \mA_k\, <\, j \in \mA_k} x_{i,j}\, \one_{i,j}(\bx).
\ee
Here, 
$\one_{i,j}(\bx) = 1$ if the number of shared partners of nodes $i \in \mA_k$ and $j \in \mA_k$ in block $\mA_k$ satisfies $\sum_{h \in \mA_k,\, h \neq i,j} x_{h,i}\, x_{h,j} > 0$ and $\one_{i,j}(\bx) = 0$ otherwise.
If $x_{i,j}\, \one_{i,j}(\bx) = 1$,
the edge between nodes $i$ and $j$ is called transitive.
We note that in recent work \citep[][]{Kr11,HuKrSc12,KrKo14,ScSt16} transitive edge terms have turned out to be attractive alternatives to the triangle terms which have been used since the classic work of \citep{FoSd86} but which possess undesirable properties \citep{Ha03,Sc09b,ChDi11}.

\subsection{Example: curved exponential families with local dependence}
\label{class1}

An example of curved exponential families with local dependence and support $\mX = \{0, 1\}^{{n \choose 2}}$ is given by exponential families with block-dependent edge and geometrically weighted edgewise shared partner terms of the form
\beno
\label{geo0}
p_{\bta}(\bx)
&\propto& \exp\left(\dsum_{k \leq l}^K \eta_{1,k,l} \dsum_{i \in \mA_k,\; j \in \mA_l} x_{i,j} + \dsum_{k=1}^K \dsum_{t=1}^{|\mA_k|-2} \eta_{2,k,k,t}\; s_{k,k,t}(\bx)\right),
\ee
where
\beno
s_{k,k,t}(\bx)
\= \dsum_{i \in \mA_k\, <\, j \in \mA_k} x_{i,j}\, \one_{i,j,t}(\bx).
\ee
Here,
$\one_{i,j,t}(\bx) = 1$ if the number of shared partners of nodes $i \in \mA_k$ and $j \in \mA_k$ in block $\mA_k$ satisfies $\sum_{h \in \mA_k,\, h \neq i,j} x_{h,i}\, x_{h,j} = t$ and $\one_{i,j,t}(\bx) = 0$ otherwise.
A curved exponential-family parameterization is given by
\be
\label{curvedpar}
\eta_{1,k,l}(\btheta) 
\= \theta_{1,k,l}\s
\\
\eta_{2,k,k,t}(\btheta)
\= \theta_{2,k} \left\{\theta_{3,k}\, \left[1 - \left(1 - \dfrac{1}{\theta_{3,k}}\right)^t\right]\right\},
& \theta_{3,k}\; >\; \dfrac12.
\ee
Such terms are called geometrically weighted edgewise shared partner terms \citep{HuHa04,HuGoHa08},
because the natural parameters $\eta_{2,k,k,t}(\btheta)$ are based on the geometric sequence $(1 - 1\, /\, \theta_{3,k})^t$, 
$t = 1, 2, \dots$
It is worth noting that the corresponding geometric series converges as long as $\theta_{3,k} > 1 / 2$ and that $\theta_{3,k} \leq 1 / 2$ is problematic on probabilistic and statistical grounds \citep{Sc09b,ScSt16}.
The parameterization is called a curved exponential-family parameterization,
because the natural parameter vector $\bta(\btheta)$ is a non-affine function of a lower-dimensional parameter vector $\btheta$;
see Remark \ref{com.1234} in Section \ref{mainresults}.
Last, 
but not least, 
note that in the special case $\theta_{3,k} = 1$ ($k = 1, \dots, K$) the curved exponential family reduces to the canonical exponential family described in Section \ref{class2}.
\hide{
We note that $\dim(\bta) \tends \infty$ as $n \tends \infty$ provided that the number of blocks $K$ increases as a function of the number of nodes $n$.
}

\subsection{Example: stochastic block models}
\label{class3}

A well-known, 
but restrictive special case of exponential families with local dependence and support $\mX = \{0, 1\}^{{n \choose 2}}$ are stochastic block models \citep[][]{NkSt01}.
Stochastic block models assume that all edge variables $X_{i,j}$ are independent given the block structure,
which implies that $p_{\bta}(\bx)$ can be written as
\beno
\label{geo00}
p_{\bta}(\bx)
&\propto& \exp\left(\dsum_{k \leq l}^K \eta_{1,k,l} \dsum_{i \in \mA_k,\; j \in \mA_l} x_{i,j}\right),
\ee
where $\eta_{1,k,l}$ is the log odds of the probability of an edge between nodes in blocks $\mA_k$ and $\mA_l$.
\hide{
While stochastic block models can be used to detect communities in networks,
such models fail to capture a wide range of dependencies encountered in networks---such as transitivity \citep[][]{ergm.book}---and are therefore misspecified models:
see, e.g., the discussion of \citep{Sn07a}.
}

\subsection{Added value of exponential families with local dependence}
\label{sec:comparison}

\alert{Exponential families with local dependence can capture many features of random graphs within blocks,
in contrast to stochastic block models,
and can therefore be worth the additional costs in terms of model complexity.

To demonstrate the added value of exponential families with local dependence compared with stochastic block models,
first note that many network data sets show evidence of systematic deviations from models which assume that edges are independent,
as has been well-documented since the 1970s \citep[see, e.g.,][]{Ra53a,Ra53b,HpLs76}.
Stochastic block models assume that edges are independent within and between blocks and hence cannot capture such systematic deviations from independence.
For example,
suppose that $\bx \in \mX$ is observed and the block structure is known,
and let $s_{1,k,k}(\bx)$ be the observed number of edges and $s_{2,k,k}(\bx)$ be the observed number of transitive edges in block $\mA_k$ ($k = 1, \dots, K$). 
A helpful observation for comparing exponential families with local dependence and stochastic block models is that stochastic block models are special cases of exponential families with local dependence and natural parameter vectors $\bta_{k,k} = (\eta_{1,k,k},\, \eta_{2,k,k}) = (\eta_{1,k,k},\, 0)$---as described in Section \ref{class2}---where $\eta_{1,k,k}$ and $\eta_{2,k,k}$ are the natural edge and transitive edge parameter of block $\mA_k$,
respectively.
If the natural parameter vector $\bta_{k,k} = (\eta_{1,k,k},\, 0)$ of block $\mA_k$ is estimated by the maximum likelihood estimator $\widehat\bta_{k,k} = (\widehat\eta_{1,k,k},\, 0)$ under stochastic block models with known block structure,
then the maximum likelihood estimator solves
\beno
\mbE_{\widehat\eta_{1,k,k},\, \eta_{2,k,k}=0}\; s_{1,k,k}(\bX)
&=& s_{1,k,k}(\bx),
& k = 1, \dots, K,
\ee
provided the maximum likelihood estimator exists \citep{Ha03,fienberg-2008}.
However,
network data sets may have many more transitive edges within blocks than expected under stochastic block models.
In other words,
we may observe that
\beno
s_{2,k,k}(\bx)
&\gg& 
\mbE_{\widehat\eta_{2,k,k},\, \eta_{2,k,k}=0}\; s_{2,k,k}(\bX)
& \mbox{ for some or all } 
& k \in \{1, \dots, K\}.
\ee
To capture such systematic deviations from stochastic block models,
exponential families with local dependence can be useful.
To see that,
note that classic exponential-family theory \citep[][Corollary 2.5, p.\ 37]{Br86} implies that,
for any $\eta_{2,k,k}>0$,
\beno
\mbE_{\eta_{1,k,k},\, \eta_{2,k,k}>0}\; s_{2,k,k}(\bX)
&>& \mbE_{\eta_{1,k,k},\, \eta_{2,k,k}=0}\; s_{2,k,k}(\bX),
& k = 1, \dots, K.
\ee
In other words,
the expected number of transitive edges in block $\mA_k$ is greater under exponential families with local dependence with $\eta_{2,k,k} > 0$ than under stochastic block models with $\eta_{2,k,k} = 0$,
assuming that both have the same edge parameters $\eta_{1,k,k}$ ($k = 1, \dots, K$).
As a consequence,
exponential families with local dependence can capture an excess in the expected number of transitive edges within blocks,
relative to stochastic block models.
In fact,
the maximum likelihood estimator $\widehat\bta_{k,k} = (\widehat\eta_{1,k,k},\, \widehat\eta_{2,k,k})$ of block $\mA_k$ under exponential families with local dependence and known block structure solves
\beno
\mbE_{\widehat\eta_{1,k,k},\, \widehat\eta_{2,k,k}}\; s_{1,k,k}(\bX)
&=& s_{1,k,k}(\bx),
& k = 1, \dots, K,\s
\\
\mbE_{\widehat\eta_{1,k,k},\, \widehat\eta_{2,k,k}}\; s_{2,k,k}(\bX)
&=& s_{2,k,k}(\bx),
& k = 1, \dots, K,
\ee
provided the maximum likelihood estimator exists \citep{Ha03,fienberg-2008}.
Thus,
exponential families with local dependence can match both the observed number of edges and transitive edges within blocks,
in contrast to stochastic block models.
As a consequence,
exponential families with local dependence can outperform stochastic block models in terms of transitivity \citep[see, e.g., the empirical results of ][where the blocks are known and correpond to school classes in schools]{StScBoMo18}.

More generally,
exponential families with local dependence can capture many features of random graphs that induce dependence among edges within blocks,
including---but not limited to---transitivity.
The flexibility of the exponential-family framework and its ability to capture many features of random graphs within blocks is one of its greatest advantages.
However,
it is worth noting that not all specifications of exponential-family models with local dependence are equally useful:
e.g.,
it is well-known that exponential-family models with $k$-star and triangle terms can induce undesirable behavior in large random graphs,
such as model near-degeneracy \citep{Jo99,Ha03,Sc09b,ChDi11}.
Thus,
within-block $k$-star and triangle terms can be used as long as the blocks are not too large,
but should not be used when the blocks are large.
Other specifications of exponential-family models are more appropriate for large blocks,
e.g.,
the specifications described in Sections \ref{class2} and \ref{class1}:
each of them implies that the value added by additional triangles to the log odds of the conditional probability of an edge,
given all other edges,
decays \citep[see, e.g.,][]{SnPaRoHa04,HuHa04,Hu08,ScKrBu17}.
By contrast,
models with triangle terms make the implicit assumption that the added value of additional triangles does not decay,
which can lead to undesirable behavior in large random graphs and hence large within-block subgraphs \citep{Jo99,Ha03,Sc09b,ChDi11}.
But the restriction that blocks cannot be too large---which we discuss in Remark \ref{com.size} in Section \ref{sec:structure}---ensures that the effect of less appropriate within-block specifications (such as within-block triangle or $k$-star terms) on the random graph remains limited. 
}

\subsection{Notation}
\label{sec:notation}

Throughout,
$\mbE\, f(\bX)$ denotes the expectation of a function $f: \mX \mapsto \mR$ of a random graph with respect to exponential-family distributions $\mbP$ admitting densities of the form \eqref{example1}.
We write $\mbP \equiv \mbP_{\bta^\star}$ and $\mbE \equiv \mbE_{\bta^\star}$,
where $\bta^\star \in \etaspace \subseteq \mbox{int}(\fullspace)$ denotes the data-generating natural parameter vector and $\etaspace \subseteq \mbox{int}(\fullspace)$ denotes a subset of the interior $\mbox{int}(\fullspace)$ of the natural parameter space $\fullspace = \{\bta \in \mR^{\dim(\bta)}: \psi(\bta) < \infty\}$.
We assume that $\bta: \bTheta \times \mZ \mapsto \etaspace$ is a function of $(\btheta, \bz) \in \bTheta \times \mZ$,
where 
\beno
\bTheta \times \mZ 
\= \{(\btheta, \bz)\; \in\; \mR^{\dim(\btheta)}\, \times\, \{1, \dots, K\}^n:\; \psi(\bta(\btheta, \bz)) < \infty\}.
\ee
Here,
$\btheta$ is a vector of block-dependent parameters of dimension $\dim(\btheta) \leq \dim(\bta)$ while $\bz$ is a vector of block memberships of nodes.
Observe that the natural parameter vectors of the canonical and curved exponential families described in Sections \ref{class2} and \ref{class1} can be represented in this form.
The data-generating values of $(\btheta, \bz) \in \bTheta \times \mZ$ are denoted by $(\bthetas, \zs)$.
The $\ell_1$-, $\ell_2$-, and $\ell_\infty$-norm of vectors are denoted by $\norm{.}_1$, $\norm{.}_2$, and $\norm{.}_\infty$,
respectively.
Uppercase letters $A, B, C > 0$ denote unspecified constants,
which may be recycled from line to line.

\section{Consistent estimation of block structure}
\label{sec:structure}

\alert{We present here the first consistency results which show that it is possible to recover the block structure with high probability under weak dependence and smoothness conditions.
These consistency results are non-trivial,
because we cover exponential families with (a) countable support;
(b) a wide range of dependencies within blocks;
and (c) a wide range of canonical and curved exponential-family parameterizations.
}

To recover the block structure along with the parameters given {an observation} $\bx$ of $\bX$,
we consider the following restricted maximum likelihood estimator:
\beno
\label{reg.ll}
\estimate
&\in& \argmax\limits_{(\btheta, \bz)\, \in\, \bTheta_0 \times \mZ_0} \ell(\btheta, \bz; s(\bx)),
\ee
where
\beno
\ell(\btheta, \bz; s(\bx))
\= \langle\bta(\btheta, \bz),\, s(\bx)\rangle - \psi(\bta(\btheta, \bz))
\ee
denotes the loglikelihood function of $(\btheta, \bz) \in \bTheta_0 \times \mZ_0$ and $\bTheta_0 \times \mZ_0$ is a subset of $\bTheta \times \mZ$ to be specified.
Computational implications are discussed in Section \ref{sec:discussion}.
We assume that the number of blocks $K$ is known 
and that both $\btheta$ and $\bz$ are parameters,
which is commonplace in the special case of stochastic block models \citep[e.g.,][]{BiCh09,ChWoAi12,AmChBiLe13}.
It is worth noting that the maximum likelihood estimator $\estimate$ is not unique,
because the likelihood function is invariant to the labeling of blocks.
All following statements are therefore understood as statements about equivalence classes of block structures.

We call the maximum likelihood estimator $\estimate$ restricted,
because we restrict maximum likelihood estimation to a subset $\bTheta_0 \times \mZ_0$ of $\bTheta \times \mZ$.
We need to do so,
because without additional restrictions exponential families with local dependence can induce strong dependence and smoothness problems.
To motivate the restrictions on $\bTheta \times \mZ$,
it is instructive to discuss the following concentration result,
which is instrumental to deriving the main consistency results of the paper.

\begin{lemma}
\label{proposition.concentration}
Suppose that a random graph is governed by an exponential family with local dependence and countable support $\mX$.
Let $f: \mX \mapsto \mbR$ be Lipschitz with respect to the Hamming metric $d: \mX \times \mX \mapsto \mbR_0^+$ defined by
\beno
d(\bx_1, \bx_2)
\= \dsum_{i<j}^n \one_{x_{1,i,j} \neq x_{2,i,j}},
& (\bx_1, \bx_2) \in \mX \times \mX,
\ee
with Lipschitz coefficient $\norm{f}_{\lip} > 0$ and expectation $\mbE\, f(\bX) < \infty$.
Then there exists $C > 0$ such that,
for all $n > 0$ and all $t > 0$,
\beno
\label{ratio}
\mbP(|f(\bX) - \mbE\, f(\bX)|\; \geq\; t)
\lte 2\, \exp\left(- \dfrac{t^2}{C\, n^2\, \norm{\mA}_\infty^4\, \norm{f}_{\lip}^2}\right),
\ee
where $\norm{\mA}_\infty = \max_{1 \leq k \leq K} |\mA_k| > 0$ denotes the size of the largest data-generating block.
\end{lemma}

\s

The proof of Lemma \ref{proposition.concentration} can be found in the supplementary materials.
The proof relies on concentration of measure inequalities for dependent random variables \citep[][]{KoRa08} and bounds mixing coefficients---which quantify the strength of dependence induced by exponential families with local dependence---in terms of $\norm{\mA}_\infty$.

\alert{Lemma \ref{proposition.concentration} demonstrates that the probability mass of a function $f(\bX)$ of a random graph concentrates around the corresponding expectation $\mbE\, f(\bX)$ as long as the data-generating exponential family induces weak dependence and the function $f(\bX)$ satisfies smoothness conditions.
}
We are interested in applying Lemma \ref{proposition.concentration} to concentrate exponential-family loglikelihood functions of the form $\ell(\btheta, \bz; s(\bX)) = \log p_{\bta(\btheta,\bz)}(\bX)$.
To make sure that the probability mass of $\log p_{\bta(\btheta,\bz)}(\bX)$ concentrates around the expectation $\mbE\, \log p_{\bta(\btheta,\bz)}(\bX)$,
we need to impose additional restrictions on $\mZ$ for at least two reasons.
First of all,
large blocks can induce strong dependence,
which weakens concentration results---as can be seen from the term $\norm{\mA}_\infty$ in Lemma \ref{proposition.concentration}.
Second,
changes of edges in large blocks can give rise to large changes of $\log p_{\bta(\btheta,\bz)}(\bx)$, 
which weakens concentration results as well---as can be seen from the Lipschitz coefficient $\norm{f}_{\lip}$ in Lemma \ref{proposition.concentration}.
Thus,
to deal with strong dependence and smoothness problems, 
restrictions need to be imposed on the sizes of blocks in $\mZ$.
An additional issue is that the unrestricted maximum likelihood estimator fails to exist with non-negligible probability \citep{Ha03,fienberg-2008}.
These observations motivate the following assumptions.

\subsection{Assumptions}
\label{parameterizations}

We assume that the data-generating natural parameter vector $\bta^\star \in \etaspace \subseteq \mbox{int}(\fullspace)$ is in the interior $\mbox{int}(\fullspace)$ of the natural parameter space $\fullspace$,
which implies that the expectation $\mbE\, s(\bX)$ exists \citep[][Theorem 2.2, pp.\ 34--35]{Br86} and so does the expectation $\mbE\; \ell(\btheta, \bz; s(\bX))$,
because
\beno
\mbE\; \ell(\btheta, \bz; s(\bX)) 
\= \langle\bta(\btheta, \bz),\, \mbE\, s(\bX)\rangle - \psi(\bta(\btheta, \bz))
\= \ell(\btheta, \bz; \mbE\, s(\bX)).
\ee
Let $\bmu(\bta) = \mbE_{\bta}\, s(\bX)$ be the mean-value parameter vector of an exponential family with natural parameter vector $\bta \equiv \bta(\btheta, \bz)$ and let $\mM = \rint(\mC)$ be the mean-value parameter space,
where $\rint(\mC)$ is the relative interior of the convex hull $\mC = \conv\{s(\bx): \bx \in \mX\}$ of the set $\{s(\bx): \bx \in \mX\}$.
It is well-known that in minimal exponential families the mapping between the relative interior of the mean-value and natural parameter space is one-to-one \citep[][Theorem 3.6, p.\ 74]{Br86} and that all non-minimal exponential families can be reduced to minimal exponential families \citep[][Theorem 1.9, p.\ 13]{Br86}.
Denote by $\bmu^\star \equiv \bmu(\bta^\star)$ the data-generating mean-value parameter vector.
For any $\alpha > 0$,
let
\beno
\mMs
\= \left\{\bmu \in \mM:\; |\ell(\bthetas, \bzs; \bmu) -  \ell(\bthetas, \bzs; \bmu^\star)|\; <\; \alpha\; |\ell(\bthetas, \bzs; \bmu^\star)|\right\}
\ee
be the subset of mean-value parameter vectors $\bmu \in \mM$ that are close to the data-generating mean-value parameter vector $\bmu^\star \in \mM$ in the sense that\linebreak
$|\ell(\bthetas, \bzs; \bmu) -  \ell(\bthetas, \bzs; \bmu^\star)|\, <\, \alpha\; |\ell(\bthetas, \bzs; \bmu^\star)|$.
The advantage of introducing the subset $\mMs$ of $\mM$ is that the main assumptions stated below can be weakened, because some of them need to hold on $\mMs$,
but need not hold on $\mM \setminus \mMs$.
\hide{
It is worth noting that $\norm{\mA}_\infty \leq \size$ provided $\zs \in \mZ_0$,
but in applications blocks may be small, 
so that $\norm{\mA}_\infty$ may be much smaller than $\size$.
}

The main assumptions can be stated as follows;
note that conditions \cone{ }and \ctwo{ }are assumed to hold on $\mMs$,
but need not hold on $\mM \setminus \mMs$.
\begin{enumerate}
\item[\czero\hspace{-.15cm}] For any fixed $\bz \in \mZ$,
the map $\bta: \bTheta \times \mZ \mapsto \etaspace$ is one-to-one and continuous on $\bTheta$.
\item[\cone\hspace{-.15cm}] For any fixed $\bz \in \mZ$ and any fixed $\bmu \in \mMs$,
the loglikelihood function $\ell(\btheta, \bz; \bmu)$ is upper semicontinuous on $\bTheta$.
\item[\ctwo\hspace{-.15cm}] There exist $A_1 > 0$ and $n_1 > 0$ such that,
for all $n > n_1$,
all $(\btheta_1, \btheta_2) \in \bTheta\times\bTheta$,
all $\bz \in \mZ$,
and all $\bmu \in \mMs$,
\beno
\label{smoothness1}
|\langle \bta(\btheta_1, \bz)-  \bta(\btheta_2, \bz),\, \bmu\rangle|
\lte A_1\, \norm{\btheta_1-\btheta_2}_2\; \uuu.
\ee
\item[\cthree\hspace{-.15cm}] There exist $A_2 > 0$ and $n_2 > 0$ such that,
for all $n > n_2$,
all $(\btheta, \bz) \in \bTheta \times \mZ$,
and all $(\bx_1, \bx_2) \in \mX\times\mX$,
\beno
\label{smoothness2}
|\langle \bta(\btheta, \bz),\, s(\bx_1) - s(\bx_2)\rangle|
\lte A_2\; d(\bx_1, \bx_2)\; L(\bz),
\ee
where $L(\bz)$ is the size of the largest block under $\bz$.
\item[\cfour\hspace{-.15cm}] The data-generating parameters $\truth$ are contained in $\bTheta_0 \times \mZ_0 \subseteq \bTheta \times \mZ$,
where 
\bi
\item[(a)] $\bTheta_0$ has dimension $\dim(\btheta) \leq A\, n$ and can be covered by $\exp(C\, n)$ closed balls $\mB(\btheta_q,\, B)$ with centers $\btheta_q \in \bTheta$ and radius $B > 0$,
i.e.,\linebreak
$\bTheta_0 \subseteq \bigcup_{1 \leq q \leq \exp(C\, n)} \mB(\btheta_q,\, B)$,
where $A, B, C > 0$.
\item[(b)] $\mZ_0$ consists of all block structures for which the size of each of the $K$ blocks is bounded above by $\size$,
where $K$ and $\size$ can increase as a function of the number of nodes $n$.
\ei
\end{enumerate}
Corollaries \ref{c.canonical} and \ref{c.curved} in Section \ref{mainresults} show that conditions \czero---\cthree{ }are satisfied by a wide range of canonical and curved exponential families with local dependence.
\alert{Condition \czero{ }along with the assumption that the exponential family is minimal ensures that $\mbP_{\bta(\btheta_1,\bz)} \neq \mbP_{\bta(\btheta_2,\bz)}$ for all $\btheta_1 \neq \btheta_2$ given $\bz \in \mZ$.}
Conditions \cone---\cthree{ }are smoothness conditions.
Condition \cone{ }is a weak assumption: 
it is well-known that canonical exponential-family loglikelihood functions are upper semicontinuous \citep[][Lemma 5.3, p. 146]{Br86} and it turns out that that the most interesting curved exponential-family loglikelihood functions are upper semicontinuous as well,
which is verified by Corollaries \ref{c.canonical} and \ref{c.curved} in Section \ref{mainresults}.
Condition \ctwo{ }imposes restrictions on how much $\log p_{\bta(\btheta,\bz)}(\bx)$ can change as a function of $\bta(\btheta,\bz)$,
whereas condition \cthree{ }imposes restrictions on how much $\log p_{\bta(\btheta,\bz)}(\bx)$ can change as a function of $\bx$.
Condition \ctwo{ }is stated in terms of $\uuu$ to accomodate both sparse and dense random graphs;
we discuss the notion of sparse and dense random graphs in more detail in Remark \ref{com.sparse} in Section \ref{mainresults}.
Condition \cfour(a) allows the dimension $\dim(\btheta)$ of the parameter space $\bTheta_0$ to increase as a function of the number of nodes $n$
and hence allows the model to be flexible while ensuring that $\bTheta_0$ cannot be too large.
We need these conditions,
because we have {an observation} and therefore cannot use conventional arguments to prove that estimators fall with high probability into compact subsets of the parameter space when the number of observations $N$ is large \citep[e.g.,][]{Be72}.
Condition \cfour(b) complements condition \cthree{ }and helps ensure that $\log p_{\bta(\btheta,\bz)}(\bx)$ is not too sensitive to changes of $\bx$ by restricting the set of block structures to blocks whose size is bounded above by $\size$.
The main consistency results of the paper,
Proposition \ref{p.z.1} and Theorem \ref{theorem.step3} in Section \ref{mainresults},
impose restrictions on $\size$.

\subsection{Main consistency results}
\label{mainresults}

We discuss the main consistency results concerning the recovery of the block structure given {an observation} of a random graph with complex dependence.

The recovery of the block structure is made possible by the following fundamental concentration result.
The concentration result shows that with high probability the distribution parameterized by the restricted maximum likelihood estimator $\estimate$ is close to the distribution parameterized by the data-generating parameters $\truth$ in terms of Kullback-Leibler divergence $KL(\bthetas, \zs;\; \bthetah, \zh) = \ell(\bthetas, \zs; \bmu^\star) - \ell(\bthetah, \zh; \bmu^\star)$ provided that the number of nodes $n$ is sufficiently large.
The result covers a wide range of canonical and curved exponential families with local dependence.

\begin{proposition}
\label{p.z.1}
Suppose that {an observation} of a random graph is generated by an exponential family with local dependence and countable support $\mX$ satisfying conditions \czero---\cfour.
Assume that,
for all $C_1 > 0$,
however large, 
there exists $n_1 > 0$ such that,
for all $n > n_1$,
\be
\label{a.condition}
\uuu \gte C_1\; n^{3/2}\; \norm{\mA}_\infty^2\; \size\, \sqrt{\log n},
\ee
\alert{where $L = \max_{\bz \in \mZ_0} L(\bz)$.}
Then there exist $C > 0$, $C_2 > 0$, and $n_2 > 0$ such that,
for all $n > n_2$, 
with at least probability $1 - 2\, \exp\left(- \alpha^2\, C_2\, n \log n\right)$,
the restricted maximum likelihood estimator $\estimate \in \bTheta_0 \times \mZ_0$ exists and,
for all $\epsilon > 0$,
\beno 
\mbP(KL(\bthetas, \zs; \bthetah, \zh) < \epsilon\, \uuu)
\geq 1 - 4 \exp\left(- \min(\alpha^2, \epsilon^2)\, C\, n \log n\right),
\ee 
where $\alpha > 0$ is identical to the constant $\alpha$ used in the construction of the subset $\mMs$ of the mean-value parameter space $\mM$.
\end{proposition}

\s

\hide{
The proof of Proposition \ref{p.z.1} in Section \ref{sec:proofs} exploits the fact that under conditions \czero---\cfour{ }exponential families with local dependence satisfy weak dependence and smoothness conditions that enable concentration results.
}

The concentration result in Proposition \ref{p.z.1} paves the ground for the main consistency result.
The consistency result is generic and covers a wide range of canonical and curved exponential families with local dependence.
It states that the discrepancy between the estimated and data-generating block structure is small with high probability given {an observation} of a random graph with complex dependence  provided that the number of nodes $n$ is sufficiently large.
To define the discrepancy between the estimated and data-generating block structure,
let $\delta: \mZ \times \mZ \mapsto [0, n]$ be a discrepancy measure that is invariant to the labeling of blocks.
An example is given by $\delta(\zs, \zh) = \min_{\pi} \sum_{i=1}^n \one_{z_i^\star \neq \pi(\widehat{z}_i)}$,
the minimum Hamming distance between $\zs$ and $\zh$,
where the minimum is taken with respect to all possible permutations $\pi$ of $\zh$.
The following consistency result holds for all discrepancy measures $\delta: \mZ \times \mZ \mapsto [0, n]$ satisfying assumption \eqref{id} of the following result.

\begin{theorem}
\label{theorem.step3}
Suppose that {an observation} of a random graph is generated by an exponential family with local dependence and countable support $\mX$ satisfying conditions \czero---\cfour.
If the random graph satisfies assumption \eqref{a.condition} and there exist $C_1 > 0$ and $n_1 > 0$ such that,
for all $n > n_1$ and all $(\btheta, \bz) \in \bTheta_0 \times \mZ_0$,
\be
\label{id}
KL(\bthetas, \zs;\; \btheta, \bz)
\gte \dfrac{\delta(\zs, \bz)\; C_1\; \uuu}{n},
\ee
then there exist $C > 0$, $C_2 > 0$, and $n_2 > 0$ such that,
for all $n > n_2$, 
with at least probability $1 - 2\, \exp\left(- \alpha^2\, C_2\, n \log n\right)$,
the restricted maximum likelihood estimator $\estimate \in \bTheta_0 \times \mZ_0$ exists and,
for all $\epsilon > 0$, 
\beno
\mbP\left(\dfrac{\delta(\zs, \zh)}{n} < \epsilon\right)
\gte 1 - 4\, \exp\left(- \min(\alpha^2, \epsilon^2)\; C\; n \log n\right),
\ee
where $\alpha > 0$ is identical to the constant $\alpha$ used in the construction of the subset $\mMs$ of the mean-value parameter space $\mM$.
\end{theorem}

\s

We discuss implications of Proposition \ref{p.z.1} and Theorem \ref{theorem.step3},
starting with a short comparison with stochastic block models (Remark \ref{com.sbm}) and then discussing assumption \eqref{a.condition} (Remark \ref{com.sparse}) and its implications in terms of the sizes of blocks (Remark \ref{com.size}) and the number of blocks (Remark \ref{com.number}).
We then proceed with a discussion of conditions \czero---\cthree{ }(Remark \ref{com.1234}) and assumption \eqref{id} (Remark \ref{com.id}) and conclude with some comments on parameter estimation (Remark \ref{com.estimation}).
\alert{Last,
but not least,
we discuss the sharpness of the results (Remark \ref{com:sharpness}).
}
\com {\em Comparison with stochastic block models.}
\label{com.sbm}
There is a large and growing body of consistency results on stochastic block models \citep[e.g.,][]{BiCh09,ChWoAi12,CeDaLa11,RoChYu11,BiChLe11,AmChBiLe13,LeRi13,RoQiFa13,Gao15}.
In the language of stochastic block models,
the consistency result in Theorem \ref{theorem.step3} is a weak consistency result in the sense that the discrepancy between the estimated and data-generating block structure is small with high probability.
In contrast to stochastic block models,
we cover exponential families with (a) countable support;
(b) a wide range of dependencies within blocks;
and (c) a wide range of canonical and curved exponential-family parameterizations.
These dependencies and parameterizations make theoretical results more challenging from a statistical point of view,
but more relevant from a scientific point of view.
However,
these results come at a cost: 
in contrast to stochastic block models,
we need to restrict the sizes of blocks from above to deal with strong dependence and smoothness problems,
as we pointed out in the discussion of Lemma \ref{proposition.concentration}.
The restrictions on the sizes of blocks are detailed in Remark \ref{com.size}.
\alert{
\com {\em Assumption \eqref{a.condition}: sparse and dense random graphs.}
\label{com.sparse}
Assumption \eqref{a.condition} of Proposition \ref{p.z.1} and Theorem \ref{theorem.step3} is stated in terms of the absolute value of the expected loglikelihood function $|\mbE\, \ell(\bthetas, \zs; s(\bX))| = \uuu$ to accomodate sparse and dense random graphs.
We first explain why $\uuu$ may be interpreted as the level of sparsity of a random graph,
and then return to assumption \eqref{a.condition}.

To demonstrate that $\uuu$ may be interpreted as the level of sparsity of a random graph,
consider the classic Bernoulli$(\omega)$ random graphs,
under which edges $X_{i,j}$ are independent Bernoulli$(\omega)$ random variables \citep{ErRe59,ErRe60}.
It is natural, and conventional, to use the expected number of edges $\mbE\, \sum_{i<j}^n X_{i,j}$ to quantify the sparsity of Bernoulli$(\omega)$ random graphs,
because $\sum_{i<j}^n X_{i,j}$ is a sufficient statistic for the natural parameter $\theta = \logit(\omega)$ of the canonical exponential-family representation of Bernoulli$(\omega)$ random graphs.
If an exponential family contains more than one natural parameter and one sufficient statistic,
it makes sense to quantify the sparsity of a random graph based on all sufficient statistics:
in fact,
in many applications,
the sufficient statistics are of substantive interest,
because researchers specify exponential-family models of random graphs by specifying sufficient statistics that capture features of random graphs considered relevant (e.g., the number of edges and transitive edges, see Section \ref{sec:comparison}).
The question,
then,
is how the sparsity of a random graph can be quantified based on all sufficient statistics,
i.e.,
all relevant features of the random graph.
The absolute value of the expected loglikelihood function $\uuu$ is a simple choice,
because it is a function of all sufficient statistics and the key to likelihood-based inference.
In the special case of Bernoulli$(\omega)$ random graphs,
$\uuu$ agrees with $\mbE\, \sum_{i<j}^n X_{i,j}$ on the level of sparsity (ignoring logarithmic factors).
If a Bernoulli$(\omega)$ random graph is dense in the sense that $\omega$ does not depend on $n$,
then both $\mbE\, \sum_{i<j}^n X_{i,j}$ and $\uuu$ are of order $n^2$ and hence agree on the level of sparsity.
If a Bernoulli$(\omega_n)$ random graph is sparse in the sense that $\omega_n \tends 0$ as $n \tends \infty$,
then both quantities are smaller:
e.g.,
if $\omega_n = \log n / n$ \citep[the threshold for connectivity of Bernoulli random graphs,][]{Bo98},
then $\mbE\, \sum_{i<j}^n X_{i,j}$ and $\uuu$ are of order $n \log n$ and $n\, (\log n)^2$,
respectively,
so both quantities agree on the level of sparsity up to a logarithmic factor.
We therefore interpret $\uuu$ as the level of sparsity of a random graph,
but note that the mathematical results in Proposition \ref{p.z.1} and Theorem \ref{theorem.step3} hold regardless of how $\uuu$ is interpreted.

To return to assumption \eqref{a.condition},
the above considerations suggest that the random graph can be sparse, but cannot be too sparse in the sense that $\uuu$ cannot be too small.
If,
e.g.,
$\norm{\mA}_\infty$ and $\size$ grow as fast as $(\log n)^{\gamma_1}$ ($\gamma_1 > 0$) and $(\log n)^{\gamma_2}$ ($\gamma_2 > 0$),
respectively,
then $\uuu$ must grow faster than $n^{3/2}\, (\log n)^{2 \gamma_1 + \gamma_2 + 1/2}$.


\hide{
Remark:
First,
note that
\beno
\left|\log \dfrac{\omega}{1-\omega}\right|
\= \left|\log \dfrac{\log n / n}{1 - \log n / n}\right|
\= \left|\log \dfrac{\log n}{n - \log n}\right|\s
\\
&&\= \left|\log\log n - \log (n - \log n)\right|\s
\\
&&\= \log (n - \log n) - \log\log n\s
\\
&&\lte \log (n - \log n)\s
\\
&&\lte \log n
\ee
for all sufficiently large $n$.
Second,
\beno
\uuu
\= \mbE \dsum_{i<j}^n \left[\log \dfrac{\omega}{1-\omega} X_{i,j} + \log(1 - \omega)\right]\s
\\
\= \dfrac{n\, (n-1)}{2}\, \left[\omega\, \log \dfrac{\omega}{1-\omega} + \log(1 - \omega)\right]\s
\\
\lte \dfrac{n\, (n-1)}{2}\, \left[\omega\, \log \dfrac{\omega}{1-\omega} - \omega\right]\s
\\
\= \dfrac{n\, (n-1)}{2}\, \omega\, \left[\log \dfrac{\omega}{1-\omega} - 1\right]\s
\\
\lte \dfrac{n\, (n-1)}{2}\, \omega\, \log \dfrac{\omega}{1-\omega}\s
\\
\lte \dfrac{n\, (n-1)}{2}\, \omega \log n\s
\\
\= \dfrac{n\, (n-1)}{2}\, \dfrac{\log n}{n}\, \log n\s
\\
\= \dfrac{(n-1)}{2}\, (\log n)^2\s
\\
\lte n\, (\log n)^2
\ee
for all sufficiently large $n$.
}

\hide{
As an aside,
$1 - \omega < 1$ implies $\log(1 - \omega) < 0$,
so that $\dots + \log(1 - \omega) < \dots$
}

\hide{
\beno
\uuu
\= \left|\dsum_{i<j}^n \eta_{i,j}\, \mu_{i,j} - \dsum_{i<j}^n\, \log\left(1 + \exp(\eta_{i,j})\right)\right|\s
\\
\lte \dsum_{i<j}^n \left|\eta_{i,j}\, \mu_{i,j}\right| + \dsum_{i<j}^n\left|\log\left(1 + \exp(\eta_{i,j})\right)\right|\s
\\
\lte \dsum_{i<j}^n \left|\eta_{i,j}\, \mu_{i,j}\right| + \dsum_{i<j}^n \left|\exp(\eta_{i,j})\right|\s
\\
\lte \dsum_{i<j}^n \dfrac{C_1 \log n}{n} + \dsum_{i<j}^n \dfrac{C_2}{n}\s
\\
\lte C\, n\, \log n
\ee
}

}

\com {\em Sizes of blocks.}
\label{com.size}
The sizes of blocks in $\mZ_0$ cannot be too large,
because changes of edges in large blocks can give rise to large changes of $\ell(\btheta, \bz; s(\bx)) = \log p_{\bta(\btheta, \bz)}(\bx)$,
which weakens concentration results,
as we pointed out in the discussion of Lemma \ref{proposition.concentration}.
In fact,
assumption \eqref{a.condition} implies that the size $L$ of the largest possible block in $\mZ_0$ must satisfy
\beno
\size
\lte \dfrac{\uuu}{C_1\; n^{3/2}\, \norm{\mA}_\infty^2\, \sqrt{\log n}}.
\ee
Thus,
in the best-case scenario when $\norm{\mA}_\infty$ is small in the sense that $\norm{\mA}_\infty$ grows at most as fast as $(\log n)^{\gamma}$ ($\gamma > 0$),
$\size$ can grow at most as fast as $n^{1/2} /\, (\log n)^{2 \gamma + 1/2}$,
assuming that the random graph is dense.
In the worst-case scenario when $\norm{\mA}_\infty$ grows as fast as $\size$,
$\size$ can grow at most as fast as $(n / \log n)^{1/6}$.

\com {\em Number of blocks.}
\label{com.number}
The fact that the sizes of blocks in $\mZ_0$ are bounded above by $\size$ implies that the number of blocks $K$ is bounded below by $K \geq n\, /\, \size$.
If,
e.g.,
$\size\, \leq\, n^{1/2}\, /\, (\log n)^{2 \gamma + 1/2}$ ($\gamma > 0$),
then $K\, \geq\, n^{1/2}\, (\log n)^{2 \gamma + 1/2}$.
\alert{Compared with stochastic block models,
the number of blocks $K$ needs to grow at least as fast as in the high-dimensional stochastic block models of \citet{ChWoAi12} (ignoring polylogarithmic terms),
where the rate of growth of $K$ is $n^{1/2}$ \citep{ChWoAi12}, 
but $K$ needs not grow as fast as in the highest-dimensional stochastic block models of \citet{RoQiFa13},
where the rate of growth of $K$ is as high as $n$ (ignoring polylogarithmic terms) \citep{RoQiFa13}.
It is worth noting that allowing $K$ to increase as a function of $n$ makes sense in applications:
\citet{leskovec2008community} and others have observed that many real-world networks have small communities,
which suggests that $K$ should increase as a function of $n$,
as \citet*[][p.\ 1883]{RoChYu11} and others have pointed out.
}

\hide{
Compared with stochastic block models,
the number of blocks needs to grow at least as fast as in the high-dimensional stochastic block model of \citep{ChWoAi12},
where $K$ can grow as fast as $n^{1/2}$,
and may have to grow as fast as in the highest-dimensional stochastic block model of \citep{RoQiFa13},
where $K$ grows as fast as $n$ (ignoring polylogarithmic terms).
}
\hide{ 
Remark:
\beno
n 
\= \dsum_{k=1}^K L_k(\bz)
\lte \dsum_{k=1}^K \size
\= K\, \size,
\ee
where $L_k(\bz)$ denotes the size of block $k$ under $\bz \in \mZ_0$.
Hence
\beno
K 
\gte \dfrac{n}{\size}.
\ee
If,
e.g.,
$\size\, \leq\, n^{1/2}\, /\, (\log n)^{\gamma + 1/2}$,
then 
\beno
K 
\gte \dfrac{n}{\size}
\gte \dfrac{n}{n^{1/2} /\, (\log n)^{\gamma + 1/2}}
\= n^{1/2}\, (\log n)^{\gamma + 1/2}.
\ee
}

\com {\em Conditions \czero---\cthree.}
\label{com.1234}
We show that conditions \czero---\cthree{ }are satisfied by a wide range of canonical and curved exponential families with local dependence.
To ease the presentation,
we consider dense random graphs,
but the following results can be extended to sparse random graphs as long as the random graphs are not too sparse;
see Remark \ref{com.sparse}.

We assume here that $\bta: \bTheta \times \mZ \mapsto \etaspace$ is separable in the sense that $\bta(\btheta, \bz) = \bA(\bz)\, \bm{b}(\btheta)$,
where $\bA: \mZ \mapsto \mR^{\dim(\bta) \times \dim(\bm{b})}$ and $\bm{b}: \bTheta \mapsto \mR^{\dim(\bm{b})}$;
note that,
e.g.,
the curved exponential-family parameterization described in Section \ref{class1} is separable,
and so are many other canonical and curved exponential-family parameterizations.
Since $\bta: \bTheta \times \mZ \mapsto \etaspace$ is separable,
$\bA(\bz)$ can be absorbed into the sufficient statistics vector,
so that $\bta: \bTheta \mapsto \etaspace$ can be considered as a function of $\btheta$ and $s: \mX \times \mZ \mapsto \mR^{\dim(\bta)}$ can be considered as a function of $\bx$ and $\bz$.
As a result,
we can write
\beno
\langle\bta(\btheta, \bz),\, \bmu\rangle
\= \langle\bta(\btheta),\, \bmu(\bz)\rangle
\= \dsum_{k \leq l}^K \langle\bta_{k,l}(\btheta),\, \bmu_{k,l}(\bz)\rangle\s
\\
\langle\bta(\btheta, \bz),\, s(\bx)\rangle
\= \langle\bta(\btheta),\, s(\bx, \bz)\rangle
\= \dsum_{k \leq l}^K \langle\bta_{k,l}(\btheta),\, s_{k,l}(\bx, \bz)\rangle,
\ee
where---in an abuse of notation---we write $\bmu(\bz) = \bA(\bz)^\top \bmu$ ($\bmu \in \mMs$) and $s(\bx, \bz) = \bA(\bz)^\top s(\bx)$ ($s(\bx) \in \mMs$).
If,
in addition, 
$\bm{b}(\btheta)$ is an affine function of $\btheta$,
then $\bta(\btheta)$ can be reduced to $\bta(\btheta) = \btheta$ and $\bta_{k,l}(\btheta)$ can be reduced to $\bta_{k,l}(\btheta) = \btheta_{k,l}$ ($k \leq l = 1, \dots, K$),
in which case we call the exponential family canonical,
otherwise we call the exponential family curved.
In the following,
we denote by $L_k(\bz)$ the number of nodes in block $k$ under block structure $\bz \in \mZ_0$.

The following result shows that conditions \czero---\cthree{ }are satisfied by all canonical exponential families with local dependence satisfying reasonable scaling and smoothness conditions.

\begin{corollary}
\label{c.canonical}
Consider canonical exponential families with local dependence and countable support $\mX$.
Assume that $\bta: \bTheta \times \mZ \mapsto \etaspace$ is separable with $\dim(\btheta_{k,l}) < \infty$ ($k \leq l = 1, \dots, K$) and that the random graph is dense.
If there exist $C_1 > 0$, $C_2 > 0$, and $n_0 \geq 1$ such that, for all $n > n_0$,
\bi
\item[{\em [C.$3^\star$]}\hspace{-.1cm}]\; $\norm{\bmu_{k,l}(\bz)}_\infty \leq C_1\, L_k(\bz)\, L_l(\bz)$ for all $\bz \in \mZ_0$ and all $\bmu \in \mMs$ ($k \leq l = 1, \dots, K$);\s
\item[{\em [C.$4^\star$]}\hspace{-.1cm}]\; $\sum_{k \leq l}^K \norm{s_{k,l}(\bx_1, \bz) - s_{k,l}(\bx_2, \bz)}_\infty \leq C_2\, d(\bx_{1}, \bx_{2})\, L(\bz)$\, for all $(\bx_1, \bx_2) \in \mX \times \mX$ and all $\bz \in \mZ_0$;
\ei
then conditions \czero---\cthree{ } are satisfied.
If conditions \cfour{ }and \eqref{id} are satisfied as well,
then the conclusions of Theorem \ref{theorem.step3} hold.
\end{corollary}

\s

\alert{Condition [C.$3^\star$] is satisfied by all between- and within-block sufficient statistics for which the absolute value of the expectation is bounded above by a constant multiple of the number of pairs of nodes between blocks and within blocks,
respectively:
e.g.,
the number of edges and transitive edges within blocks satisfy condition [C.$3^\star$] and so do all other sufficient statistics that count the number of pairs of nodes within blocks having specified properties or being related to other nodes in the same block in some specified form.
}
Condition [C.$4^\star$] is satisfied by most sufficient statistics,
including the number of edges and transitive edges.

We turn to curved exponential families with local dependence.
We consider curved exponential families of densities of the form
\be
\label{geo}
p_{\bta(\btheta, \bz)}(\bx)
&\propto& \exp\left(\langle\bta(\btheta),\, s(\bx, \bz)\rangle\right),
\ee
where
\beno
\langle\bta(\btheta),\, s(\bx, \bz)\rangle
= \dsum_{k \leq l}^K \eta_{1,k,l}(\btheta) \dsum_{i, j:\, z_i = k,\, z_j = l} x_{i,j} + \dsum_{k=1}^K \dsum_{t=1}^{T_k} \eta_{2,k,k,t}(\btheta)\, s_{k,k,t}(\bx, \bz),
\ee
where $s_{k,k,t}(\bx, \bz)$ are sufficient statistics that induce dependence within blocks (e.g., in case $\mX = \{0, 1\}^{{n \choose 2}}$, $s_{k,k,t}(\bx, \bz)$ may be the number of pairs of nodes with $t$ edgewise shared partners in block $k$).
Here,
the natural parameters are given by
\beno
\eta_{1,k,l}(\btheta)
\= \theta_{1,k,l}\s
\\
\eta_{2,k,k,t}(\btheta)
\= \theta_{2,k} \left\{\theta_{3,k}\, \left[1 - \left(1 - \dfrac{1}{\theta_{3,k}}\right)^t\right]\right\},
& \theta_{3,k}\; >\; \dfrac12,
& T_k \geq 2.
\ee
The following result shows that as long as the underlying geometric series converges,
i.e.,
as long as $\theta_{3,k} > 1 / 2$ ($k = 1, \dots, K$),
conditions \czero---\cthree{ }are satisfied.
The result can be extended to other model terms,
e.g., 
covariate terms.

\begin{corollary}
\label{c.curved}
Consider curved exponential families of the form \eqref{geo} with local dependence and countable support $\mX$.
Assume that $\bta: \bTheta \times \mZ \mapsto \etaspace$ is separable and that there exists $B > 1 / 2$ such that $1 / 2 < \theta_{3,k} < B$ ($k = 1, \dots, K$) and that the random graph is dense.
If there exist $C_1 > 0$, $C_2 > 0$, and $n_0 \geq 1$ such that,
for all $n > n_0$,
\bi
\item[{\em [C.$3^{\star\star}$]}\hspace{-.15cm}] $\sum_{t = 1}^{T_k} |\mu_{k,k,t}(\bz)| \leq C_1\, {L_k(\bz) \choose 2}$ for all $\bz \in \mZ_0$,
where $\mu_{k,k,t}(\bz) = \mbE\, s_{k,k,t}(\bX, \bz)$;\s
\item[{\em [C.$4^{\star\star}$]}\hspace{-.15cm}] $|\sum_{t=1}^{T_k} s_{k,k,t}(\bx_1, \bz) - \sum_{t=1}^{T_k} s_{k,k,t}(\bx_2, \bz)| \leq C_2\, d(\bx_{1,k,k}, \bx_{2,k,k})\, L(\bz)$ for all $(\bx_{1,k,k}, \bx_{2,k,k})$ $\in$ $\mX_{k,k}(\bz) \times \mX_{k,k}(\bz)$ and all $\bz \in \mZ_0$,
where $\mX_{k,k}(\bz)$ denotes the set of all possible within-block subgraphs of block $k$ under $\bz \in \mZ_0$ ($k = 1, \dots, K$);
\ei
then conditions \czero---\cthree{ }are satisfied.
If conditions \cfour{ }and \eqref{id} are satisfied as well,
then the conclusions of Theorem \ref{theorem.step3} hold.
\end{corollary}

\s

The most popular curved exponential families with geometrically weighted terms \citep{SnPaRoHa04,HuHa04,HuGoHa08} satisfy conditions [C.$3^{\star\star}$] and [C.$4^{\star\star}$] of Corollary \ref{c.curved}.
Consider,
e.g.,
geometrically weighted edgewise shared partner terms.
In the case of geometrically weighted edgewise shared partner terms,
$T_k = L_k(\bz) - 2$ and $\sum_{t=1}^{T_k} s_{k,k,t}(\bx, \bz)$ is the number of transitive edges in block $k$,
hence conditions [C.$3^{\star\star}$] and [C.$4^{\star\star}$] are satisfied.

\com {\em Assumption \eqref{id}.}
\label{com.id}
Assumption \eqref{id} of Theorem \ref{theorem.step3} states that the Kullback-Leibler divergence of the distribution parameterized by $(\btheta, \bz)$ from the distribution parameterized by $\truth$ must increase with the discrepancy measure $\delta(\zs, \bz)$.
In the special case of stochastic block models,
\citep{ChWoAi12} and \citep{RoQiFa13} verified identifiability assumption \eqref{id} using the number of misclassified nodes as defined by \citep{ChWoAi12} as a discrepancy measure,
where the number of blocks can grow as fast as $n^{1/2}$ \citep{ChWoAi12} and as fast as $n$ (ignoring polylogarithmic terms) \citep{RoQiFa13},
respectively.
In general,
an application of the mean-value theorem to the expected loglikelihood function 
$\ell(\bta^\star; \bmu^\star) = \langle\bta^\star,\, \bmu^\star\rangle - \psi(\bta^\star)$ shows that,
for all $\bta \in \etaspace \subseteq \mbox{int}(\fullspace)$,
\beno
KL(\bta^\star;\, \bta)
\= \ell(\bta^\star;\, \bmu^\star) - \ell(\bta;\, \bmu^\star)
\= \langle\bta^\star - \bta,\, \bmu(\bta^\star) - \bmu(\dot\bta)\rangle,
\ee
where 
$\dot\bta = \lambda\, \bta^\star + (1 - \lambda)\, \bta$\, ($0 \leq \lambda \leq 1$);
note that $\dot\bta \in \mbox{int}(\fullspace)$ since $\bta^\star \in \mbox{int}(\fullspace)$ and $\bta \in \mbox{int}(\fullspace)$ and the natural parameter space $\fullspace$ is convex.
Therefore, 
assumption \eqref{id} is satisfied as long as changes of blocks give rise to large enough changes of mean-value and natural parameter vectors.
\hide{
Remark:
\beno
\ell(\bta^\star;\, \bmu^\star)
\= \ell(\bta;\, \bmu^\star) + \langle\bta^\star-\bta,\, \nabla_{\bta}\, \ell(\bta;\, \bmu^\star)\rangle + \dfrac12\, (\bta^\star - \bta)^\top\, \nabla_{\dot\bta}^2\, \ell(\bta;\, \bmu^\star)|_{\bta=\dot\bta}\; (\bta^\star - \bta)\s
\\
\= \ell(\bta;\, \bmu^\star) + \langle\bta^\star - \bta,\, \bmu^\star - \bmu\rangle - \dfrac12\, (\bta^\star - \bta)^\top\, \mbV_{\dot\bta}\, s(\bX)\; (\bta^\star - \bta)
\ee
using
\beno
\nabla_{\bta}\, \ell(\bta;\, \bmu^\star)
\= \nabla_{\bta}\, \left[\langle\bta,\, \bmu^\star\rangle - \psi(\bta)\right]\s
\\
\= \nabla_{\bta}\, \langle\bta,\, \bmu^\star\rangle - \nabla_{\bta}\,  \psi(\bta)\s
\\
\= \bmu^\star - \bmu
\ee
and
\beno
\nabla_{\bta}^2\, \ell(\bta;\, \bmu^\star)|_{\bta=\dot\bta}
\= \nabla_{\bta}^2\, \left[\langle\bta,\, \bmu^\star\rangle - \psi(\bta)\right]\s
\\
\= -\nabla_{\bta}^2\, \psi(\bta)\s
\\
\= -\mbV_{\dot\bta}\, s(\bX)
\ee
}

\com {\em Estimation of parameters.}
\label{com.estimation}
The restricted maximum likelihood estimator estimates the parameter vector $\btheta$ along with the block structure $\bz$.
We leave the study of the properties of estimators of $\btheta$ to future research,
but it is worth noting the following.
If the blocks are known \citep[e.g., in multilevel networks,][]{multilevelnetwork},
$M$-estimators of canonical and curved exponential-family random graph models with local dependence and growing blocks are consistent under weak conditions \citep{ScSt16}.
If the blocks are unknown,
$M$-estimators may not be consistent estimators of the data-generating parameters.
Indeed,
it is not too hard to see that,
for any $\bz \neq \zs$---where $\bz \in \mZ_0$ may be an estimate of $\zs \in \mZ_0$---the estimator
\beno
\bthetah(\bz)
\= \argmax\limits_{\btheta\, \in\, \bTheta_0} \left[\ell(\btheta, \bz; s(\bx)) - \ell(\bthetas, \bzs; s(\bx))\right]
\ee
estimates
\beno
\bthetad(\bz)
\hide{
\= \argmax\limits_{\btheta\, \in\, \bTheta} \mbE\, \log p_{\bta(\btheta)}(\bX)
}
\= \argmax\limits_{\btheta\, \in\, \bTheta_0} \left[\ell(\btheta, \bz; \bmu^\star) - \ell(\bthetas, \bzs; \bmu^\star)\right],
\ee
which is equivalent to minimizing the Kullback-Leibler divergence\linebreak 
$KL(\bthetas, \bzs;\, \btheta, \bz) = \ell(\bthetas, \zs; \bmu^\star) - \ell(\btheta, \bz; \bmu^\star)$ 
with respect to $\btheta$ given $\bz \in \mZ_0$.
In other words,
$\bthetah(\bz)$ is an estimator of the parameter vector $\bthetad(\bz)$ that is as close as possible to the data-generating parameter vector $\bthetas$ in terms of Kullback-Leibler divergence given $\bz \in \mZ_0$.
These considerations suggest that $\bthetah(\bz)$ may be a consistent estimator of $\bthetad(\bz)$, 
but in general $\bthetah(\bz)$ is not a consistent estimator of $\bthetas$ unless $\bz = \bzs$ \citep{ScSt16}.
\alert{
\com {\em Sharpness.}
\label{com:sharpness}
The results in Proposition \ref{p.z.1} and Theorem \ref{theorem.step3} are not,
and cannot be expected to be as sharp as results based on stochastic block models \citep[e.g.,][]{BiCh09,ChWoAi12,CeDaLa11,RoChYu11,BiChLe11,zhao2012,PrSuTaVo12,AmChBiLe13,MoNeSl15,LeRi13,RoQiFa13,Gao15,Jin15,ZhZh16,BiVoRo17},
for at least three reasons:
\bi
\item {\em Dependence.}
We are concerned with random graphs with dependent edges within blocks,
and concentration results for dependent random variables tend to be weaker than concentration results for independent random variables.
\item {\em The results cover many models with many possible forms of dependence.}
One of the greatest advantages of exponential-family models of random graphs---which can be viewed as generalizations of Erd\H{o}s and R\'enyi random graphs,
GLMs,
and Markov random fields for dependent network data---is the flexibility of the exponential-family framework and its ability to model many dependencies within blocks.
As a consequence,
we do not focus on sharp results in special cases,
but on results that cover many models with many possible forms of dependence.
Indeed,
our concentration results are worst-case results and therefore are not,
and cannot be expected to be sharp in special cases.
\item {\em The combination of dependence and sparsity.}
Many papers concerned with stochastic block models focus on sparse random graphs for which the expected number of edges grows slower than the number of possible edges $\binom{n}{2}$.
While studying random graphs under sparsity assumption makes sense and has a long tradition in classic random graph theory \citep[e.g.,][]{AlSp92,MoRe02,FrKa16},
it requires sharp concentration results for sparse random graphs.
Such results are available for sparse random graphs with independent edges based on, e.g., clever applications of Bernstein's and Talagrand's concentration inequalities \citep{AlSp92,MoRe02,FrKa16}:
e.g.,
\citet{ChWoAi12} used Bernstein's concentration inequality to obtain concentration results for sparse random graphs with independent edges and the expected number of edges growing faster than $n\, (\log n)^{3 + \beta}$ ($\beta > 0$).
But Bernstein's and Talagrand's concentration inequalities are limited to random graphs with independent edges.
To the best of our knowledge,
no sharp concentration results have been developed for sparse random graphs with dependence among edges induced by transitivity or other network phenomena.
While developing sharp concentration results for sparse random graphs with dependent edges would doubtless be an important contribution to the literature,
it is beyond the scope of our paper.
\ei
In short,
the sharpest results can be obtained when edges within and between blocks are independent \citep[e.g.,][]{BiCh09,ChWoAi12,CeDaLa11,RoChYu11,BiChLe11,zhao2012,PrSuTaVo12,AmChBiLe13,MoNeSl15,LeRi13,RoQiFa13,Gao15,Jin15,ZhZh16,BiVoRo17},
but those results come at a cost:
the assumption that edges are independent within and between blocks may be violated in applications,
because network data are dependent data \citep[e.g.,][]{HpLs76,WsFk94,ergm.book}.
We remove the assumption that edges are independent within blocks.
It comes at the cost of less sharp results,
but the benefit is that exponential families with local dependence can capture many features of random graphs that induce dependence among edges within blocks,
including---but not limited to---transitivity,
as explained in Section \ref{sec:comparison}.

}

\section{Simulation results}
\label{sec:simulations}

To demonstrate that the block structure can be recovered in practice,
we simulate data from exponential families with block-dependent edge and transitive edge terms as described in Section \ref{class2}.
\alert{To estimate the block structure,
note that (restricted) maximum likelihood estimators are intractable,
because maximization over (as many as) $\exp(n \log K)$ possible partitions of a set of $n$ nodes into $K$ blocks is infeasible unless $n$ is small.
The same issue arises in stochastic block models,
despite the simplifying assumption that edges are independent conditional on the block structure:
see,
e.g.,
\citet{ChWoAi12} and \citet{RoQiFa13}.
Both of these papers are concerned with theoretical results for (restricted) maximum likelihood estimators, 
but base simulation results on approximate methods,
because (restricted) maximum likelihood estimators are intractable:
\citet{ChWoAi12} use Markov chain Monte Carlo methods,
whereas \citet{RoQiFa13} use pseudolikelihood methods.
We likewise have to resort to approximate methods,
and use Bayesian auxiliary-variable methods for exponential families with local dependence \citep{ScHa13},
as implemented in {\tt R} package {\tt hergm} \citep{ScLu15}.
}

\hide{
Remark: To ensure that the expected number of edges of all nodes in the smallest block---with the most between-block edges---is 3,
we take
\beno
(n - \min(\mA_k,\, \mA_l))\, p \lte 3\s
\\
\dfrac{n - \min(\mA_k,\, \mA_l)}{1 + \exp(-\eta)} \lte 3\s
\\
\dfrac{n - \min(\mA_k,\, \mA_l)}{3} \lte 1 + \exp(-\eta)\s
\\
\dfrac{n - \min(\mA_k,\, \mA_l)}{3} - 1 \lte \exp(-\eta)\s
\\
\log \left(\dfrac{n - \min(\mA_k,\, \mA_l)}{3} - 1\right) \lte -\eta\s
\\
- \log \left(\dfrac{n - \min(\mA_k,\, \mA_l)}{3} - 1\right) \gte \eta\s
\\
\eta \lte - \log \left(\dfrac{n - \min(\mA_k,\, \mA_l)}{3} - 1\right)
\ee
}

\hide{
\subsection{Small networks}
\label{sec:small.networks}

For small networks ($n \leq 100$),
Bayesian auxiliary-variable Markov chain Monte Carlo methods can be used to recover the block structure \citep{ScHa13,ScLu15}.
}
We consider networks with $n = 50$, $n = 75$, and $n = 100$ nodes and $K = 5$ blocks $\mA_1, \dots, \mA_K$ of equal size.
The data-generating natural parameters are given by 
\beno
\eta_{1,k,l} \= - \log \left(\dfrac{n - \min(\mA_k,\, \mA_l)}{3} - 1\right), & k < l = 1, \dots, K,\s
\\
\eta_{1,k,k} \= -1,\hspace{1.5cm} \eta_{2,k,k} \;=\; 1, & k = 1, \dots, K,
\ee
where the between-block natural parameters $\eta_{1,k,l}$ have been chosen to ensure that, for each node, the expected number of edges between blocks is $3$.
To deal with the so-called label-switching problem of Bayesian Markov chain Monte Carlo methods---which arises from the invariance of the likelihood function to the labeling of blocks---we follow the Bayesian decision-theoretic approach of \citep{St00} and estimate block memberships by assigning each node to its maximum-posterior-probability block \citep[][]{ScHa13,ScLu15}.

\begin{figure}
\begin{center}
\includegraphics[scale=0.4]{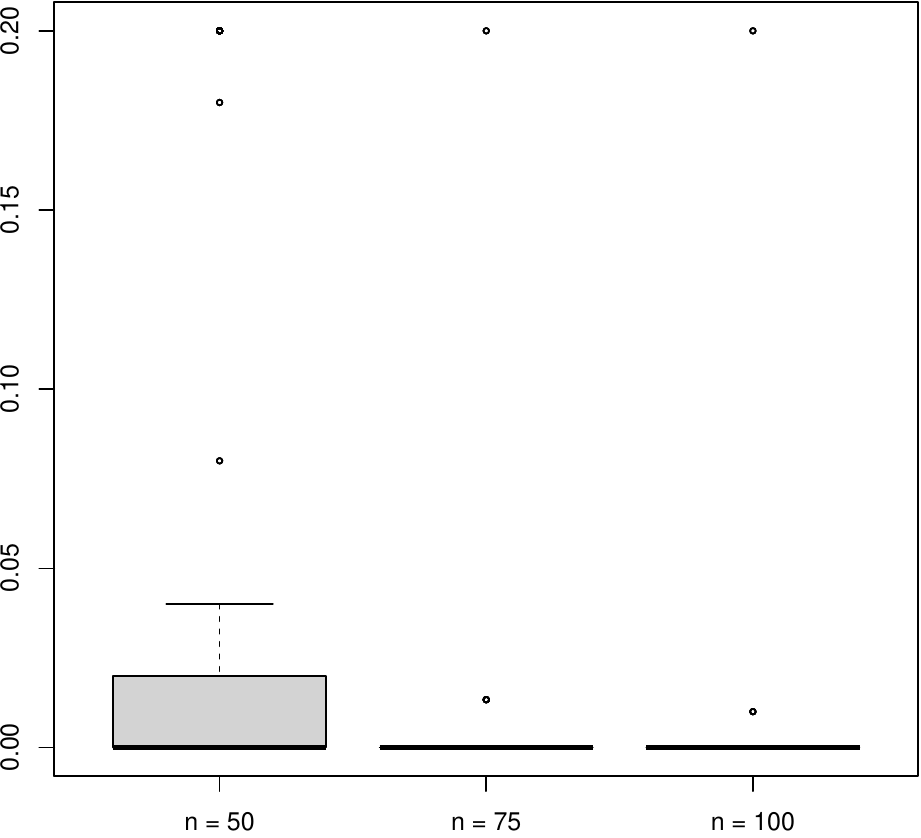}
\end{center}
\caption{\label{figurefigure}Fraction of misclassified nodes based on 100 simulated data sets with $n = 50$, $n = 75$, and $n = 100$ nodes and $K = 5$ blocks of equal size,
where the model is estimated by Bayesian methods.}
\end{figure}

Figure \ref{figurefigure} shows the fraction of misclassified nodes in terms of the normalized minimum Hamming distance $\delta(\zs, \widehat\bz)\, /\, n = \min_{\pi} \sum_{i=1}^n \one_{z_i^\star \neq \pi(\hat{z}_i)}\, /\, n$ based on 100 simulated data sets with $n = 50$, $n = 75$, and $n = 100$ nodes and $K = 5$ blocks of equal size;
we note that Bayesian methods are too time-consuming to be applied to more than 100 simulated data sets.
Figure \ref{figurefigure} suggests that the fraction of misclassified nodes is small in most data sets and decreases as the number of nodes increases from $n = 50$ to $n = 100$ and hence the sizes of the blocks increases from 10 to 20.

\hide{
For each simulated data set,
we used Bayesian auxiliary-variable Markov chain Monte Carlo methods \citep{ScHa13,ScLu15} to recover block memberships of random graphs.
generating a Markov chain Monte Carlo sample of size 100,000 from the posterior,
discarding the first 20,000 sample points as burn-in,
and keeping every $10$-th post-burn-in sample point,
which gives a sample of size 8,000.
}

\hide{
\subsection{Large networks}
\label{sec:large.networks}

For large networks ($n \gg 100$),
Bayesian methods are too time-consuming and approximate methods have to be used.
We demonstrate them here and elaborate on them elsewhere.
The approximate methods are based on the idea that as long as the blocks are not too large and the random graph is not too sparse,
the ${K \choose 2}$ between-block subgraphs dominate the random graph.
Therefore,
despite the fact that exponential families induce dependence within blocks and hence have an added value relative to stochastic block models---as pointed out in Section \ref{sec:comparison}---most of the random graph is governed by the same probability law as random graphs governed by stochastic block models.
Thus,
one can estimate the block structure by using approximate methods based on stochastic block models.
Stochastic block models admit the estimation of block structure from large networks \citep[e.g.,][]{RoChYu11,BiChChZh13,VuHuSc12}.
To demonstrate,
we consider approximate methods based on the following two-step estimation approach.
In the first step,
we estimate the block structure by assuming that $\eta_{2,k,k} = 0$ ($k = 1, \dots, K$)---in which case the exponential family with local dependence reduces to stochastic block models---and estimating the block structure by using variational methods for stochastic block models described in \citep{VuHuSc12}.
In the second step,
we estimate parameters under the assumption that $\eta_{2,k,k} \neq 0$ ($k = 1, \dots, K$) conditional on the estimated block structure by using Monte Carlo maximum likelihood methods described by \citep{HuHa04}.

\begin{figure}
\begin{center}
\includegraphics[scale=0.4]{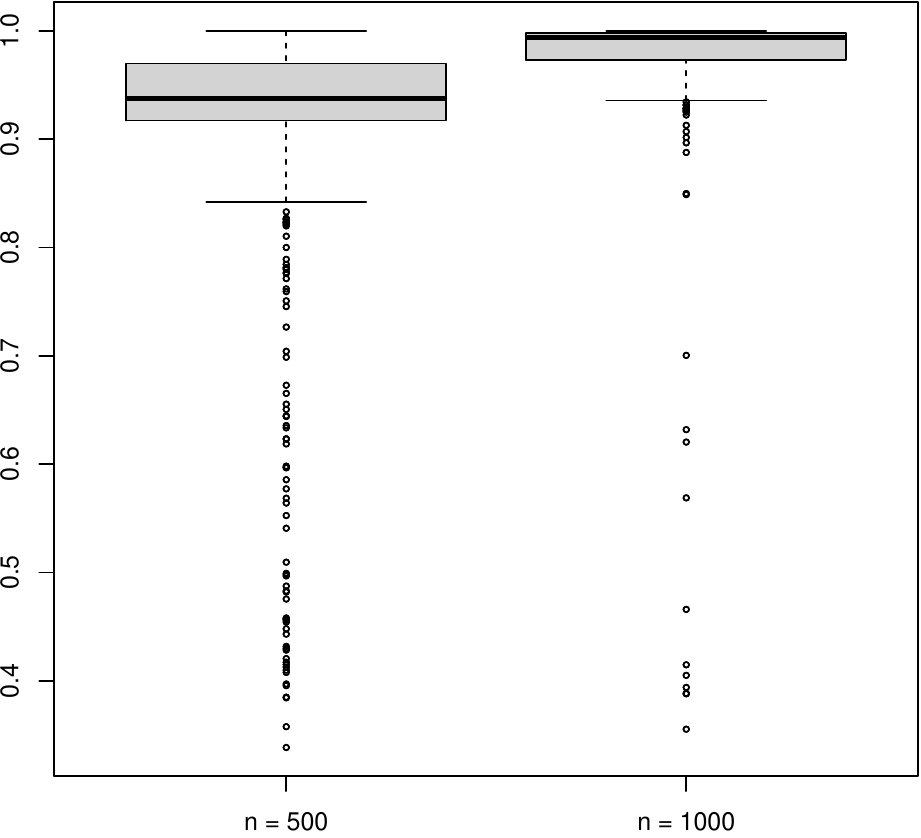}
\end{center}
\caption{\label{figurefigure2}Agreement of estimated and data-generating block structure in terms of Yule's $\phi$-coefficient based on 1,000 simulated data sets with $n = 500$ and $n = $ 1,000 nodes with $K = 50$ blocks of equal size,
where the model is estimated by approximate maximum likelihood methods.}
\end{figure}

We consider networks with $K = 50$ blocks of equal size,
where the size is either 10 or 20, 
so that $n = $ 500 or $n = $ 1,000.
The data-generating natural parameters are given by
\beno
\eta_{1,k,l} \;=\; - \log \left(\dfrac{n - \min(\mA_k,\, \mA_l)}{2.5} - 1\right), & k < l = 1, \dots, K,\s
\\
\eta_{1,k,k} \;=\; - \log \left(\dfrac{n - \min(\mA_k,\, \mA_l)}{5} - 1\right), & k = 1, \dots, K,\s
\\
\eta_{2,k,k} = 1, & k = 1, \dots, K.
\ee

\alert{Figure \ref{figurefigure2} shows the agreement of the estimated and data-generating block structure in terms of Yule's $\phi$-coefficient \citep{phi} based on 1,000 simulated data sets with $n = 500$ and $n = $ 1,000 nodes with $K = 50$ blocks of equal size.
Note that the minimum normalized Hamming distance cannot be computed when $K \gg 5$,
because the minimization over all possible $K!$ permutations is infeasible when $K \gg 5$.
We therefore use Yule's $\phi$-coefficient,
defined by
\beno
\phi(\zs, \hat\bz)
&=& \dfrac{n_{0,0}\, n_{1,1} - n_{0,1}\, n_{1,0}}{\sqrt{(n_{0,0} + n_{0,1})\, (n_{1,0} + n_{1,1})\, (n_{0,0} + n_{1,0})\, (n_{0,1} + n_{1,1})}},
\ee
where
\beno
n_{a,b}
\;=\; n_{a,b}(\zs, \hat\bz)
\;=\; \dsum_{i<j}^{n} \one(\one(\bz_i^\star = \bz_j^\star)=a)\; \one(\one(\hat\bz_i=\hat\bz_j)=b),
& a, b \in \{0, 1\}.
\ee
By definition,
Yule's $\phi$-coefficient is invariant to the labeling of blocks,
so that its calculation does not require computations involving all $K!$ permutations.
}
Figure \ref{figurefigure2} demonstrates that the agreement is high in most data sets and increases as the number of nodes increases from $n = 500$ to $n = $ 1,000 and hence the sizes of the blocks increases from 10 to 20.
}

\section{Proofs of main consistency results}
\label{sec:proofs}

We prove the main consistency results,
Proposition \ref{p.z.1} and Theorem \ref{theorem.step3}.
To prove them,
we need two additional lemmas,
Lemmas \ref{good.event} and \ref{mle.exists}.
The proofs of Lemmas \ref{proposition.concentration},
\ref{good.event}, 
and \ref{mle.exists}
are delegated to the supplementary materials along with the proofs of Corollaries \ref{c.canonical} and \ref{c.curved}.

To state Lemmas \ref{good.event} and \ref{mle.exists},
note that the data-generating natural parameter vector $\bta^\star \in \etaspace \subseteq \mbox{int}(\fullspace)$ is in the interior $\mbox{int}(\fullspace)$ of the natural parameter space $\fullspace$.
Therefore,
the expectation $\mbE\, s(\bX)$ exists \citep[][Theorem 2.2, pp.\ 34--35]{Br86} and so does the expectation $\mbE\, \ell(\btheta, \bz; s(\bX)) = \ell(\btheta, \bz; \mbE\, s(\bX))$.
\hide{
because
\beno
\mbE\; \ell(\btheta, \bz; s(\bX)) 
\hide{
\= \langle\bta(\btheta, \bz),\, \mbE\, s(\bX)\rangle - \psi(\bta(\btheta, \bz))
}
\= \langle\bta(\btheta, \bz),\, \bmu^\star\rangle - \psi(\bta(\btheta, \bz))
\= \ell(\btheta, \bz; \bmu^\star),
\ee
where $\bmu^\star \equiv \bmu(\bta^\star) = \mbE_{\bta^\star}\, s(\bX)$.
}
Let 
\beno
\mG
\= \left\{\bx \in \mX:\; |\ell(\bthetas, \bzs; s(\bx)) -  \ell(\bthetas, \bzs; \bmu^\star)|\; <\; \alpha\; |\ell(\bthetas, \bzs; \bmu^\star)|\right\}
\ee
be the subset of $\bx \in \mX$ such that $s(\bx) \in \mMs$,
where $\alpha > 0$ is identical to the constant $\alpha$ used in the construction of the subset $\mMs$ of the mean-value parameter space $\mM$.
\hide{
Throughout,
to save space,
we write
\beno
\uuushort 
\= \uuu.
\ee
}

Lemma \ref{good.event} shows that the event $\bX \in \mG$ occurs with high probability provided that the number of nodes $n$ is sufficiently large and hence all probability statements in Proposition \ref{p.z.1} and Theorem \ref{theorem.step3} can be restricted to the high-probability subset $\mG$ of $\mX$.

\begin{lemma}
\label{good.event}
Suppose that {an observation} of a random graph is generated by an exponential family with local dependence and countable support $\mX$ satisfying condition \cthree{ }along with assumption \eqref{a.condition}.
Then there exist $C > 0$ and $n_0 > 0$ such that,
for all $n > n_0$,
\beno
\mbP\left(\bX \in \mG\right)
\gte 1 - 2\, \exp\left(- \alpha^2\; C\, n \log n\right),
\ee
where $\alpha > 0$ is identical to the constant $\alpha$ used in the construction of the subset $\mMs$ of the mean-value parameter space $\mM$.
\end{lemma}

\s

Lemma \ref{mle.exists} shows that in the event $\bX \in \mG$,
the restricted maximum likelihood estimator $\estimate$ exists,
which implies that the restricted maximum likelihood estimator $\estimate$ exists with high probability provided that the number of nodes $n$ is sufficiently large by Lemma \ref{good.event}.

\begin{lemma}
\label{mle.exists}
Suppose that {an observation} of a random graph is generated by an exponential family with local dependence and countable support $\mX$ satisfying conditions \cone{ }and \cthree{ }along with assumption \eqref{a.condition}.
Then the following statements hold:
\bi
\item[(a)] For all $\bx \in \mG$,
the restricted maximum likelihood estimator $\estimate$ exists;
\item[(b)] There exist $C > 0$ and $n_0 > 0$ such that,
for all $n > n_0$,
the restricted maximum likelihood estimator $\estimate$ exists with at least probability $1 - 2\, \exp\left(- \alpha^2\; C\, n \log n\right)$;
\ei
where $\alpha > 0$ is identical to the constant $\alpha$ used in the construction of the subset $\mMs$ of the mean-value parameter space $\mM$.
\end{lemma}

\s

Armed with Lemmas \ref{good.event} and \ref{mle.exists},
we can prove Proposition \ref{p.z.1} and Theorem \ref{theorem.step3}.

\pproof \ref{p.z.1}.
Throughout,
to ease the presentation,
we use the short-hand expression
\beno
\uuushort 
\= \uuu.
\ee
By Lemma \ref{good.event},
there exist $C_0 > 0$ and $n_0 > 0$ such that,
for all $n > n_0$,
\beno
\mbP\left(\mX \setminus \mG\right)
\lte 2\, \exp\left(- \alpha^2\; C_0\, n \log n\right).
\ee
Thus,
all following arguments can be restricted to the high-probability subset $\mG$ of $\mX$.
It is therefore convenient to bound the probability of the event $KL(\bthetas, \zs;\; \bthetah, \zh) \geq \epsilon\, \uuushort$ by using a divide- and conquer strategy based on the inequality
\be
\label{ltp}
&& \mbP\left(KL(\bthetas, \zs;\; \bthetah, \zh) \geq \epsilon\, \uuushort\right)\s
\\
\hide{
\= \mbP(\mE \cap \comp\GG) + \mbP(\GG \cap \mE)
\= \mbP(\mE \mid \comp\GG)\, \mbP(\comp\GG) + \mbP(\GG \cap \mE)
}
\lte \mbP\left(KL(\bthetas, \zs;\, \bthetah, \zh) \geq \epsilon\, \uuushort\; \cap\; \mG\right) + \mbP(\mX \setminus \mG).
\ee
The advantage of doing so is that we can confine attention to observations $s(\bx) \in \mMs$ that fall into well-behaved subsets $\mMs$ of the mean-value parameter space $\mM$ satisfying conditions \cone{ }and \ctwo.
Observe that conditions \cone{ }and \ctwo{ }are assumed to hold on $\mMs$, 
but need not hold on $\mM \setminus \mMs$.

To bound the probability of the event $KL(\bthetas, \zs;\; \bthetah, \zh)\, \geq\, \epsilon\, \uuushort \cap \mG$,
note that,
for any $\bx \in \mG$,
the restricted maximum likelihood estimator $\estimate$ exists by Lemma \ref{mle.exists} and that
\beno
KL(\bthetas, \zs;\; \bthetah, \zh)
\= \ell(\bthetas, \bzs; \bmu^\star) - \ell(\bthetah, \zh; \bmu^\star)
\gte 0.
\ee
Since $\estimate \in \bTheta_0 \times \mZ_0$ maximizes $\ell(\btheta, \bz; s(\bx))$ and $\truth \in \bTheta_0 \times \mZ_0$,
we have
\beno
&& \ell(\bthetas, \bzs; \bmu^\star) + [\ell(\bthetas, \zs; s(\bx)) - \ell(\bthetas, \bzs; \bmu^\star)]\s
\\
\lte \ell(\bthetah, \zh; \bmu^\star) + [\ell(\bthetah, \zh; s(\bx)) - \ell(\bthetah, \zh; \bmu^\star)]
\ee
and hence $KL(\bthetas, \zs;\; \bthetah, \zh)$ can be bounded above as follows:
\beno
\ell(\bthetas, \bzs; \bmu^\star) - \ell(\bthetah, \zh; \bmu^\star)
\hide{
\;\leq\; [\ell(\bthetah, \zh; s(\bx)) - \ell(\bthetah, \zh; \bmu^\star)] - [\ell(\bthetas, \zs; s(\bx)) - \ell(\bthetas, \bzs; \bmu^\star)]\s
\\
\hide{
\lte \left|[\ell(\bthetah, \zh; s(\bx)) - \ell(\bthetah, \zh; \bmu^\star)] - [\ell(\bthetas, \zs; s(\bx)) - \ell(\bthetas, \bzs; \bmu^\star)]\right|\s
\\
\lte |\ell(\bthetah, \zh; s(\bx)) - \ell(\bthetah, \zh; \bmu^\star)| + |\ell(\bthetas, \bzs; \bmu^\star) - \ell(\bthetas, \zs; s(\bx))|\s
\\
\=
\lte |\ell(\bthetah, \zh; s(\bx)) - \ell(\bthetah, \zh; \bmu^\star)| + |\ell(\bthetas, \zs; s(\bx)) - \ell(\bthetas, \bzs; \bmu^\star)|\s
\\
}
}
&\leq& 2\, \max\limits_{\bz \in \mZ_0}\, \su_{\btheta\in\bTheta_0} |\ell(\btheta, \bz; s(\bx)) - \ell(\btheta, \bz; \bmu^\star)|.
\ee
\hide{
Since $\bTheta$ is a subset of $\bTheta \subseteq \{\btheta \in \mR^{\dim(\btheta)}: c(\btheta) < \infty\} \subseteq \mR^{\dim(\bta)}$ and thus a subset of $\mR^{\dim(\btheta)}$ and $\mR^{\dim(\bta)}$ is separable,
}
Choose any $\rho > 0$ satisfying $0 < \rho < \epsilon\, /\, (12\, A_1)$,
where $A_1 > 0$ is equal to the constant $A_1 > 0$ in condition \ctwo.
By condition \cfour,
there exist $A, B, C > 0$ such that the $\dim(\btheta) \leq A\, n$-dimensional parameter space $\bTheta_0 \subseteq \bTheta$ can be covered by $\exp(C\, n)$ closed balls with centers $\btheta \in \bTheta$ and radius $B > 0$.
Each of the $\exp(C\, n)$ balls with radius $B > 0$ can be covered by
\beno
\left(\dfrac{4\, B + \rho}{\rho}\right)^{\dim(\btheta)}
\ee
balls $\mB(\btheta, \rho)$ with centers $\btheta \in \bTheta$ and radius $\rho > 0$.
Therefore,
\linebreak
$\bTheta_0 \subseteq \bigcup_{1 \leq q \leq Q}\, \BB(\btheta_q,\, \rho)$ can be covered by $Q$ balls $\mB(\btheta_q, \rho)$ with centers $\btheta_q \in \bTheta$ and radius $\rho > 0$,
where $Q$ is bounded above by
\be
\label{l.bound}
Q
\hide{
\lte N(\bTheta, B)\, \left(\dfrac{4\, B + \rho}{\rho}\right)^{\dim(\btheta)}\s
\\
\lte \exp(C\, n)\, \left(\dfrac{4\, B + \rho}{\rho}\right)^{A\, n}\s
\\
}
\lte \exp\left(A\, \log \left(\dfrac{4\, B + \rho}{\rho}\right)\, n + C\, n\right).
\ee
As a result,
we can write
\beno
\ell(\bthetas, \zs; \bmu^\star) - \ell(\bthetah, \zh; \bmu^\star)
\leq 2 \max\limits_{\bz \in \mZ_0} \max\limits_{1 \leq q \leq Q} \su_{\btheta\, \in\, \BB(\btheta_q,\, \rho)} |\ell(\btheta, \bz; s(\bx)) - \ell(\btheta, \bz; \bmu^\star)|.
\ee
Collecting terms shows that
\beno
\label{uniform.p}
\mbP\left(KL(\bthetas, \zs;\; \bthetah, \zh) \;\geq\; \epsilon\; \uuushort\; \cap\; \mG\right)\s
\\
=\;\; \mbP\left(\ell(\bthetas, \bzs; \bmu^\star) - \ell(\bthetah, \zh; \bmu^\star)\; \geq\; \epsilon\; \uuushort\; \cap\; \mG\right)\s
\\
\leq\;\; \mbP\left(\max\limits_{\bz \in \mZ_0}\, \max\limits_{1 \leq q \leq Q}\, \su_{\btheta\in\BB(\btheta_q,\, \rho)} |\ell(\btheta, \bz; s(\bX)) - \ell(\btheta, \bz; \bmu^\star)| \geq \dfrac{\epsilon\, \uuushort}{2}\; \cap\; \mG\right).
\ee
To bound the probability of the max-sup of deviations of the form $|\ell(\btheta, \bz; s(\bX)) - \ell(\btheta, \bz; \bmu^\star)|$,
observe that,
for any $\bx \in \mG$,
the deviation reduces to
\beno
|\ell(\btheta, \bz; s(\bx)) - \ell(\btheta, \bz; \bmu^\star)|
\hide{
\= \langle\bta(\btheta, \bz),\, s(\bx)\rangle - \psi(\bta(\btheta, \bz))
- \langle\bta(\btheta, \bz),\, \bmu^\star\rangle + \psi(\bta(\btheta, \bz))\s
\\
\= \langle\bta(\btheta, \bz),\, s(\bx)\rangle - \langle\bta(\btheta, \bz),\, \bmu^\star\rangle\s
\\
}
\= |\langle\bta(\btheta, \bz),\, s(\bx)\rangle - \langle\bta(\btheta, \bz),\, \bmu^\star\rangle|,
\ee
because $\psi(\bta(\btheta, \bz))$ cancels.
Consider any $\bz \in \mZ_0$ and any of the $Q$ balls $\mB(\btheta_q,\, \rho)$ that make up the cover $\bigcup_{1 \leq q \leq Q}\, \BB(\btheta_q,\, \rho)$ of $\bTheta_0$.
Let
\beno
\bthetad_q(\bz)
\= \argmax\limits_{\btheta\, \in\, \mbox{\footnotesize cl}\, \mB(\btheta_q,\, \rho)} \ell(\btheta, \bz; \bmu^\star),
\ee
where the subscript $q$ is added to indicate the closed ball $\mbox{cl}\, \mB(\btheta_q,\, \rho)$ that contains $\bthetad_q(\bz)$.
Observe that,
for any $\bz \in \mZ_0$,
$\ell(\btheta, \bz; \bmu^\star)$ is upper semicontinuous on $\mbox{cl}\, \mB(\btheta_q,\, \rho)$ by condition \cone{ }and hence assumes a maximum on $\mbox{cl}\, \mB(\btheta_q,\, \rho)$.
Thus,
for any $\bz \in \mZ_0$,
the maximizer $\bthetad_q(\bz)$ exists and is unique by condition \czero{ }and the assumption that the exponential family is minimal,
which can be assumed without loss \citep[][Theorem 1.9, p.\ 13]{Br86}.
The triangle inequality shows that,
for any $\bx \in \mG$,
any $\bz \in \mZ_0$, 
any $\btheta \in \mbox{cl}\, \mB(\btheta_q,\, \rho)$,
and any $\bthetad_q(\bz) \in \mbox{cl}\, \mB(\btheta_q,\, \rho)$,
\beno
\label{uniform.pp}
|\ell(\btheta, \bz; s(\bx)) - \ell(\btheta, \bz; \bmu^\star)|
\= |\langle\bta(\btheta, \bz),\, s(\bx)\rangle - \langle\bta(\btheta, \bz),\, \bmu^\star\rangle|\s
\\
\lte |\langle\bta(\btheta, \bz),\, s(\bx)\rangle - \langle\bta(\bthetad_q(\bz), \bz),\, s(\bx)\rangle|\s
\\
&+& |\langle\bta(\bthetad_q(\bz), \bz),\, s(\bx)\rangle - \langle\bta(\bthetad_q(\bz), \bz),\, \bmu^\star\rangle|\s
\\
&+& |\langle\bta(\bthetad_q(\bz), \bz),\, \bmu^\star\rangle - \langle\bta(\btheta, \bz),\, \bmu^\star\rangle|.
\ee
\hide{
\beno
\label{uniform.pp}
&& |\ell(\btheta, \bz; s(\bx)) - \ell(\btheta, \bz; \bmu^\star)|
\;=\; |\langle\bta(\btheta, \bz),\, s(\bx)\rangle - \langle\bta(\btheta, \bz),\, \bmu^\star\rangle|\s
\\
\lte |\langle\bta(\btheta, \bz),\, s(\bx)\rangle - \langle\bta(\bthetad_q(\bz), \bz),\, s(\bx) + \langle\bta(\bthetad_q(\bz), \bz),\, \bmu^\star\rangle - \langle\bta(\btheta, \bz),\, \bmu^\star\rangle\rangle|\s
\\
&+& |\langle\bta(\bthetad_q(\bz), \bz),\, s(\bx)\rangle - \langle\bta(\bthetad_q(\bz), \bz),\, \bmu^\star\rangle|\s
\\
\= |\langle\bta(\btheta, \bz) - \bta(\bthetad_q(\bz), \bz),\, s(\bx) - \bmu^\star\rangle| + |\langle\bta(\bthetad_q(\bz), \bz),\, s(\bx) - \bmu^\star\rangle|.
\ee
}
A union bound over the three terms on the right-hand side of the inequality above shows that
\beno
\label{uniform.p}
\mbP\left(\max\limits_{\bz \in \mZ_0}\, \max\limits_{1 \leq q \leq Q}\, \su_{\btheta\in\BB(\btheta_q,\, \rho)} |\ell(\btheta, \bz; s(\bX)) - \ell(\btheta, \bz; \bmu^\star)| \geq \dfrac{\epsilon\; \uuushort}{2} \cap \mG\right)\s
\\
\leq\; \mbP\left(\max\limits_{\bz \in \mZ_0}\, \max\limits_{1 \leq q \leq Q}\, \su_{\btheta\, \in\, \BB(\btheta_q,\, \rho)} |\langle\bta(\btheta, \bz) - \bta(\bthetad_q(\bz),\, \bz),\, s(\bX)\rangle| \geq \dfrac{\epsilon\; \uuushort}{6} \cap \mG\right)\s
\\
+\; \mbP\left(\max\limits_{\bz \in \mZ_0}\, \max\limits_{1 \leq q \leq Q}\, \su_{\btheta\, \in\, \BB(\btheta_q,\, \rho)} |\langle\bta(\bthetad_q(\bz),\, \bz),\, s(\bX) - \bmu^\star\rangle| \geq \dfrac{\epsilon\; \uuushort}{6} \cap \mG\right)\s
\\
+\; \mbP\left(\max\limits_{\bz \in \mZ_0}\, \max\limits_{1 \leq q \leq Q}\, \su_{\btheta\, \in\, \BB(\btheta_q,\, \rho)}  |\langle\bta(\bthetad_q(\bz),\, \bz) - \bta(\btheta, \bz),\, \bmu^\star\rangle| \geq \dfrac{\epsilon\; \uuushort}{6} \cap \mG\right).
\ee
We bound the last three terms on the right-hand side of the inquality above one by one.

\s

{\bf First term.}
The first term can be bounded by using condition \ctwo,
which implies that there exist $A_1 > 0$ and $n_1 > 0$ such that,
for any $n > n_1$,
any $\bx \in \mG$,
any $\bz \in \mZ_0$,  
any $\btheta \in \mbox{cl}\, \mB(\btheta_q,\, \rho)$,
and any $\bthetad_q(\bz) \in \mbox{cl}\, \mB(\btheta_q,\, \rho)$,
\beno
|\langle\bta(\btheta, \bz) - \bta(\bthetad_q(\bz), \bz),\, s(\bx)\rangle|
\lte A_1\, \norm{\btheta - \bthetad_q(\bz)}_2\; \uuushort.
\ee
\hide{
Remark:
Please note that $\btheta_q$ is the center of the ball $\mB(\btheta_q,\, \rho)$ containing both $\btheta$ and $\bthetad_q(\bz)$. 
Thus,
by the triangle inequality,
\beno
\norm{\btheta - \bthetad_q(\bz)}_2
\lte \norm{\btheta - \btheta_q}_2 + \norm{\btheta_q - \bthetad_q(\bz)}_2
\lte \rho + \rho
\= 2\, \rho.
\ee
}
Since both $\btheta$ and $\bthetad_q(\bz)$ are contained in the ball $\mbox{cl}\, \mB(\btheta_q,\, \rho)$,
an application of the triangle inequality shows that
\beno
A_1\, \norm{\btheta - \bthetad_q(\bz)}_2\; \uuushort
\,\leq\, A_1\; 2\; \rho\; \uuushort
\,<\, \dfrac{\epsilon\; \uuushort}{6},
\ee
where we used the fact that $\rho > 0$ satisfies $0 < \rho < \epsilon\, /\, (12\, A_1)$ by construction. 
As a result,
for all $n > n_1$,
we have
\beno
\mbP\left(\max\limits_{\bz \in \mZ_0} \max\limits_{1 \leq q \leq Q} \su_{\btheta\in\BB(\btheta_q,\, \rho)} |\langle\bta(\btheta, \bz) - \bta(\bthetad_q(\bz),\, \bz),\, s(\bX)\rangle| \geq \dfrac{\epsilon\, \uuushort}{6} \cap \mG\right)
= 0.
\ee

\s

{\bf Second term.}
We are interested in bounding the probability of deviations of the form $|\langle\bta(\bthetad_q(\bz), \bz),\, s(\bX) - \bmu^\star\rangle|$.
We make two observations.
First,
observe that,
for any $\bx \in \mG$,
\beno
&& \max\limits_{\bz \in \mZ_0}\, \max\limits_{1 \leq q \leq Q}\, \su_{\btheta\in\mB(\btheta_q,\, \rho)} |\langle\bta(\bthetad_q(\bz), \bz),\, s(\bx) - \bmu^\star\rangle|\s
\\
\= \max\limits_{\bz \in \mZ_0}\, \max\limits_{1 \leq q \leq Q} |\langle\bta(\bthetad_q(\bz), \bz),\, s(\bx) - \bmu^\star\rangle|,
\ee
which implies that
\beno
\mbP\left(\max\limits_{\bz \in \mZ_0}\, \max\limits_{1 \leq q \leq Q}\, \su_{\btheta\in\mB(\btheta_q,\, \rho)} |\langle\bta(\bthetad_q(\bz), \bz),\, s(\bX) - \bmu^\star\rangle|\; \geq\; \dfrac{\epsilon\; \uuushort}{6}\; \cap\; \mG\right)\s
\\
=\; \mbP\left(\max\limits_{\bz \in \mZ_0}\, \max\limits_{1 \leq q \leq Q} |\langle\bta(\bthetad_q(\bz), \bz),\, s(\bX) - \bmu^\star\rangle|\; \geq\; \dfrac{\epsilon\; \uuushort}{6}\; \cap\; \mG\right).
\ee
Second,
bounding the probability of deviations of the form $|\langle\bta(\bthetad_q(\bz), \bz),\, s(\bX) - \bmu^\star\rangle|$ is equivalent to bounding the probability of deviations of the form $|f(\bX) - \mbE\, f(\bX)|$,
where
\beno
f(\bX)
\= \langle\bta(\bthetad_q(\bz), \bz),\, s(\bX)\rangle,
&&
\mbE\, f(\bX)
\= \langle\bta(\bthetad_q(\bz), \bz),\, \bmu^\star\rangle.
\ee
Here,
$f: \mX \mapsto \mR$ is considered as a function of $\bX$ for fixed $(\bthetad_q(\bz), \bz) \in \bTheta_0 \times \mZ_0$.
To bound the probability of deviations of the form $|f(\bX) - \mbE\, f(\bX)|$,
observe that by condition \cthree{ }there exist $A_2 > 0$ and $n_2 > 0$ such that,
for all $n > n_2$, 
the Lipschitz coefficient of $f(\bX)$ satisfies $\norm{f}_{\lip} \leq A_2\, \size$.
Thus,
by applying Lemma \ref{proposition.concentration} to deviations of size $t = \epsilon\; \uuushort \,/\, 6$ along with a union bound over the $|\mZ_0|$ block structures and all $Q$ balls that make up the cover $\bigcup_{1 \leq q \leq Q}\, \BB(\btheta_q,\, \rho)$ of $\bTheta_0$,
there exists $C_1 > 0$ such that,
for all $\epsilon > 0$ and all $n > n_2$,
\beno
&& \mbP\left(\max\limits_{\bz \in \mZ_0}\, \max\limits_{1 \leq q \leq Q} |\langle\bta(\bthetad_q(\bz), \bz),\, s(\bX) - \bmu^\star\rangle|\; \geq\; \dfrac{\epsilon\; \uuushort}{6}\; \cap\; \mG\right)\s
\\
\lte \mbP\left(\max\limits_{\bz \in \mZ_0}\, \max\limits_{1 \leq q \leq Q} |\langle\bta(\bthetad_q(\bz), \bz),\, s(\bX) - \bmu^\star\rangle|\; \geq\; \dfrac{\epsilon\; \uuushort}{6}\right)\s
\\
\hide{
\lte \mbP\left(\bigcup\limits_{\b\bz\, \in\, \mZ}\, \bigcup\limits_{l=1}^L 2\, |\langle\bta(\bthetas, \bz),\, s(\bX) - \bmu^\star\rangle| \;>\; \dfrac{\epsilon\; \uuushort}{2}\; \cap\; \mG\right)\s
\\
\lte \dsum_{\b\bz\, \in\, \mZ} \dsum_{l=1}^L \mbP\left(2\, |\langle\bta(\bthetas, \bz),\, s(\bX) - \bmu^\star\rangle| \;>\; \dfrac{\epsilon\; \uuushort}{2}\; \cap\; \mG\right)\s
\\
\lte \dsum_{\b\bz\, \in\, \mZ} \dsum_{l=1}^L 2\, \exp\left(- \dfrac{\epsilon^2\; \uuushort^2}{C_2\, n^2\, \norm{\mA}_\infty^4\, (\log n)^2}\right)\s
\\
}
\lte 2\, \exp\left(- \dfrac{\epsilon^2\; \uuushort^2}{36\; C_1\, n^2\, \norm{\mA}_\infty^4\, \size^2} + \log |\mZ_0| + \log Q\right).
\ee
To bound the exponential term,
observe that by assumption \eqref{a.condition} of Proposition \ref{p.z.1} there exists,
for all $M > 0$,
however large,
$n_3 > 0$ such that,
for all $n > n_3$,
\beno
\label{a.condition.proofs}
\uuushort 
\gte M\; n^{3/2}\; \norm{\mA}_\infty^2\, \size\, \sqrt{\log n}.
\ee
Therefore,
for all $n > n_3$, 
the three terms in the exponent are bounded above by
\beno
&& - \dfrac{\epsilon^2\, \uuushort^2}{36\; C_1\, n^2\, \norm{\mA}_\infty^4\, \size^2}
+ \log |\mZ_0|
+ \log Q\s
\\
\hide{
\lte - \dfrac{\epsilon^2\, \uuushort^2}{36\; C_1\, n^2\; \size^6}
+ \size\, K \log n
+ A\, \log \left(\dfrac{4\, B + \rho}{\rho}\right)\, K + C\, K\s
\\
}
\lte - \dfrac{\epsilon^2\, \uuushort^2}{36\; C_1\, n^2\, \norm{\mA}_\infty^4\, \size^2}
+ \left[1 + A\, \log \left(\dfrac{4\, B + \rho}{\rho}\right) + C\right]\, n \log n,
\ee
where we used $\log |\mZ_0| \leq n \log K$ and $\log Q \leq (A \log(4\, B + \rho) / \rho + C)\, n$ by \eqref{l.bound}.
Since $M > 0$ can be chosen as large as desired,
we can choose
\beno
M 
&>& \sqrt{36\, C_1\, C_2\, \left[1 + A\, \log \left(\dfrac{4\, B + \rho}{\rho}\right) + C\right]},
\ee
where $C_2 > 0$ is chosen so that $C_2\, \epsilon^2 > 1$.
Hence there exists $C_3 > 0$ such that,
for all $n > n_3$,
\beno
- \dfrac{\epsilon^2\, \uuushort^2}{36\; C_1\, n^2\, \norm{\mA}_\infty^4\, \size^2}
+ \left[1 + A\, \log \left(\dfrac{4\, B + \rho}{\rho}\right) + C\right] n \log n
\leq - \epsilon^2\, C_3\, n \log n.
\ee
Collecting terms shows that,
for all $n > n_3$,
\beno
\label{uniform.p.2a}
\mbP\left(\max\limits_{\bz \in \mZ_0}\, \max\limits_{1 \leq q \leq Q}\, \su_{\btheta\in\mB(\btheta_q,\, \rho)} |\langle\bta(\bthetad_q(\bz), \bz),\, s(\bX) - \bmu^\star\rangle|\; \geq\; \dfrac{\epsilon\; \uuushort}{6}\; \cap\; \mG\right)\s
\\
\leq\; 2\, \exp\left(- \epsilon^2\, C_3\, n \log n\right).
\ee
\hide{
Remark:
\be
\label{a.condition.proofs}
\dfrac{\uuushort^2}{n^2\, \norm{\mA}^6}
&>& \log |\mZ_0|\s
\\
\uuushort^2
&>& n^2\, \norm{\mA}^6 \log |\mZ_0|\s
\\
\uuushort
&>& n\, \size^{7/2} \sqrt{\log |\mZ_0|}\s
\\
\uuushort
&>& n^{1 + \gamma\, /\, 2}\, \size^{7/2} \sqrt{\log |\mZ_0|}
\ee
}
\hide{
Remark:
\beno
\dfrac{\uuushort^2}{n^2\; \norm{\mA}_\infty^{6}}
&>& n \log K\s
\\
\dfrac{\uuushort^2}{n^3\; \log K}
&>& \norm{\mA}_\infty^{6}\s
\\
\norm{\mA}_\infty\s
&<& \left(\dfrac{\uuushort^2}{n^3\; \log K}\right)^{\frac16}
\ee
Choose $\gamma > 0$ such that
\beno
\norm{\mA}_\infty\s
&<& \left(\dfrac{\uuushort^2}{n^{3 + \beta}\; \log K}\right)^{\frac16}
&<& \left(\dfrac{\uuushort^2}{n^3\; \log K}\right)^{\frac16}
\ee
Then
\beno
\dfrac{\uuushort^2}{n^2\; \norm{\mA}_\infty^{6}}
&>& \dfrac{\uuushort^2}{n^2\; \left(\dfrac{\uuushort^2}{n^{3 + \beta}\; \log K}\right)}\s
\\
\= \dfrac{\uuushort^2\, n^{3 + \beta} \log K}{n^2\; \uuushort^2}\s
\\
\= n^{1 + \beta} \log K\s
\\
&>& n \log K
\ee
}

\s

{\bf Third term.}
The third term can be bounded along the same lines as the first term,
which implies that there exists $n_4 > 0$ such that,
for all $n > n_4$,
\beno
\hspace{-.25cm}
\mbP\left(\max\limits_{\bz \in \mZ_0}\, \max\limits_{1 \leq q \leq Q}\, \su_{\btheta\, \in\, \BB(\btheta_q,\, \rho)}  |\langle\bta(\bthetad_q(\bz),\, \bz) - \bta(\btheta, \bz),\, \bmu^\star\rangle| \geq \dfrac{\epsilon\, \uuushort}{6}\; \cap\; \mG\right)
= 0.
\ee

\s

{\bf Conclusion.}
Using \eqref{ltp} and collecting terms shows that there exists $C > 0$ such that,
for all $\epsilon > 0$ and all $n > \max(n_0, n_1, n_2, n_3, n_4)$,
\beno
&& \mbP\left(KL(\bthetas, \zs;\; \bthetah, \zh) \;\geq\; \epsilon\; \uuushort\right)
\;\leq\; 2\, \exp\left(- \alpha^2\, C_0\, n \log n\right)\s
\\
&+& 2\, \exp\left(- \epsilon^2\, C_3\, n \log n\right)
\;\leq\; 4\, \exp\left(- \min(\alpha^2, \epsilon^2)\; C\; n \log n\right).
\ee

\ttproof \ref{theorem.step3}.
By assumption \eqref{id} of Theorem \ref{theorem.step3},
there exist $C_1 > 0$ and $n_1 > 0$ such that,
for all $n > n_1$,
\beno
KL(\bthetas, \zs;\; \bthetah, \zh)
\gte \dfrac{\delta(\zs, \zh)\; C_1\; \uuu}{n}
\ee
provided $\estimate$ exists.
By Proposition \ref{p.z.1},
there exist $C_2 > 0$ and $n_2 > 0$ such that,
for all $\epsilon > 0$ and all $n > n_2$,
the event
\beno
KL(\bthetas, \zs;\; \bthetah, \zh)
&<& \epsilon\; C_1\; \uuu
\ee
occurs with at least probability
\be
\label{leastprob}
1 - 4\, \exp\left(- \min(\alpha^2, \epsilon^2)\; C_2\; n \log n\right).
\ee
Therefore,
for all $\epsilon > 0$ and all $n > \max(n_1, n_2)$,
with at least probability \eqref{leastprob},
we observe the event
\[
\begin{array}{ccccc}
\dfrac{\delta(\zs, \zh)\; C_1\; \uuu}{n}
\,\leq\, KL(\bthetas, \zs;\; \bthetah, \zh)
\,<\, \epsilon\, C_1\; \uuu,
\end{array}
\]
i.e.,
the event $\delta(\zs, \zh)\, /\, n\, <\, \epsilon$.

\section{Discussion}
\label{sec:discussion}

Here,
and elsewhere \citep{ScSt16},
we have taken first steps to demonstrate that---while statistical inference for exponential-family random graph models without additional structure is problematic \citep{Ha03,Sc09b,ChDi11,ShRi11}---statistical inference for exponential-family random graph models with additional structure in the form of block structure makes sense.
It goes without saying that numerous open problems remain,
ranging from probabilistic problems (e.g., understanding properties of probability models) and statistical problems (e.g., understanding properties of statistical methods) to computational problems (e.g., the development of computational methods for large networks).

One important problem is that the maximum likelihood estimator discussed here is at least as intractable as maximum likelihood estimators in the special case of stochastic block models \citep[][]{ChWoAi12,RoQiFa13}.
The intractability stems in part from the fact that the block structure is unknown and the number of possible block structures is large and in part from the fact that the likelihood function is intractable even when the block structure is known owing to complex dependence within blocks.
There do exist Bayesian auxiliary-variable methods for small networks \citep{ScHa13,ScLu15} and promising directions for methods for large networks \citep{Waetal18,BaSc17}.
\alert{As pointed out in the introduction,
an indepth investigation of all of these models and methods is beyond the scope of a single paper.}
However,
the main consistency results reported here suggest that statistical inference for these models and methods is possible and worth exploring in more depth.

\section*{Acknowledgements}

The author acknowledges support from the National Science Foundation (NSF awards DMS-1513644 and DMS-1812119).

\section*{Supplementary materials}

The proofs of Lemmas \ref{proposition.concentration}, \ref{good.event}, and \ref{mle.exists} and Corollaries \ref{c.canonical} and \ref{c.curved} can be found in the supplementary materials.

\hide{

\end{supplement}

}



\hide{
\begin{supplement} 
\stitle{\longtitle}
\slink[doi]{}
\sdatatype{.pdf}
\sdescription{We prove all results in \citep{Sc15}.}
}

\pagebreak

\makeatletter

\setcounter{page}{1}

\setcounter{section}{0}

\begin{frontmatter}

\title{Supplement:\\
\longtitle}
\runtitle{
\shorttitle}

\begin{aug}
\author{\fnms{Michael} \snm{Schweinberger}\ead[label=e1]{m.s@rice.edu}}
\affiliation{Rice University}
\address{
Michael Schweinberger\\
Department of Statistics\\
Rice University\\
6100 Main St, MS-138\\
Houston, TX 77005-1827\\
E-mail:\ m.s@rice.edu
}
\end{aug}

\end{frontmatter}

\begin{appendix}

\section{Proofs of auxiliary results}
\label{sec:auxiliary.proofs}

We prove Lemmas \ref{proposition.concentration}, \ref{good.event}, and \ref{mle.exists} and Corollaries \ref{c.canonical} and \ref{c.curved}.

\s

\llproof \ref{proposition.concentration}.
By assumption,
$\mbE\, f(\bX) < \infty$.
We are interested in deviations of the form $|f(\bX) - \mbE\, f(\bX)| \geq t$,
where $t > 0$.
In the following,
we denote by $\mbP$ a probability measure on $(\mX, \mS)$ with densities of the form \eqref{example1},
where $\mS$ is the power set of the countable set $\mX$.
Let $\bX = (\bX_{k,l})_{k \leq l}^K$ be a sequence of edge variables,
where $\bX_{k,k} = (X_{i,j})_{i\in\mA_k\, <\, j\in\mA_k}$ denotes the sequence of within-block edge variables of nodes in block $\mA_k$ and $\bX_{k,l} = (X_{i,j})_{i\in\mA_k,\, j\in\mA_l}$ denotes the sequence of between-block edge variables between nodes in blocks $\mA_k$ and $\mA_l$ ($k < l$).
In an abuse of notation,
we denote the elements of the sequence of edge variables $\bX$ by $X_1, \dots, X_m$ with sample spaces $\mX_1, \dots, \mX_m$,
respectively,
where $m = {n \choose 2} \leq n^2$ is the number of edge variables.
Let $\bX_{i:j} = (X_i, \dots, X_j)$ be a subsequence of edge variables with sample space $\mX_{i:j}$,
where $i \leq j$.
By applying Theorem 1.1 of \citep[][]{KoRa08} to $\norm{f}_{\lip}$-Lipschitz functions $f: \mX \mapsto \mbR$ defined on the countable set $\mX$,
\beno
\mbP(|f(\bX) - \mbE\, f(\bX)|\; \geq\; t)
\lte 2\, \exp\left(- \dfrac{t^2}{2\, m\, \norm{\Phi}_\infty^2\, \norm{f}_{\lip}^2}\right),
\ee
where $\Phi$ is the $m \times m$-upper triangular matrix with entries
\beno
\phi_{i,j} \=
\begin{cases}
\varphi_{i,j} & \mbox{if } i < j\\
1 & \mbox{if } i = j\\
0 & \mbox{if } i > j
\end{cases}
\ee
and
\beno
\norm{\Phi}_\infty
\= \max\limits_{1 \leq i \leq m} \left|1 + \dsum_{j=i+1}^m \varphi_{i,j}\right|.
\ee
The coefficients $\varphi_{i,j}$ are known as mixing coefficients and are defined by
\beno
\label{mixing.coefficients}
\varphi_{i,j}
\equiv \sup\limits_{\substack{\bx_{1:i-1} \in \mbX_{1:i-1}\\(\uu,\, \vv) \in \mX_i\times\mX_i}} \varphi_{i,j}(\bx_{1:i-1}, \uu, \vv)
\= \sup\limits_{\substack{\bx_{1:i-1} \in \mbX_{1:i-1}\\(\uu,\, \vv) \in \mX_i\times\mX_i}} \norm{\pi_{\uu} - \pi_{\vv}}_\tv,
\ee
where $\norm{\pi_{\uu} - \pi_{\vv}}_\tv$ is the total variation distance between the distributions $\pi_{\uu}$ and $\pi_{\vv}$
given by
\beno
\pi_{\uu}
\equiv \pi(\bx_{j:m} \mid \bx_{1:i-1}, \uu)
= \mbP(\bX_{j:m} = \bx_{j:m} \mid \bX_{1:i-1} = \bx_{1:i-1}, X_i = \uu)
\ee
and
\beno
\pi_{\vv}
\equiv \pi(\bx_{j:m} \mid \bx_{1:i-1}, \vv)
= \mbP(\bX_{j:m} = \bx_{j:m} \mid \bX_{1:i-1} = \bx_{1:i-1}, X_i = \vv).
\ee
Since the support
of $\pi_{\uu}$ and $\pi_{\vv}$ is countable,
\beno
\norm{\pi_{\uu} - \pi_{\vv}}_\tv
\hide{
\= \dfrac12\, \norm{\pi_{\uu} - \pi_{\vv}}_1
}
\= \dfrac12 \dsum\limits_{\bx_{j:m} \in \mbX_{j:m}} |\pi(\bx_{j:m} \mid \bx_{1:i-1}, \uu) - \pi(\bx_{j:m} \mid \bx_{1:i-1}, \vv)|.
\ee
An upper bound on $\norm{\Phi}_\infty$ can be obtained by bounding the mixing coefficients $\varphi_{i,j}$ as follows.
Consider any pair of edge variables $X_i$ and $X_j$.
If $X_i$ and $X_j$ involve nodes in more than one block,
the mixing coefficient $\varphi_{i,j}$ vanishes by the local dependence induced by exponential families with local dependence.
If the pair of nodes corresponding to $X_i$ and the pair of nodes corresponding to $X_j$ belong to the same block,
the mixing coefficient $\varphi_{i,j}$ can be bounded as follows:
\beno
\varphi_{i,j}(\bx_{1:i-1}, \uu, \vv)
\hide{
\;=\; \norm{\pi - \pi^\prime}_\tv
}
\;=\; \dfrac12 \dsum\limits_{\bx_{j:m} \in \mbX_{j:m}} |\pi(\bx_{j:m} \mid \bx_{1:i-1}, \uu) - \pi(\bx_{j:m} \mid \bx_{1:i-1}, \vv)|\s
\hide{
\\
\lte \dfrac12 \dsum\limits_{x_{j:m} \in \mbX_{j:m}} \left(|\pi(x_{j:m} \mid \bx_{1:i-1}, \uu)| + |\pi(x_{j:m} \mid \bx_{1:i-1}, \vv)|\right)\s
}
\\
\hspace{0cm}\leq\; \dfrac12 \dsum\limits_{\bx_{j:m} \in \mbX_{j:m}} \pi(\bx_{j:m} \mid \bx_{1:i-1}, \uu) + \dfrac12 \dsum\limits_{\bx_{j:m} \in \mbX_{j:m}} \pi(\bx_{j:m} \mid \bx_{1:i-1}, \vv)
\;=\; 1,
\ee
because $\pi_{\uu}$ and $\pi_{\vv}$ are conditional probability mass functions with countable support $\mbX_{j:m}$.
We note that the upper bound is not sharp,
but it has the advantage that it covers a wide range of dependencies within blocks.
As a result,
\beno
\norm{\Phi}_\infty
\= \max\limits_{1 \leq i \leq m} \left|1 + \dsum_{j=i+1}^m \varphi_{i,j}\right|
\lte \dis{\norm{\mA}_\infty \choose 2},
\ee
because each edge variable $X_i$ can depend on at most ${\norm{\mA}_\infty \choose 2}$ edge variables corresponding to pairs of nodes belonging to the same pair of blocks.
Therefore,
there exists $C > 0$ such that,
for all $K > 0$ and all $t > 0$,
\beno
\mbP(|f(\bX) - \mbE\, f(\bX)|\; \geq\; t)
\lte 2\, \exp\left(- \dfrac{t^2}{C\, n^2\, \norm{\mA}_\infty^4\, \norm{f}_{\lip}^2}\right),
\ee
where $\norm{\mA}_\infty > 0$ and $\norm{f}_{\lip} > 0$ by assumption.

\hide{

\begin{proposition}
\label{theorem.step2}
Suppose that $\truth \in \bTheta\times\mG$ and that {an observation} of a random graph is generated by an exponential family with local dependence with additional structure and that the assumptions of Theorem \ref{theorem.step3}.
If $\beta > 1 / 2 + \gamma$,
then there exist $C > 0$ and $n_0 > 0$ such that,
for all $n > n_0$,
\beno
\mbP\left(\estimate \in \bTheta\times\mG\right)
\gte 1 - 2\, \exp\left(- \dfrac{C\, n^{2\, \beta - 2\, \gamma}}{(\log n)^{2\, (1 + \alpha)\, \gamma}}\right).
\ee
\end{proposition}

The result is interesting,
as it shows that the maximum likelihood estimator $\estimate$ is in $\bTheta\times\mG$ with high probability and thus the search of the maximum likelihood estimator may be confined to $\bTheta\times\mG$,
i.e.,
to block structures where none of the blocks is too large.

\pproof \ref{theorem.step2}.
For all $z \in \mZ$,
the profile likelihood maximum likelihood estimator of $\btheta$ given $z$ defined by
\beno
\bthetah(z) \= \argmax\limits_{\btheta\in\bTheta} \ell(\btheta, \bz; s(\bx))
\ee
exists,
because $\bTheta \subset \bTheta = \mR^{\rrr}$ is a compact subset of $\mR^{\rrr}$,
$\ell$ is an upper semicontinuous function on the compact set $\bTheta$,
and therefore $\ell$ assumes a maximum on $\bTheta$ \citep[][Theorem 5.7, p.\ 152]{Br86};
note that $\mZ$ excludes all block structures with empty blocks and that $\zs \in \mZ$.
We want to bound the probability of observing $x \in \mX$ such that the maximum likelihood estimator $\estimate$ falls into $\bTheta \times \mB$.
For the maximum likelihood estimator $\estimate$ to fall into $\bTheta\times\mB$,
there must exist $(\btheta, \bz) \in \bTheta\times\mB$ such that
\beno
\ell(\btheta, \bz; s(\bx)) 
\gte \ell(\bthetah, \zh; s(\bx))
\gte \ell(\bthetas, \zs; s(\bx)), 
\ee
which can be expressed in terms of the expected loglikelihood functions and the deviations of the loglikelihood functions from the expected loglikelihood functions:
\be
\label{inequality0}
\ell(\btheta, \bz; s(\bx)) - \ell(\bthetas, \zs; s(\bx))
\= \ell(\btheta, \bz; \bmu^\star) - \mbE\;\ell(\bthetas, \zs)
+ \Delta(\btheta, \bz;\, \bthetas, \zs)
\gte 0,
\ee
where
\beno
\Delta(\btheta, \bz;\, \bthetas, \zs) 
\= [\ell(\btheta, \bz; s(\bx)) - \ell(\btheta, \bz; \bmu^\star)]
- [\ell(\bthetas, \zs; s(\bx)) - \mbE\;\ell(\bthetas, \zs)].
\ee
By \eqref{inequality0} and Proposition \ref{p.z.2},
\beno
\Delta(\btheta, \bz;\, \bthetas, \zs)
\gte \ell(\bthetas, \bzs; \bmu^\star) - \ell(\btheta, \bz; \bmu^\star)
\gte \dfrac{\delta(\zs, \bz)\, \uuu}{K}.
\ee
To lower bound $\delta(\zs, \bz)$,
observe that $z \in \mB$ and $\zs \in \mGs$ by assumption.
Therefore,
there existst least one block $k$ such that
\beno
\delta(\zs, \bz)
&>& |\mA_k| - |\mA_k|
\gte \ccc - \aaa.
\ee
Since $\aaa = A\, (\log n)^{1 + \alpha}$ and $\ccc = \ccclong$ with $0 < \gamma < 1$,
there exist $C_1 > 0$ and $n_0 > 0$ such that,
for all $n > n_0$,
\beno
\delta(\zs, \bz)
\gte \ccc - \aaa
\gte C_1\, n^{1 - \gamma}.
\ee
Therefore,
the deviations of the loglikelihood functions from the corresponding expectations is bounded below by
\be
\label{inequality1}
\Delta(\btheta, \bz;\, \bthetas, \zs)
\hide{
\gte \dfrac{\delta(\zs, \bz)\, \uuu}{K}
}
\gte \dfrac{C_1\, n^{1-\gamma}\, \uuu}{K}
\hide{
\= \dfrac{C_1\, n^{-\gamma}\, \uuu}{K}
\;\;\geq\;\; C_2\, K^{1-\gamma}\, K^\beta\, (\log n)^\gamma
}
\gte C_2\, K^{1+\beta-\gamma}\, (\log n)^\gamma
\ee
using $\uuu = \uuulong$ and $n \geq K$.
By Jensen's inequality and the compactness of $\bTheta$,
the deviations of the loglikelihood functions from the corresponding expectations are bounded above by
\be
\label{inequality2}
\Delta(\btheta, \bz;\, \bthetas, \zs)
\;\leq\; \left|\Delta(\btheta, \bz;\, \bthetas, \zs)\right|\s
\hide{
\\
\= \left|[\ell(\btheta, \bz; s(\bx)) - \ell(\btheta, \bz; \bmu^\star)] - [\ell(\bthetas, \zs; s(\bx)) - \ell(\btheta, \bz; \bmu^\star)]\right|\s
\\
\lte \left|\ell(\btheta, \bz; s(\bx)) - \ell(\btheta, \bz; \bmu^\star)\right| + \left|\ell(\bthetas, \zs; s(\bx)) - \ell(\btheta, \bz; \bmu^\star)\right|\s
\\
}
\;\leq\; 2\, \max\limits_{\bz\, \in\, \mZ}\, \max\limits_{1 \leq q \leq Q} \su_{\bthetas\in\BB(\bthetas,\, \rho)} \left|\ell(\btheta, \bz; s(\bx)) - \ell(\btheta, \bz; \bmu^\star)\right|.
\ee
Combining the lower bound on $\Delta$ in \eqref{inequality1} and the upper bound on $\Delta$ in \eqref{inequality2} shows that,
for the maximum likelihood estimator $\estimate$ to fall into $\bTheta\times\mB$,
we must have
\beno
2\, \max\limits_{\bz\, \in\, \mZ}\, \max\limits_{1 \leq q \leq Q} \su_{\bthetas\in\BB(\bthetas,\, \rho)} \left|\ell(\btheta, \bz; s(\bx)) - \ell(\btheta, \bz; \bmu^\star)\right|
\gte \Delta(\btheta, \bz;\, \bthetas, \zs)
\gte C_2\, K^{1+\beta-\gamma}\, (\log n)^\gamma.
\ee
\hide{
In other words,
we have
\beno
2\, \max\limits_{\bz\, \in\, \mZ}\, \max\limits_{1 \leq q \leq Q} \su_{\bthetas\in\BB(\bthetas,\, \rho)} \left|\ell(\btheta, \bz; s(\bx)) - \ell(\btheta, \bz; \bmu^\star)\right|
\gte C_2\, K^{1+\beta-\gamma}\, (\log n)^\gamma.
\ee
}
Therefore,
for the maximum likelihood estimator $\estimate$ to fall into $\bTheta\times\mB$,
the deviations of the loglikelihood functions from the corresponding expectations must be $\Omega(K^{1+\beta-\gamma}\, (\log n)^\gamma)$.
We show that,
if $\beta > 1 / 2 + \gamma$,
then deviations of order $\Omega(K^{1+\beta-\gamma}\, (\log n)^\gamma)$ are improbable.
To do so,
observe that,
by \ctwo{ }and the boundedness of $\bTheta$,
there exists $C_3 > 0$ such that the Lipschitz coefficient of $f$ is bounded above by $\norm{f}_{\lip} \leq C_3 \log n$.
Thus,
by applying Lemma \ref{proposition.concentration} to deviations of the form $t = C_2\, K^{1+\beta-\gamma}\, (\log n)^\gamma$ along with a union bound over all $z \in \mZ$ and all $L \leq C \exp(n \log K)$ open balls that make up the finite cover of $\bTheta$,
we have,
for all $\epsilon > 0$,
\beno
&& 
\hide{
\mbP\left(|[\ell(\btheta, \bz; s(\bx)) - \ell(\btheta, \bz; \bmu^\star)] + [\ell(\bthetas, \zs; s(\bX)) - \ell(\btheta, \bz; \bmu^\star)]| \;>\; C_2\, K^{1+\beta-\gamma}\, (\log n)^\gamma\right)\s
\\
\lte 
}
\mbP\left(\max\limits_{\bz\, \in\, \mZ}\, \max\limits_{1 \leq q \leq Q} \su_{\bthetas\in\BB(\bthetas,\, \rho)}|\ell(\btheta, \bz; s(\bX)) - \ell(\btheta, \bz; \bmu^\star)|  \;>\; \dfrac{C_2\, K^{1+\beta-\gamma}\, (\log n)^\gamma}{2}\right)\s
\\
\hide{
\lte \mbP\left(\bigcup\limits_{z\in\mG}\, \bigcup\limits_{l=1}^L \su_{\btheta\in\BB(\bthetas, \rho)}|\langle\bta(\btheta, \bz),\, s(\bX) - \bmu^\star\rangle| \;>\; \dfrac{C_2\, K^{1+\beta-\gamma}\, (\log n)^\gamma}{2}\right)\s
\\
\lte \dsum_{z\in\mG} \dsum_{l=1}^L \mbP\left(\su_{\btheta\in\BB(\bthetas, \rho)}|\langle\bta(\btheta, \bz),\, s(\bX) - \bmu^\star\rangle| \;>\; \dfrac{C_2\, K^{1+\beta-\gamma}\, (\log n)^\gamma}{2}\right)\s
\\
\lte \dsum_{z\in\mG} \dsum_{l=1}^L 2\, \exp\left(- \dfrac{C_4\, K^{2 + 2\, \beta - 2\, \gamma}\, (\log n)^\gamma}{n^2\, \norm{\zs)}_\infty^4\, (\log n)^2}\right)\s
\\
\lte \dsum_{\bz\, \in\, \mZ} \dsum_{l=1}^L 2\, \exp\left(- \dfrac{C_4\, K^{2 + 2\, \beta - 2\, \gamma}\, (\log n)^\gamma}{n^2\, \norm{\zs)}_\infty^4\, (\log n)^2}\right)\s
\\
\= 2\; |\mZ_0|\, \exp\left(- \dfrac{C_5\, K^{2 + 2\, \beta - 2\, \gamma}\, (\log n)^\gamma}{n^2\, \norm{\zs)}_\infty^4\, (\log n)^2} + \log L\right)\s
\\
\= 2\; \exp\left(- \dfrac{C_3\, K^{2 + 2\, \beta - 2\, \gamma}\, (\log n)^\gamma}{n^2\, \norm{\zs)}_\infty^4\, (\log n)^2} + \log |\mZ_0| + \log L\right)\s
\\
\= 2\; \exp\left(- \dfrac{C_3\, K^{2 + 2\, \beta - 2\, \gamma}\, (\log n)^\gamma}{n^2\, \norm{\zs)}_\infty^4\, (\log n)^2} + \log |\mZ_0| + \log L\right)\s
\\
}
\lte 2\, \exp\left(- \dfrac{C_4\, K^{2 + 2\, \beta - 2\, \gamma}\, (\log n)^{2\, \gamma}}{n^2\, \norm{\zs)}_\infty^4\, (\log n)^2} + \log |\mZ_0| + \log L\right)\s
\\
\lte 2\, \exp\left(- \dfrac{C_5\, n^{2 + 2\, \beta - 2\, \gamma}\, (\log n)^{2\, \gamma}}{n^2\, \norm{\zs)}_\infty^4\, (\log n)^2\, (\log n)^{(2 + 2\, \beta - 2\, \gamma)\, (1 + \alpha)}} + \log |\mZ_0| + \log L\right)\s
\\
\lte 2\, \exp\left(- \dfrac{C_6\, n^{2\, \beta - 2\, \gamma}}{(\log n)^{2\, (1 + \alpha)\, \gamma}} + \log |\mZ_0| + \log L\right)
\ee
using $\norm{\zs)}_\infty \leq \aaalong$ since $\zs \in \mGs$ and thus 
\beno
\aaa^4\, (\log n)^2\, (\log n)^{(1 + \alpha)\, (2 + 2\, \beta - 2\, \gamma)}
\hide{
\= A^4\, (\log n)^{2 + 4\, (1 + \alpha) + (1 + \alpha)\, (2 + 2\, \beta - 2\, \gamma)}\s
\\
\= 
A^4\, (\log n)^{2 + (1 + \alpha)\, (6 + 2\, \beta - 2\, \gamma)}
\;\;=\;\; 
}
\= A^4\, (\log n)^{2 + 2\, (1 + \alpha)\, (3 + \beta - \gamma)},
\ee
which,
together with the assumption $\gamma \geq 1 + (1 + \alpha)\, (3 + \beta)$ stated in \cthree,
implies that
\beno
\dfrac{(\log n)^{2\, \gamma}}{\norm{\mA}_\infty^4\, (\log n)^2\, (\log n)^{(1 + \alpha)\, (2 + 2\, \beta - 2\, \gamma)}}
\hide{
\gte 2 + 2\, (1 + \alpha)\, (3 + \beta) - [2 + 2\, (1 + \alpha)\, (3 + \beta - \gamma)]\s
\\
\= 2\, (1 + \alpha)\, (3 + \beta - 3 - \beta - \gamma)\s
\\
}
\gte \dfrac{1}{(\log n)^{2\, (1 + \alpha)\, \gamma}}.
\ee
If $\beta > 1 / 2 + \gamma$,
then $2\, \beta - 2\, \gamma > 1$ and therefore there exists $n_0 > 0$ such that,
for all $n > n_0$,
\beno
&& \mbP\left(\max\limits_{\bz\, \in\, \mZ}\, \max\limits_{1 \leq q \leq Q} \su_{\bthetas\in\BB(\bthetas,\, \rho)}|\ell(\btheta, \bz; s(\bX)) - \ell(\btheta, \bz; \bmu^\star)|  \;>\; \dfrac{C_2\, K^{1+\beta-\gamma}\, (\log n)^\gamma}{2}\right)\s
\\
\lte 2\, \exp\left(- \dfrac{C_6\, n^{2\, \beta - 2\, \gamma}}{(\log n)^{2\, (1 + \alpha)\, \gamma}} + \log |\mZ_0| + \log L\right)\s
\;\;\leq\;\; 2\, \exp\left(- \dfrac{C_7\, n^{2\, \beta - 2\, \gamma}}{(\log n)^{2\, (1 + \alpha)\, \gamma}}\right).
\ee
As a result,
there exist $C > 0$ and $n_0 > 0$ such that,
for all $n > n_0$,
the maximum likelihood estimator $\estimate$ is contained in $\bTheta\times\mG$ with at least probabibility 
\beno
1 - 2\, \exp\left(- \dfrac{C_4\, n^{2\, \beta - 2\, \gamma}}{(\log n)^{2\, (1 + \alpha)\, \gamma}}\right)
\ee
provided $\beta > 1 / 2 + \gamma$.

}

\hide{

\s

We prove Corollary \ref{corollary.edges}.

\ccproof \ref{corollary.edges}.
By assumption,
the data-generating exponential family with local dependence contains indicators of all possible edges $x_{i,j}$ and transitive edges $x_{i,j}\, \max_{k \in \mA_{\zs_i} \cup \mA_{\zs_j},\, k \neq i,j} x_{i,k}\, x_{j,k}$ as sufficient statistics,
the natural parameters corresponding to edges and transitive edges satisfy $\bta_{i,j,1}\truth < 0$ and $\bta_{i,j,2}\truth > 0$,
and the marginal distributions of within-block subgraphs satisfy weak sparsity conditions while between-block subgraphs satisfy strong sparsity assumptions as defined in Appendix \ref{size.dependent.parameterizations}.
Therefore,
\beno
\bta_{i,j,1}\truth 
\= 
\begin{cases}
\bthetas_{\zs_i,\zs_j,1}\, \log n_{\zs_i}(\zs) & \mbox{if } \zs_i = \zs_j\s
\\
\bthetas_{\zs_i,\zs_j,1}\, \log n & \mbox{if } \zs_i \neq \zs_j
\end{cases}
\ee
and
\beno
\bta_{i,j,2}\truth 
\= 
\begin{cases}
\bthetas_{\zs_i,\zs_j,2}\, \dfrac{\log n_{\zs_i}(\zs)}{n_{\zs_i}(\zs)} & \mbox{if } \zs_i = \zs_j\s
\\
\bthetas_{2,\zs_i,\zs_j}\, \dfrac{\log n}{\max(n_{\zs_i}(\zs),\, n_{\zs_j}(\zs))} & \mbox{if } \zs_i \neq \zs_j,
\end{cases}
\ee
where $\max(n_k(\zs),\, n_l(\zs)) > 0$ for all $k \leq l$ since $\zs \in \mGs \subseteq \mZ$ and where we used the fact that a change in one edge can change the number of transitive edges within and between blocks $\mA_k$ and $\mA_l$ in the worst case by $\max(n_k(\zs),\, n_l(\zs)) - 2$ ($k \leq l$).
In addition,
note that $\bthetas_{k,l,1} < 0$ and $\bthetas_{k,l,2} > 0$ by the assumption $\bta_{i,j,1}\truth < 0$ and $\bta_{i,j,2}\truth > 0$ ($k \leq l$).
\hide{
If,
e.g.,
$X_{1,k} = X_{2,k} = 1$ for $k = 3, \dots, n$ while all other edges are absent,
then adding edge $X_{1,2} = 1$ transforms the existing edges $X_{1,k}$ for $k = 3, \dots, n$ from intransitive into transitive edges and thus increases the number of transitive edges from $0$ to $n - 2$.
}
Let $\bta \equiv \bta\truth$ be the natural parameter vector of the exponential family with local dependence and $\bta_0$ be the parameter vector of the same dimension as $\bta$ that has the same edge parameters and the same zero transitive edge parameters as $\bta$ but replaces all non-zero transitive edge parameters of $\bta$ by zeroes.
We first bound the expectation $\langle v,\, \bmu^\star\rangle$ and then turn to the concentration of $\langle v,\, \widehat{\bmu^\star}\rangle$ around its expectation $\langle v,\, \bmu^\star\rangle$.

\s

Concerning the expection $\langle v,\, \bmu^\star\rangle$,
note that by definition $v$ eliminates all mean-value parameters corresponding to transitive edges,
therefore $\langle v,\, \bmu^\star\rangle$ is the sum of mean-value parameters corresponding to edges, 
i.e.,
the expected number of edges.
Observe that it is not straightforward to bound the expected number of edges,
because the data-generating model induces dependence among edges.
We therefore follow an indirect approach to bounding the expected number of edges as follows.
Since $X_{i,j} \geq X_{i,j}\, \max_{k \in \mA_{\zs_i} \cup \mA_{\zs_j},\, k \neq i,j} X_{i,k}\, X_{j,k}$ almost surely and $0 < \bta_{i,j,2}\truth \leq 1$ provided $n$ is sufficiently large,
there exists $n_0 > 0$ such that,
for all $n > n_0$,
\beno
&&\langle v,\, \bmu^\star\rangle
\;\;=\;\; \dsum_{i<j}^n \mbE(X_{i,j})
\;\;\geq\;\; \dsum_{i<j}^n \bta_{i,j,2}\truth\; \mbE(X_{i,j})\s
\\
\gte \dsum_{i<j}^n \bta_{i,j,2}\truth\; \mbE\left(X_{i,j}\, \max\limits_{k \in \mA_{\zs_i} \cup \mA_{\zs_j},\, k \neq i,j} X_{i,k}\, X_{j,k}\right).
\ee
Let $f_{i,j}(\bx) = x_{i,j}\, \max_{k \in \mA_{\zs_i} \cup \mA_{\zs_j},\, k \neq i,j} x_{i,k}\, x_{j,k}$,
where $f_{i,j}$ is considered as a function of $x$ for fixed $\zs \in \mZ$.
Due to the dependence among edges,
the expectations $\mbE\; f_{i,j}$ are not available in closed form,
but can be bounded as follows.
By using well-known exponential-family results \citep[e.g., Corollary 2.5,][p.\ 37]{Br86} along with the definition of $\bta$ and $\bta_0$,
we obtain
\beno
\langle\bta - \bta_0,\, \bmu(\bta) - \bmu(\bta_0)\rangle
\= \dsum_{i<j}^n \bta_{i,j,2}\truth\; \left(\mbE\; f_{i,j} - \mbE_0^\star\, f_{i,j}\right) 
\gte 0,
\ee
where $\bta_{i,j,2}\truth > 0$ by assumption and $\mbE\; f_{i,j}$ and $\mbE_0^\star\, f_{i,j}$ denote the expectation of $f_{i,j}$ under $\bta$ and $\bta_0$,
respectively.
Thus,
\beno
\dsum_{i<j}^n \bta_{i,j,2}\truth\; \mbE\; f_{i,j} 
\gte \dsum_{i<j}^n \bta_{i,j,2}\truth\; \mbE_0^\star\, f_{i,j}.
\ee
Since,
for all $i, j$ and $k^\prime \in \mA_{\zs_i} \cup \mA_{\zs_j}$ such that $k^\prime \neq i, j$, 
we have $X_{i,j}\, \max_{k \in \mA_{\zs_i} \cup \mA_{\zs_j},\, k \neq i,j} X_{i,k}\, X_{j,k} \geq X_{i,j}\, X_{j,k^\prime}\, X_{i,k^\prime}$ almost surely,
we obtain
\beno
\mbE_0^\star\, f_{i,j}
\= \mbE_0^\star\left(X_{i,j}\, \max\limits_{k \in \mA_{\zs_i} \cup \mA_{\zs_j},\, k \neq i,j} X_{i,k}\, X_{j,k}\right)
\gte \mbE_0^\star\left(X_{i,j}\, X_{j,k^\prime}\, X_{i,k^\prime}\right).
\ee
To bound the expectation $\mbE_0^\star(X_{i,j}\, X_{j,k^\prime}\, X_{i,k^\prime})$,
observe that by setting all transitive edge parameters to zero,
$\bta_0$ reduces the model to an exponential family with edge indicators as sufficient statistics,
which implies that edge variables $X_{i,j}$ are independent under $\bta_0$.
Therefore,
\beno
\mbE_0^\star\left(X_{i,j}\, X_{j,k^\prime}\, X_{i,k^\prime}\right)
\= \mbE_0^\star(X_{i,j})\, \mbE_0^\star(X_{j,k^\prime})\, \mbE_0^\star(X_{i,k^\prime}).
\ee
To bound $\mbE_0^\star(X_{i,j})$,
note that,
using $n_k(\zs) > 0$ since $\zs \in \mG^\star \subset \mZ$,
we obtain
\beno
\mbE_0^\star\, (X_{i,j})
\gte 
\begin{cases}
\dfrac{1}{2\, n_{\zs_i}(\zs)^{- \min_{1 \leq k < l \leq K} (\bthetas_{k,l,1})}} & \mbox{if } \zs_i = \zs_j\s
\\
\dfrac{1}{2\, n^{- \min_{1 \leq k < l \leq K} (\bthetas_{k,l,1})}} & \mbox{if } \zs_i \neq \zs_j,\s
\end{cases}
\ee
where we used the fact that under $\bta_0$ the edge variables $X_{i,j} \in \{0, 1\}$ are independent and therefore the expectation $\mbE_0^\star(X_{i,j})$ is given by $\mbE_0^\star(X_{i,j}) = \mbox{logit}^{-1}(\bta_{i,j,1}\truth)$.
\hide{
Remark:
For all $n \geq 1$,
\beno
\mbE\; X 
\= p
\= \dfrac{\exp(\btheta \log n)}{1 + \exp(\btheta \log n)}
\= \dfrac{1}{\exp(-\btheta \log n) + 1}
\= \dfrac{1}{n^{-\btheta} + 1}
\gte \dfrac{1}{2\, n^{-\btheta}}.
\ee
}
As a result,
\beno
&& \dsum_{i<j}^n \bta_{i,j,2}\truth\; \mbE\; f_{i,j}
\;\;\geq\;\; \dsum_{i<j}^n \bta_{i,j,2}\truth\; \mbE_0^\star\, f_{i,j}\s
\\
\gte \dsum_{k=1}^K C_1\, n_k(\zs)^3\, \left[\bthetas_{2,k,k}\, \dfrac{\log n_k(\zs)}{n_k(\zs)}\right]\, \left[\dfrac{1}{8\, n_k(\zs)^{- 3\,  \min_{1 \leq k < l \leq K} (\bthetas_{k,l,1})}}\right]\s
\\
&+& \dsum_{k<l}^K C_2\, \min(n_k(\zs),\, n_l(\zs))^3\, \left[\bthetas_{k,l,2}\, \dfrac{\log n}{\max(n_k(\zs),\, n_l(\zs))}\right]\, \left[\dfrac{1}{8\, n^{- 3\,  \min_{1 \leq k < l \leq K} (\bthetas_{k,l,1})}}\right]\s
\\
\gte \dsum_{k<l}^K C_2\, \min(n_k(\zs),\, n_l(\zs))^3\, \left[\bthetas_{k,l,2}\, \dfrac{\log n}{\max(n_k(\zs),\, n_l(\zs))}\right]\, \left[\dfrac{1}{8\, n^{- 3\,  \min_{1 \leq k < l \leq K} (\bthetas_{k,l,1})}}\right]\s
\hide{
\\
\= \dsum_{k<l}^K C_4\, \bthetas_{k,l,2}\, \min(|\mA_k|,\, |\mA_l|)^2 \dfrac{\log n}{n^{- 3\, \min_{1 \leq k < l \leq K} (\bthetas_{k,l,1})}}\s
\\
\gte C_5\, K^2\, \min\limits_{1 \leq k < l \leq K}(\bthetas_{k,l,2})\, \dfrac{(\log n)^{3 + 2\, \alpha}}{n^{- 3\, \min_{1 \leq k < l \leq K} (\bthetas_{k,l,1})}}\s
\\
\gte C_5\, n^2\, \min\limits_{1 \leq k < l \leq K}(\bthetas_{k,l,2})\, \dfrac{\log n}{n^{-3\, \min_{1 \leq k < l \leq K} (\bthetas_{k,l,1})}}\s
}
\\
\gte C_5\, \min\limits_{1 \leq k < l \leq K} (\bthetas_{k,l,2})\; n^{2 + 3\, \min\limits_{1 \leq k < l \leq K} (\bthetas_{k,l,1})} \log n
\ee
using $C_6\, (\log n)^{1 + \alpha} \leq n_k(\zs) \leq C_7\, (\log n)^{1 + \alpha}$ since $\zs \in \mGs$ and $K \geq C\, n\, /\, (\log n)^{1 + \alpha}$. 
Thus,
there exists $n_0 > 0$ such that,
for all $n > n_0$,
\beno
\langle v,\, \bmu^\star\rangle
\;\geq\; \dsum_{i<j}^n \bta_{i,j,2}\truth\; \mbE\; f_{i,j}
\;\geq\; C_5\, \min\limits_{1 \leq k < l \leq K} (\bthetas_{k,l,2})\; n^{2 + 3\, \min\limits_{1 \leq k < l \leq K} (\bthetas_{k,l,1})} \log n.
\ee

Concerning the concentration of $\langle v,\, \widehat{\bmu^\star}\rangle$ around its expectation $\langle v,\, \bmu^\star\rangle$,
observe that $\vartheta = 0$,
because $v \in \{0, 1\}^{\dim(\bta)}$ reduces the vector of edge and transitive edge indicators $s(\bx)$ to the vector of edge indicators,
implying
\beno
|\langle v,\, s(\bx) - s(y)\rangle|
\= d(x, y) & \mbox{for all} & (x, y) \times \mX\times\mX.
\ee
By applying Theorem \ref{theorem.mean.value.parameters} with $\vartheta = 0$,
there exists $n_0 > 0$ such that,
for all $n > n_0$,
$\delta > 0$,
and $\epsilon > 0$,
\beno
\mbP\left(\left|\langle v,\, \widehat{\bmu^\star} - \bmu^\star\rangle\right| \;\geq\; \delta\, n^{1 + \epsilon}\right)
\lte 2\, \exp\left(- \dfrac{\delta^2\, n^{2\, \epsilon}}{C\, (\log n)^{4 + 4\, \alpha}}\right).
\ee

\hide{
\com An alternative approach based on the Law of Total Probability gives rise to more general results by extending the results from $\mG$ to $\mZ$:
\beno
&& \mbP\left(\left|\langle v(\zh),\, \widehat{\bmu^\star} - \bmu^\star\rangle\right| \;>\; \delta_3\, n^{1+\epsilon_3}\right)\s
\\
\= \mbP\left(\left|\langle v(\zh),\, \widehat{\bmu^\star} - \bmu^\star\rangle\right| \;>\; \delta_3\, n^{1+\epsilon_3} \mid \zh \in \mG\right)\, \mbP\left(\zh \in \mG\right)\s
\\
&+& \mbP\left(\left|\langle v(\zh),\, \widehat{\bmu^\star} - \bmu^\star\rangle\right| \;>\; \delta_3\, n^{1+\epsilon_3} \mid \zh \in \mB\right)\, \mbP\left(\zh \in \mB\right)\s
\\
\lte \mbP\left(\left|\langle v(\zh),\, \widehat{\bmu^\star} - \bmu^\star\rangle\right| \;>\; \delta_3\, n^{1+\epsilon_3} \mid \zh \in \mG\right) + \mbP\left(\zh \in \mB\right)\s
\\
\lte 2\, \exp\left(- \dfrac{\delta_3^2\, n^{2\, \epsilon_3 - 2 + 2\, \beta - 2\, \gamma}}{C_4\, (\log n)^{4 + 4\, \alpha}}\right) + \exp\left(- C_5\, n^{1 + \ithree}\, (\log n)^{1 + (1 + \alpha)\, (3 + \beta)}\right)
\ee
by using Proposition \ref{theorem.step2}.
}

}

\hide{
\lemma
\label{lemma.z.1}
Let $\truth \in \bTheta_0 \times \mZ_0$ and $(\btheta, \bz) \in \bTheta_0 \times \mZ_0$.
Then
\beno
\ell(\bthetas, \bzs; \bmu^\star) - \ell(\btheta, \bz; \bmu^\star)
\= KL(\bthetas, \zs;\, \btheta, \bz)
\gte 0
\ee
and $\mbE\;\ell$ is maximized by the truth $\truth$:
\beno
\ell(\bthetas, \bzs; \bmu^\star) \gte \ell(\btheta, \bz; \bmu^\star) & \mbox{for all} & (\btheta, \bz) \in \bTheta_0 \times \mZ_0.
\ee

\llproof \ref{lemma.z.1}.
By definition,
for all $(\btheta, \bz) \in \bTheta_0 \times \mZ_0$,
\beno
\ell(\btheta, \bz; \bmu^\star)
\= \mbE[\langle\bta(\btheta, \bz),\, s(\bX)\rangle - \psi(\bta(\btheta, \bz))]
\= \langle\bta(\btheta, \bz),\, \bmu^\star\rangle - \psi(\bta(\btheta, \bz)),
\ee
where $\bmu^\star = \mbE\; s(\bX)$.
Thus,
\beno
\ell(\bthetas, \bzs; \bmu^\star) - \ell(\btheta, \bz; \bmu^\star)
\= \langle\bta\truth - \bta(\btheta, \bz),\, \bmu^\star\rangle - \psi\truth + \psi(\bta(\btheta, \bz)),
\ee
which is equal to the Kullback-Leibler divergence from the exponential-family distribution with natural parameter vector $\bta\truth$ to the exponential-family distribution with natural parameter vector $\bta(\btheta, \bz)$.
Observe that both distributions have the same support $\conv(\mS)$.
By Jensen's inequality,
\beno
\ell(\bthetas, \bzs; \bmu^\star) - \ell(\btheta, \bz; \bmu^\star)
\= KL(\bthetas, \zs;\, \btheta, \bz)
\gte 0,
\ee
implying that $\mbE\;\ell$ is maximized by the truth $\truth$:
\beno
\ell(\bthetas, \bzs; \bmu^\star) \gte \ell(\btheta, \bz; \bmu^\star) & \mbox{for all} & (\btheta, \bz) \in \bTheta_0 \times \mZ_0.
\ee
}

\hide{

\ttproof \ref{theorem.selection}.
By C.\ref{a.graph},
the conditional independence graph is fully specified when the block structure $\zs$ is known.
If the block structure $\zs$ is unknown,
it can be estimated by $\zh$,
which produces an estimator of the conditional independence graph as a by-product.
Suppose $\delta(\zs, \zh) = M$.
Then the number of false-negative and false-positive edges in the estimated conditional independence graph can be bounded by bounding the number of false-negative and false-positive edges due to moving $M$ nodes from true to false blocks.

\s

{\bf (a) Fraction of false-negative edges.}
Moving node $i$ from its true block $\mA_{\zs_i}$ to a false block can remove up to $C_1\, (\log n)^{3 + 3\, \alpha}$ edges in the conditional independence graph,
because node $i$ is involved in $n_k(\zs_i) - 1 \leq n_k(\zs_i)$ edge variables $X_{i,j}$ in block $\mA_{\zs_i}$ and each edge variable $X_{i,j}$ can have up to ${n_k(\zs_i) \choose 2} - 1 \leq n_k(\zs_i)^2$ neighbors in the true conditional independence graph within block $\mA_{\zs_i}$,
and the size $n_k(\zs_i)$ of block $\mA_{\zs_i}$ under $\zs$ is bounded above by $n_k(\zs_i) \leq C_2\, (\log n)^{1 + \alpha}$ since $\zs \in \mGs$.
Therefore,
the number of false-negative edges in the estimated conditional independence graph based on the estimated block structure $\zh$ at distance $\delta(\zs, \zh) = M$ from the truth $\zs$ can be bounded above by $C_1\, M\, (\log n)^{3 + 3\, \alpha}$. 
By C.\ref{a.graph},
each block $\mA_k$ under $\zs$ contributes at least $C_3\, (\log n)^{3 + 3\, \alpha}$ edges to the true conditional independence graph,
because there are ${|\mA_k| \choose 2} \geq C_4\, (\log n)^{2 + 2\, \alpha}$ edge variables in block $\mA_k$ and each edge variable is connected to $C_5\, n(\log n)^{1 + \alpha}$ edges in the conditional independence graph.
Therefore,
the total number of possible false-negative edges in the estimated conditional independence graph is at least $C_3 \sum_{k=1}^K (\log n)^{3 + 3\, \alpha} = C_3\, K\, (\log n)^{3 + 3\, \alpha} \geq C_6\, n\, (\log n)^{2 + 2\, \alpha}$ using $K \geq C\, n / (\log n)^{1 + \alpha}$.
As a result,
the fraction of false-negative edges $f_N$ in the estimated conditional independence graph is bounded above by
\beno
f_N
\lte \dfrac{C_1\, M\, (\log n)^{3 + 3\, \alpha}}{C_3\, n\, (\log n)^{2 + 2\, \alpha}}
\= \dfrac{C_7\, M\, (\log n)^{1 + \alpha}}{n}.
\ee
\hide{
Let $\epsilon > 0$.
We are interested in the event
\beno
\dfrac{C_7\, M\, (\log n)^{1 + \alpha}}{n}
\gte f_N 
\gte \epsilon.
\ee
}
An application of Theorem \ref{theorem.step3} shows that there exist $C > 0$ and $n_0 > 0$ such that,
for all $n > n_0$,
\beno
\mbP\left(f_N \;\geq\; \epsilon\right)
\lte \mbP\left(\dfrac{C_7\, M\, (\log n)^{1 + \alpha}}{n} \;\geq\; \epsilon\right)
\hide{
\lte \mbP\left(\dfrac{M}{n} \;\geq\; \dfrac{\epsilon}{C_7\, (\log n)^{1 + \alpha}}\right)
}
\lte 2\, \exp\left(- \dfrac{\epsilon^2\, n^{2\, \beta}}{(\log n)^{2 + 2\, \alpha}}\right).
\ee

\s

{\bf (b) Fraction of false-positive edges.}
\hide{
Starting at $\zs$,
move nodes from true blocks to false blocks.
}
Consider the $m$-th move of the $m$-th node $i$ from its true block $\mA_{\zs_i}$ to some false block $\mA_l$.
The false block $\mA_l$ contains up to 
\beno
|\mA_l| + (m - 1) 
\lte C_1\, (\log n)^{1 + \alpha} + m
\ee
nodes before the $m$-th move and up to 
\beno
\dis{|\mA_l| + (m - 1) \choose 2}
\hide{
\lte \dfrac{(|\mA_l| + m)^2}{2}
\lte \dfrac{|\mA_l|^2 + 2\, |\mA_l|\, m + m^2}{2}\s
\\
&&\lte \dfrac{2\, m^2\, |\mA_l|^2 + 2\, m^2\, |\mA_l|^2 + 2\, m^2\, |\mA_l|}{2}\s
\\
&&
}
\lte C_2\, m^2\, (\log n)^{2 + 2\, \alpha}
\ee
edge variables before the $m$-th move using $\norm{n(z)}_\infty \leq C\, (\log n)^{1 + \alpha}$ since $\zs \in \mGs$.
The $m$-th move adds up to 
\beno
|\mA_k| + m - 1 
\lte C_3\, (\log n)^{1 + \alpha} + m  
\ee
edge variables $X_{i,j}$ to the false block $\mA_l$.
Therefore,
the $m$-th move adds up to 
\beno
C_2\, m^2\, (\log n)^{2 + 2\, \alpha}\, \left[C_3\, (\log n)^{1 + \alpha} + m\right] 
\hide{
\= C_4\, m^2\, (\log n)^{3 + 3\, \alpha}  + C_5\, m^3\, (\log n)^{2 + 2\, \alpha}\s
\\
}
\lte C_5\, m^3\, (\log n)^{3 + 3\, \alpha}
\ee
false-positive edges in the conditional independence graph between old and new edge variables in the false block $\mA_l$ and up to
\beno
\dis{C_3\, (\log n)^{1 + \alpha} + m  \choose 2}
\hide{
\lte \dfrac{(C_3\, (\log n)^{1 + \alpha} + m)^2}{2}
\vspace{0.125cm}
\\
\lte C_4\, \left[(\log n)^{2 + 2\, \alpha} + 2\, m\, (\log n)^{1 + \alpha} + m^2\right]\s
\\
}
\lte C_6\, m^2\, (\log n)^{2 + 2\, \alpha}
\ee
false-positive edges in the conditional independence graph between new and new edge variables in the false block $\mA_l$.
Thus,
the total number of false-positive edges in the conditional independence graph due to the $m$-th move is bounded above by
\beno
C_5\, m^3\, (\log n)^{3 + 3\, \alpha}  + C_6\, m^2\, (\log n)^{2 + 2\, \alpha}
\lte C_7\, m^3\, (\log n)^{3 + 3\, \alpha}
\ee
and the total number of false-positive edges in the conditional independence graph due to $M$ moves is bounded above by
\beno
\dsum_{m=1}^M C_7\; m^3\, (\log n)^{3 + 3\, \alpha}
\= C_7\, \left(\dfrac{M\, (M + 1)}{2}\right)^2\, (\log n)^{3 + 3\, \alpha}
\lte C_8\, M^4\, (\log n)^{3 + 3\, \alpha}
\ee
provided $n$ is sufficiently large.
Observe that under C.\ref{a.graph} the number of edges in the true conditional independence graph is bounded above by $C\, K\, (\log n)^{4 + 4\, \alpha}$,
which implies that the total number of possible false-positive edges is at least
\beno
\dis{{n \choose 2} \choose 2} - C\, K\, (\log n)^{4 + 4\, \alpha}
\hide{
\gte \dis{{n \choose 2} \choose 2} - C\, n\, (\log n)^{4 + 4\, \alpha}
}
\gte \dfrac{C\, n^4}{4}
\ee
provided $n$ is sufficiently large.
Therefore,
the fraction of false-positive edges in the conditional independence graph is bounded above by
\beno
f_P 
\lte 
\hide{
\dfrac{C_8\, M^4\, (\log n)^{3 + 3\, \alpha}}{\dfrac{n^4}{4}}.
\= 
}
\dfrac{4\, C_9\, M^4\, (\log n)^{3 + 3\, \alpha}}{n^4}.
\ee
\hide{
Let $\epsilon > 0$.
We are interested in the event
\beno
\dfrac{4\, C_9\, M^4\, (\log n)^{3 + 3\, \alpha}}{n^4}\s
\gte f_P 
\gte \epsilon.
\ee
}
An application of Theorem \ref{theorem.step3} shows that there exist $C > 0$ and $n_0 > 0$ such that,
for all $n > n_0$,
\beno
\mbP\left(f_P \;\geq\; \epsilon\right)
\lte \mbP\left(\dfrac{4\, C_9\, M^4\, (\log n)^{3 + 3\, \alpha}}{n^4} \;\geq\; \epsilon\right)
\hide{
\= \mbP\left(\dfrac{M}{n} \;>\; \dfrac{\epsilon^{1/4}}{4\, C_8\, (\log n)^{(3 / 4)\, (1 + \alpha)}}\right)
}
\lte 2\, \exp\left(- \dfrac{\epsilon^{1/2}\, n^{2\, \beta}}{C\, (\log n)^{(3 / 2)\, (1 + \alpha)}}\right).
\ee

}

\hide{

\section{Size-dependent parameterizations}
\label{size.dependent.parameterizations}

We first discuss size-dependent parameterizations and then show that these parameterizations satisfy \cthree.

\s

There is a long tradition of size-dependent parameterizations in random graph models:
size-dependent parameterizations have been used in classic Bernoulli random graph models since the 1950s \citep[e.g.,][]{ErRe59,AlSp92,Bo98,Jo99,MoRe02,BoMu08,Lo12} and were introduced to the domain of exponential-family random graph models by \citep{Jo99}, \citep{KrHaMo11}, and \citep{KrKo14}.
The size-dependent parameterizations are motivated by the fact that real-world networks tend to be sparse in the sense that the observed number of edges is much smaller than the number of possible edges ${n \choose 2}$.
A possible explanation is that in practice nodes may not be able to maintain $n-1$ edges and therefore the expected number of edges of nodes $i$ should satisfy $\mbE(\sum_{j:\, j \neq i}^n X_{i,j}) \ll n-1$.
In the simplest case where edge variables $X_{i,j}$ are independent and identically distributed Bernoulli$(p)$ random variables---which is an exponential-family random graph model with the number of edges as a sufficient statistic,
the classic notion of sparsity postulates that $\mbE(\sum_{j:\, j \neq i}^n X_{i,j}) = (n - 1)\, p < C$ for some $C > 0$,
reflecting the idea that nodes cannot maintain an unbounded number of edges (e.g., friendships in a social network).
Therefore,
the probability of an edge $p$ should be of order $n^{-1}$ and the natural parameter of the exponential family,
which is the log odds of $p$, 
should be of order $\log n$:
\be
\label{sparse.edge.term}
\bta_1(\btheta) 
\lte \theta_1 \log n, 
&& \theta_1 \;<\; 0,
&& n \;\geq\; 1.
\ee
\citep{Jo99} extended the sparse edge model to the sparse edge-and-triangle model by assuming that the natural edge parameter satisfies \eqref{sparse.edge.term} while the natural triangle parameter satisfies
\be
\label{sparse.triangle.term}
\bta_2(\btheta) 
\lte \log\left(1 + \theta_2\, \dfrac{\log n}{n}\right)
\lte \theta_2\, \dfrac{\log n}n,
&& \theta_2 \;>\; 0,
&& n \;\geq\; 1.
\ee
As a result of \eqref{sparse.edge.term} and \eqref{sparse.triangle.term},
both the edge term $\bta_1(\btheta)\, \sum_{i<j}^n x_{i,j}$ and triangle term
\linebreak
$\bta_2(\btheta)\, \sum_{i<j<k}^n x_{i,j}\, x_{j,t}\, x_{i,t}$ are of order $O(n^2 \log n)$.
\citep{KrHaMo11}, and \citep{KrKo14} used related parameterizations motivated by sparsity considerations and \citep{ChDi11} used size-dependent parameterizations to ensure the existence of dense-graph limits (without the $\log n$-term, which is motivated by sparsity)---note that size-dependent weights can be pulled into either the natural parameter vector or the sufficient statistics vector due to the fact that the canonical form of exponential families is not unique, and some of the cited authors have followed the former approach while others have followed the latter approach.
The following assumption extends the ideas of \citep{Jo99} to finite-dimensional exponential families governing within- and between-block subgraphs.


\assumption
\label{size}
{\em Size-dependent parameterizations.}
Consider an exponential family with local dependence satisfying \ctwo.
Without loss,
assume that the marginal exponential families governing within- and between-block subgraphs are minimal in the sense of \citep{BN78} and \citep{Br86} and let $\bta_{\min,k,l}(\btheta,z)$ be the natural parameter vector of the marginal minimal exponential family governing the subgraph of blocks $k \leq l$ and $s_{\min,k,l}(x_{k,l})$ be the corresponding sufficient statistics vector.
Denote by $\bta_{\min,k,l,i}(\btheta, \bz)$ and $s_{\min,k,l,i}(x_{k,l})$ coordinates $i$ of $\bta_{\min,k,l}(\btheta,z)$ and $s_{\min,k,l}(x_{k,l})$,
respectively.
Assume that there exist $B_i \geq 0$ and $C_i \geq 0$ such that,
for all pairs of blocks $k \leq l$ and all $x \in \mX$,
the change in sufficient statistics satisfies,
for all $(x, y) \in \mX \times \mX$,
\beno
|s_{\min,k,l,i}(\bx) - s_{\min,k,l,i}(y)|
\lte C_i\; d(x_{k,l}, y_{k,l})\, \max(|\mA_k|,\, |\mA_l|)^{B_i},
\ee
where $d: \mX_{k,l} \times \mX_{k,l}$ denotes the Hamming metric,
where $\mX_{k,l}$ denotes the sample space of the subgraph corresponding to blocks $\mA_k$ and $\mA_l$ and $n_k(\zs)$ and $n_l(\zs)$ denote the sizes of blocks $\mA_k$ and $\mA_l$ ($k \leq l$).
Then the natural parameter $\bta_{\min,k,l,i}(\btheta, \bz)$ corresponding to $s_{\min,k,l,i}(x_{k,l})$ is said to satisfy weak sparsity assumptions if
\beno
\label{consistent.scaling}
|\bta_{\min,k,l,i}(\btheta, \bz)| 
\lte |\btheta_{i,k,l}|\, \dfrac{\log \max(|\mA_k|,\, |\mA_l|)}{\max(|\mA_k|,\, |\mA_l|)^{B_i}}
\ee
and strong sparsity assumptions if
\beno
\label{consistent.scaling2}
|\bta_{\min,k,l,i}(\btheta, \bz)|
\lte |\btheta_{i,k,l}|\, \dfrac{\log n}{\max(|\mA_k|,\, |\mA_l|)^{B_i}},
\ee
where $\max(|\mA_k|,\, |\mA_l|) > 0$ since $\bz \in \mZ_0$.

\s

Condition C.\ref{size} implies that all model terms corresponding to blocks $k \leq l$ are of order $O(\max(|\mA_k|,\, |\mA_l|)^2 \log \max(|\mA_k|,\, |\mA_l|))$ in the case of weak sparsity or $O(\max(|\mA_k|,\, |\mA_l|)^2 \log n)$ in the case of strong sparsity.
Examples of parameterizations inducing weak sparsity are within- and between-block edge and triangle terms of the form
\beno
\label{weights}
\bta_{\min,k,l,1}(\btheta, \bz)
= 
\begin{cases}
\btheta_{k,l,k} \log |\mA_k| & k = l \mbox{  (edges)}\s
\\
\btheta_{k,l,1} \log \max(|\mA_k|,\, |\mA_l|) & k < l \mbox{  (edges)}
\end{cases}
\ee
and
\beno
\bta_{\min,k,l,2}(\btheta, \bz)
=
\begin{cases}
\btheta_{k,k,2}\, \dfrac{\log |\mA_k|}{|\mA_k|} & k = l \mbox{  (triangles)}\s
\\
\btheta_{k,l,2}\, \dfrac{\log \max(|\mA_k|,\, |\mA_l|)}{\max(|\mA_k|,\, |\mA_l|)} & k < l \mbox{  (triangles).}
\end{cases}
\ee
The example includes sparse edge term of order $O(\norm{n(z)}_\infty^2 \log \norm{n(z)}_\infty)$ and ensures that all other terms are of the same order.

\s

Condition C.\ref{size} implies that the natural parameter vector $\bta$ can be written as
\beno
\bta(\btheta, \bz) 
\= \bA(\bz)\, \bm{b}(\btheta),
\ee
where $A: \mZ \mapsto \mR^{\dim(\bta) \times \rrr}$ denotes a $\dim(\bta) \times \rrr$-matrix-valued function and $\bm{b}: \bTheta \mapsto \mR^{\rrr}$ may be a linear or non-linear function of parameter vector $\btheta \in \bTheta$.
We note that $\bA(\bz)$ contains the size-dependent weights and satisfies $\norm{\bA(\bz)}_\infty \leq C_1 \log n$ whereas $\btheta = (\btheta_{k,l})$ contains the size-independent parameters satisfying $\rrr \leq C_2\, K^2$,
thus $\bA(\bz)$ and $\btheta$ satisfy \ctwo.

\hide{
As a result,
the exponential-family densities \eqref{local.ergm} can be written as
\be
\label{local.ergm.2}
p_{\bta}(\bx)
&=& \exp\left(\langle\bta,\, s(\bx)\rangle - \psi(\bta)\right)
&=& \exp\left(\langle\btheta,\, t(x, \bz)\rangle - \psi(\bta(\btheta, \bz))\right),
\ee
where $t(x, \bz) = \langle \bA(\bz),\, s(\bx)\rangle$.
We assume without loss that the exponential family of densities on the right-hand side of \eqref{local.ergm.2} is a minimal representation,
because if it was non-minimal,
it could be reduced to a minimal representation \citep[e.g.,][Theorem 1.9, p.\ 13]{Br86}.
}

\s

We show that size-dependent parameterizations along these lines satisfy \ctwo.

\lemma
\label{lemma.lipschitz}
The size-dependent parameterizations described in C.\ref{size},
both the weak and strong form of sparsity,
satisfy \ctwo,
i.e.,
there exists $C > 0$ such that
\beno
\left|\langle\bta(\btheta, \bz),\, s(\bx) - s(y)\rangle\right|
\lte C\, d(x, y) \log n 
&& \mbox{for all} & (x, y) \in \mX\times\mX.
\ee

\llproof \ref{lemma.lipschitz}.
If the exponential family is not minimal in the sense of \citep{BN78} and \citep{Br86},
we reduce it to a minimal exponential family;
note that all non-minimal exponential families can be reduced to minimal exponential families 
\citep[e.g.,][Theorem 1.9, p.\ 13]{Br86}.
Let $L$ be the dimension of the finite-dimensional marginal minimal exponential family governing the subgraph of blocks $k \leq l$.
Denote by $\bta_{\min,k,l}(\btheta,z)$ the natural parameter vector of the minimal exponential family governing the subgraph of blocks $k \leq l$ and denote by $s_{\min,k,l}(\bx)$ the corresponding sufficient statistics vector.
Let $\bta_{\min,k,l,i}(\btheta, \bz)$ and $s_{\min,k,l,i}(\bx)$ be the $i$-th coordinate of $\bta_{\min,k,l}(\btheta,z)$ and $s_{\min,k,l}(\bx)$,
respectively.
By Jensen's inequality and the assumption that $\bTheta$ is a compact subset of $\bTheta = \mR^{\rrr}$,
\beno
\label{lipschitz}
\left|\langle\bta(\btheta, \bz),\, s(\bx) - s(y)\rangle\right|
\= \left|\dsum_{k \leq l}^K \langle\bta_{\min,k,l}(\btheta, \bz),\, s_{\min,k,l}(\bx)\rangle - \dsum_{k \leq l}^K \langle\bta_{\min,k,l}(\btheta, \bz),\, s_{\min,k,l}(y)\rangle\right|\s
\\
\lte \dsum_{k \leq l}^K \dsum_{i=1}^L \left|\bta_{\min,k,l,i}(\btheta, \bz)\right|\, \left|s_{\min,k,l,i}(\bx) - s_{\min,k,l,i}(y)\right|\s
\\
\lte \dsum_{k < l}^K C_2\, d(x_{k,l}, y_{k,l}) \log n
\;\;=\;\; C_2\, d(x, y) \log n
\ee
using $|\mA_k| > 0$ for all $\bz \in \mZ_0$ and $k = 1, \dots, K$ and
\beno
&& |\bta_{\min,k,l,i}(\btheta, \bz)|\, |s_{\min,k,l,i}(\bx) - s_{\min,k,l,i}(y)|\s
\\
\lte |\btheta_{i,k,l}|\, \dfrac{\log n}{\max(|\mA_k|,\, |\mA_l|)^{B_i}}\, C_1\, d(x_{k,l}, y_{k,l})\, \max(|\mA_k|,\, |\mA_l|)^{B_i}
\;\;\leq\;\; C_2\, d(x_{k,l}, y_{k,l})\, \log n,
\ee 
which covers both weak and strong sparsity.

}

\hide{
As a result,
we can write
\beno
\langle\bta(\btheta, \bz),\, \bmu\rangle
\= \langle\btheta,\, \bmu(\bz)\rangle
\= \dsum_{k \leq l}^K \langle\btheta_{k,l},\, \bmu_{k,l}(\bz)\rangle\s
\\
\langle\bta(\btheta, \bz),\, s(\bx)\rangle
\= \langle\btheta,\, s(\bx, \bz)\rangle
\= \dsum_{k \leq l}^K \langle\btheta_{k,l},\, s_{k,l}(\bx, \bz)\rangle,
\ee
where $\bmu(\bz) = \bA(\bz)^\top \bmu$ ($\bmu \in \mMs$) and $s(\bx, \bz) = \bA(\bz)^\top s(\bx)$ ($s(\bx) \in \mMs$).
}

\s

\llproof \ref{good.event}.
Since the data-generating natural parameter vector $\bta^\star \in \etaspace \subseteq \mbox{int}(\fullspace)$ is in the interior $\mbox{int}(\fullspace)$ of the natural parameter space $\fullspace$,
the expectation $\mbE\, s(\bX)$ exists \citep[][Theorem 2.2, pp.\ 34--35]{Br86} and so does the expectation
\linebreak
$\mbE\, \ell(\btheta, \bz; s(\bX)) = \ell(\btheta, \bz; \mbE\, s(\bX))$.
We want to bound
\beno
\mbP\left(\bX \in \mX \setminus \mG\right)
\= \mbP\left(|\ell(\bthetas, \bzs; s(\bX)) -  \ell(\bthetas, \bzs; \bmu^\star)|\; \geq\; \alpha\; \uuushort\right),
\ee
where 
\beno
\uuushort 
\= \uuu.
\ee
Bounding the probability of deviations of the form $|\ell(\bthetas, \bzs; s(\bX)) -  \ell(\bthetas, \bzs; \bmu^\star)|$ is equivalent to bounding the probability of deviations of the form $|f(\bX) - \mbE\, f(\bX)|$,
where 
\beno
\label{ff}
f(\bX) 
\= \langle\bta\truth,\, s(\bX)\rangle,
&& \mbE\, f(\bX) 
\= \langle\bta\truth,\, \bmu^\star\rangle.
\ee
We note that $f: \mX \mapsto \mR$ is considered as a function of $\bX$ for fixed $\truth \in \bTheta_0 \times \mZ_0$ and that $\psi(\bta\truth)$ cancels.
Observe that by condition \cthree{ }there exist $A_2 > 0$ and $n_0 > 0$ such that,
for all $n > n_0$,
the Lipschitz coefficient of $f(\bX)$ satisfies $\norm{f}_{\lip} \leq A_2\; \size$.
Thus,
by applying Lemma \ref{proposition.concentration} to deviations of size $t = \alpha\, \uuushort$,
there exist $C_0 > 0$ and $n_0 > 0$ such that,
for all $n > n_0$,
\beno
\mbP\left(\bX \in \mX \setminus \mG\right)
\hide{
\= \mbP\left(|\ell(\bthetas, \bzs; s(\bX)) -  \ell(\bthetas, \bzs; \bmu^\star)|\; \geq\; \alpha\; \uuushort\right)\s
\\
}
\lte 2\, \exp\left(- \dfrac{\alpha^2\; \uuushort^2}{C_0\, n^2\, \norm{\mA}_\infty^4\, \size^2}\right).
\ee
By assumption \eqref{a.condition} of Proposition \ref{p.z.1},
there exists,
for all $C_1 > 0$,
however large,
$n_1 > 0$ such that,
for all $n > n_1$,
\beno
\uuushort
\gte C_1\; n^{3/2}\; \norm{\mA}_\infty^2\, \size\, \sqrt{\log n}.
\ee
Therefore,
there exists $C > 0$ such that,
for all $n > \max(n_0, n_1)$,
\beno
\mbP\left(\bX \in \mX \setminus \mG\right)
\hide{
\= \mbP\left(|\ell(\bthetas, \bzs; s(\bX)) -  \ell(\bthetas, \bzs; \bmu^\star)|\; \geq\; \alpha\; \uuushort\right)\s
\\
}
\lte 2\, \exp\left(- \alpha^2\; C\, n \log n\right).
\ee

\s

\llproof \ref{mle.exists}.
In the following,
we confine attention to $\bx \in \mG$,
because we are interested in the existence of the restricted maximum likelihood estimator $\estimate$ in the event $\mG$.
For any $\bx \in \mG$ and any $\bz \in \mZ_0$,
let
\beno
\bthetah(\bz) 
\= \argmax\limits_{\btheta\in\bTheta_0} \ell(\btheta, \bz; s(\bx)).
\ee
Observe that,
for any $\bx \in \mG$ and any $\bz \in \mZ_0$,
the loglikelihood function $\ell(\btheta, \bz; s(\bx))$ is upper semicontinuous on $\bTheta_0$ by condition \cone.
In addition,
by condition \cfour{ }there exist $A, B, C > 0$ such that the $\dim(\btheta) \leq A\, n$-dimensional parameter space $\bTheta_0$ can be covered by $\exp(C\, n)$ closed balls with centers $\btheta \in \bTheta$ and radius $B > 0$.
As a result, 
for any $\bx \in \mG$ and any $\bz \in \mZ_0$,
$\ell(\btheta, \bz; s(\bx))$ assumes a maximum on $\bTheta_0$ and hence the maximizer $\bthetah_l(\bz)$ exists and is unique by condition \czero{ }and the assumption that the exponential family is minimal,
which can be assumed without loss \citep[][Theorem 1.9, p.\ 13]{Br86}.
Since,
for any $\bz \in \mZ_0$,
$\widehat\btheta(\bz)$ exists,
so does $\estimate$.

Last, 
but not least,
since $\estimate$ exists for all $\bx \in \mG$,
$\estimate$ exists with at least probability $\mbP\left(\bX \in \mG\right)$.
By Lemma \ref{good.event},
there exist $C_0 > 0$ and $n_0 > 0$ such that,
for all $n > n_0$,
\beno
\mbP\left(\bX \in \mG\right)
\gte 1 - 2\, \exp\left(- \alpha^2\; C_0\, n \log n\right).
\ee
Therefore,
for all $n > n_0$,
$\estimate$ exists with at least probability\linebreak $1 - 2\, \exp\left(- \alpha^2\; C_0\, n \log n\right)$.

\ccproof \ref{c.canonical}.
To show that conditions \czero---\cthree{ }are satisfied,
note that $\bta: \bTheta \times \mZ \mapsto \etaspace$ is separable in the sense that $\bta(\btheta, \bz) = \bA(\bz)\, \bm{b}(\btheta)$ and hence $\bta(\btheta, \bz)$ can be reduced to $\bta(\btheta)$ by absorbing $\bA(\bz)$ into the sufficient statistics vector.
In addition,
since the exponential family is canonical,
$\bta(\btheta)$ can be reduced to $\bta(\btheta) = \btheta$.
Condition \czero{ }is satisfied because $\bta(\btheta) = \btheta$.
Condition \cone{ }follows from $\bta(\btheta) = \btheta$ and the upper semicontinuity of canonical exponential-family loglikelihood functions \citep[][Lemma 5.3, p.\ 146]{Br86}.
To show that condition \ctwo{ }holds,
observe that
\beno
|\langle\bta(\btheta_1, \bz) - \bta(\btheta_2, \bz),\, \bmu\rangle|
\= |\langle\btheta_1 - \btheta_2,\, \bmu(\bz)\rangle|,
\ee
where $\bmu(\bz) = \bA(\bz)^\top \bmu$ ($\bmu \in \mMs$).
We can therefore write
\beno
|\langle\btheta_1 - \btheta_2,\, \bmu(\bz)\rangle|
\= \dsum_{k \leq l}^K |\langle\btheta_{1,k,l} - \btheta_{2,k,l},\, \bmu_{k,l}(\bz)\rangle|\s
\\
\lte \dsum_{k \leq l}^K \norm{\btheta_{1,k,l}-\btheta_{2,k,l}}_1\, \norm{\bmu_{k,l}(\bz)}_\infty\s
\\
\lte \dsum_{k \leq l}^K \sqrt{\dim(\btheta_{k,l})}\; \norm{\btheta_{1,k,l}-\btheta_{2,k,l}}_2\, \norm{\bmu_{k,l}(\bz)}_\infty.
\ee
Since the parameter vectors $\btheta_{k,l}$ are finite-dimensional,
the parameter space $\bTheta_0$ is compact, 
and the random graph is dense in the sense that\linebreak
$\uuu = C_0\, {n \choose 2}$ ($C_0 > 0$),
condition \ctwo{ }is satisfied as long as $\norm{\bmu_{k,l}(\bz)}_\infty \leq C_1\, L_k(\bz)\, L_l(\bz)$ ($C_1 > 0$) for all $\bmu \in \mMs$ and all $\bz \in \mZ_0$.
The same argument shows that
\beno
|\langle\bta(\btheta, \bz),\, s(\bx_1) - s(\bx_2)\rangle|
\;\;=\;\; |\langle\btheta,\, s(\bx_1, \bz) - s(\bx_2, \bz)\rangle|\s
\\
\hspace{3cm}=\;\; \dsum_{k \leq l}^K |\langle\btheta_{k,l},\, s_{k,l}(\bx_1, \bz) - s_{k,l}(\bx_2, \bz)\rangle|\s
\\
\hspace{3cm}\leq\;\; \dsum_{k \leq l}^K \sqrt{\dim(\btheta_{k,l})}\, \norm{\btheta_{k,l}}_2\, \norm{s_{k,l}(\bx_1, \bz) - s_{k,l}(\bx_2, \bz)}_\infty.
\ee
As a result,
condition \cthree{ }is satisfied as long as
$\sum_{k \leq l}^K \norm{s_{k,l}(\bx_1, \bz)-s_{k,l}(\bx_2, \bz)}_\infty$\linebreak $\leq C_2\, d(\bx_{1}, \bx_{2})\, L(\bz)$ for all $(\bx_1, \bx_2) \in \mX\times\mX$ and all $\bz \in \mZ_0$.

\s

\ccproof \ref{c.curved}.
To streamline the presentation,
we assume the following:
\bi
\item We take advantage of the fact that $\bta: \bTheta \times \mZ \mapsto \etaspace$ is separable in the sense that $\bta(\btheta, \bz) = \bA(\bz)\, \bm{b}(\btheta)$ and reduce $\bta(\btheta, \bz)$ to $\bta(\btheta)$ by absorbing $\bA(\bz)$ into the sufficient statistics vector.
\item Since,
under the curved exponential-family random graph model \eqref{geo} described in Section \ref{mainresults}, 
between- and within-block edge terms cannot violate conditions \czero---\cthree,
we assume that there is a single block without edge terms but with geometrically weighted model terms of the form \eqref{geo},
so that we can write
\beno
\langle\bta(\btheta, \bz),\, \bmu\rangle
\= \langle\bta(\btheta),\, \bmu(\bz)\rangle
\= \dsum_{t=1}^T \eta_t(\btheta)\, \mu_t(\bz)\s
\\
\langle\bta(\btheta, \bz),\, s(\bx)\rangle
\= \langle\bta(\btheta),\, s(\bx, \bz)\rangle
\= \dsum_{t=1}^T \eta_t(\btheta)\, s_t(\bx, \bz),
\ee
where $\bmu(\bz) = \bA(\bz)^\top \bmu$ ($\bmu \in \mMs$) and $s(\bx, \bz) = \bA(\bz)^\top s(\bx)$ ($s(\bx) \in \mMs$).
Throughout,
we drop the subscript $k$---which indexes\linebreak 
blocks---from all block-dependent quantities,
because there is a single block.
\hide{
The coordinates $\eta_t(\btheta)$ of the within-block parameter vector $\bta(\btheta)$ are then given by
\beno
\eta_{t}(\btheta)
\= \theta_1\, \left\{\theta_2\, \left[1 - \left(1 - \dfrac{1}{\theta_2}\right)^t\right]\right\},
&& t = 1, \dots, T,
&& T \geq 2,
\ee
where we dropped the subscript $k$ of all quantities,
because there is one block.
}
\item We assume that the parameter $\theta_1$ of the within-block edge term and the base parameter $\theta_2$ of the within-block geometrically weighted model term are given by $\theta_1 = 0$ and $\theta_2 = 1$, respectively, and drop the subscript of $\theta_3$, 
i.e.,
we write $\theta$ rather than $\theta_3$.
\ei
The extension to more than one block and $(\theta_1, \theta_2) \in \mR \times \mR$ is straightforward.

\hide{
Remark:
Let $\btheta = (\theta_1, \theta_2)$.
Then,
by the triangle inequality,
\beno
\norm{\hat\btheta - \btheta}_2
\= \norm{(\hat\theta_1, \hat\theta_2) - (\theta_1, \theta_2)}_2\s
\\
\= \norm{((\hat\theta_1, 0) - (\theta_1, 0)) + ((0, \hat\theta_2) - (0, \theta_2))}_2\s
\\
\lte \norm{(\hat\theta_1, 0) - (\theta_1, 0)}_2 + \norm{(0, \hat\theta_2) - (0, \theta_2)}_2\s
\\
\= \norm{\hat\theta_1 - \theta_1}_2 + \norm{\hat\theta_2 - \theta_2}_2.
\ee
In general,
\beno
\norm{\hat\bta - \bta}_2
\lte \dsum_{k=1}^K \norm{\hat\bta_k - \bta_k}_2.
\ee
}

Under the assumptions outlined above,
the coordinates $\eta_t(\theta)$ of the single within-block natural parameter vector $\bta(\theta)$ can be written as 
\be
\label{reduced.form}
\eta_{t}(\theta)
\= \theta\, \left[1 - \left(1 - \dfrac{1}{\theta}\right)^t\right]
\= \theta - \theta\, \beta(\theta)^t,
&& t = 1, \dots, T,
\ee
where 
\beno
\beta(\theta)
\= 1 - \dfrac{1}{\theta}.
\ee
The parameter space $\Theta$ is given by
\beno
\Theta
\= \left\{\theta \in \mR:\;\; \dfrac12 < \theta < B,\;\; \psi(\bta(\theta, \bz)) < \infty\right\}, 
& B > \dfrac12.
\ee
A helpful observation is that the coordinates $\eta_t(\theta)$ of $\bta(\theta)$ are continuously differentiable on $(1/2, B)$ with derivatives
\beno
\nabla_\theta\, \eta_t(\theta)
\hide{
\= \nabla_\theta \left[\theta - \theta\, \beta(\theta)^t\right]
}
\= 1 - \beta(\theta)^t - \dfrac{t}{\theta}\, \beta(\theta)^{t-1},
&& \theta \in (1/2, B).
\ee
We check conditions \czero---\cthree{ }one by one.

\s

\underline{Condition \czero.}
To show that the map $\bta: \Theta \mapsto \etaspace$ is one-to-one on $\Theta$,
we show that at least one coordinate of $\bta(\theta+\delta)$ must deviate from $\bta(\theta)$ for all $\theta \in (1/2, B)$ and all $\delta > 0$.
To do so,
note that $\bta(\theta)$ has at least two coordinates,
denoted by $\eta_1(\theta)$ and $\eta_2(\theta)$,
because $T \geq 2$ by assumption.
The first coordinate $\eta_1(\theta)$ of $\bta(\theta)$ is constant on $(1/2, B)$:
\beno
\eta_1(\theta)
\= 1,
&& \theta \in (1/2, B).
\ee
The second coordinate $\eta_2(\theta)$ of $\bta(\theta)$ is continuously differentiable on $(1/2, B)$ with derivative
\beno
\nabla_\theta\, \eta_2(\theta)
\= 1 - \beta(\theta)^2 - \dfrac{2}{\theta}\, \beta(\theta)
\= \dfrac{1}{\theta^2}
&>& 0,
&& \theta \in (1/2, B).
\ee
By the mean-value theorem,
\beno
\eta_2(\theta + \delta) - \eta_2(\theta)
\gte \dfrac{\delta}{(\theta + \delta)^2}
&>& 0,
&& \theta \in (1/2, B),
&& \delta > 0.
\ee
Thus,
$\eta_2(\theta)$ is strictly increasing on $(1/2, B)$ and at least one coordinate of $\bta(\theta+\delta)$ must deviate from $\bta(\theta)$ for all $\theta \in (1/2, B)$ and all $\delta > 0$. 
As a result,
the map $\bta: \Theta \mapsto \etaspace$ is one-to-one and continuous on $\Theta$.
Thus condition \czero{ }is satisfied.

\s

\underline{Condition \cone.}
Condition \cone{ }follows from the continuity of $\bta: \Theta \mapsto \etaspace$ and the upper semicontinuity of exponential-family loglikelihood functions \citep[][Lemma 5.3, p.\ 146]{Br86}.

\s

\underline{Condition \ctwo.}
Choose any $\theta \in \Theta$ and $\theta^\prime \in \Theta$ and let $\bmu(\bz) = \bA(\bz)^\top \bmu$ ($\bmu \in \mMs$).
By the triangle inequality,
we obtain,
for all $\theta \in \Theta$ and $\theta^\prime \in \Theta$ and all $\bmu \in \mMs$,
\be
\label{ccurved1}
\left|\langle\bta(\theta^\prime) - \bta(\theta),\, \bmu(\bz)\rangle\right|
\hide{
\= \left|\dsum_{k=1}^K \langle\bta_k(\theta^\prime) - \bta_k(\theta),\, \bmu_k(\bz)\rangle\right|
}
\= \left|\dsum_{t=1}^T \left[\eta_t(\theta^\prime) - \eta_t(\theta)\right] \mu_t(\bz)\right|\s
\\
\hide{
\= \left|\dsum_{t=1}^T \left[\eta_t(\theta^\prime) - \eta_t(\theta)\right] \dsum_{k=1}^K \mu_{k,t}(\bz)\right|
}
\lte \dsum_{t=1}^T |\eta_t(\theta^\prime) - \eta_t(\theta)| \left|\mu_t(\bz)\right|.
\ee
It can be shown that there exists $C > 2$ such that,
for all $\theta \in \Theta$ and all $t \in \{1, 2, \dots\}$,
\beno
\left|\nabla_\theta\, \eta_t(\theta)\right|
\lte \max(3, C),
&& t \in \{1, 2, \dots\},
\ee
which,
by the mean-value theorem,
implies that
\beno
|\eta_t(\theta^\prime) - \eta_t(\theta)|
\hide{
\= |\theta^\prime - \theta|\, |\nabla_\theta\, \eta_t(\theta)|_{\theta=\dot\theta}|
}
\lte |\theta^\prime - \theta|\, \max(3, C),
&& t \in \{1, 2, \dots\}.
\ee
Using \eqref{ccurved1} along with condition [C.$3^{\star\star}$] shows that there exist $C_1 > 0$ and $n_1 \geq 1$ such that,
for all $n > n_1$,
\beno
\label{ccurved0}
&& \left|\langle\bta(\theta^\prime) - \bta(\theta),\, \bmu(\bz)\rangle\right|
\;\leq\; \dsum_{t=1}^T |\eta_t(\theta^\prime) - \eta_t(\theta)| \left|\mu_t(\bz)\right|\s
\\
\lte |\theta^\prime - \theta|\, \max(3, C) \dsum_{t=1}^T \left|\mu_t(\bz)\right|
\;\leq\; C_1\, \norm{\theta^\prime - \theta}_2\, \dis{n \choose 2}.
\ee
Hence condition \ctwo{ }is satisfied,
because $\uuu = C\, {n \choose 2}$ ($C > 0$) in dense random graphs and because we assume that there is a single block.

\s

\underline{Condition \cthree.}
Using $|\beta(\theta)| < 1$ for all $\theta \in \Theta$,
\beno
\label{bounded.eta}
|\eta_t(\theta)|
\hide{
\= |\theta - \theta\, \beta(\theta)^t|
}
\lte |\theta| + |\theta|\, |\beta(\theta)|^t
\lte 2\, B,
&& \theta\in\Theta.
\ee
By condition [C.$4^{\star\star}$],
there exist $C_2 > 0$ and $n_2 \geq 1$ such that,
for all $n > n_2$,
\beno
\label{ccurved2}
\left|\langle\bta(\theta),\, s(\bx_1, \bz) - s(\bx_2, \bz)\rangle\right|
\;=\; \left|\dsum_{t=1}^T \eta_t(\theta)\, \left[s_t(\bx_{1}, \bz) - s_t(\bx_{2}, \bz)\right]\right|\s
\\
\leq\; 2\, B\, \left|\dsum_{t=1}^T  s_t(\bx_1, \bz) - \dsum_{t=1}^T  s_t(\bx_2, \bz)\right|
\;\leq\;
C_2\, d(\bx_1, \bx_2)\, L(\bz).
\ee
Thus condition \cthree{ }is satisfied.

\end{appendix}

\end{document}